\documentclass[11pt]{article}
\title{The largest eigenvalue of \\rank one deformation of large Wigner matrices}
\author{Delphine F\'eral \footnote{Institut de Math\'ematiques, Laboratoire de
Statistique et Probabilit\'es, Universit\'e Paul Sabatier, 31062 Toulouse Cedex 9, France; E-mail:
dferal@math.ups-tlse.fr} \, and Sandrine P\'ech\'e \footnote{Institut Fourier BP 74, 100 Rue des
    maths, 38402 Saint Martin d'Heres, France; E-mail: Sandrine.Peche@ujf-grenoble.fr}}
\usepackage{epsfig}
\usepackage{amsfonts}
\usepackage{amssymb}
\usepackage{amsthm}
\usepackage{amstext}
\usepackage{amscd}
\usepackage{amsmath}
\usepackage{delarray}
\begin{document}
\maketitle
\newtheorem{theo}{Theorem}[section]
\newtheorem{prop}[theo]{Proposition}
\newtheorem{lemme}[theo]{Lemma}
\newtheorem{conjecture}[theo]{Conjecture}
\newtheorem{definition}[theo]{Definition}
\newtheorem{fact}[theo]{Fact}
\newtheorem{hyp}[theo]{Assumption}
\theoremstyle{remark}
\newtheorem{rem}[theo]{Remark}
\newtheorem{remark}[theo]{Remark}
\newtheorem{Remark}[theo]{Remark}
\newtheorem{Notationnal remark}[theo]{Remark}
\newcommand{\bremnot}{\begin{Notationnal remark}}
\newcommand{\eremnot}{\end{Notationnal remark}}
\newcommand{\brem}{\begin{remark}}
\newcommand{\erem}{\end{remark}}
\newcommand{\bconj}{\begin{conjecture}}
\newcommand{\econj}{\end{conjecture}}
\newcommand{\bdefi}{\begin{definition}}
\newcommand{\edefi}{\end{definition}}
\newcommand{\bt}{\begin{theo}}
\newcommand{\bfa}{\begin{fact}}
\newcommand{\efa}{\end{fact}}
\newcommand{\Si}{\Sigma}
\newcommand{\et}{\end{theo}}
\newcommand{\bp}{\begin{prop}}
\newcommand{\ep}{\end{prop}}
\newcommand{\bl}{\begin{lemme}}
\newcommand{\el}{\end{lemme}}
\newcommand{\be}{\begin{equation}}
\newcommand{\ee}{\end{equation}}
\newcolumntype{L}{>{$}l<{$}}
\newenvironment{Cases}{\begin{array}\{{lL.}}{\end{array}}

\begin{abstract}
The purpose of this paper is to establish universality of the fluctuations of the largest
eigenvalue of some non necessarily Gaussian complex Deformed Wigner Ensembles. The real model is
also considered. Our approach is close to the one used by A. Soshnikov (c.f. \cite{Sos}) in the
investigations of classical real or complex Wigner Ensembles. It is based on the computation of
moments of traces of high powers of the random matrices under consideration.
\end{abstract}

\section{Introduction: model and results}
The scope of this paper is to study the spectral properties of some well chosen rank one
perturbation of classical complex or real large Wigner matrices. Our model can matricially be
described by a sequence $(M_N)_N$ of some complex or real Deformed Wigner matrices given by
$$M_N=\frac{1}{\sqrt{N}} W_N +A_N$$ where $A_N=(A_{i,j})_{1 \leq i,j \leq N}$ is the $N \times N$
deterministic real
matrix defined by $A_{i,j}=\frac{\theta}{N}$, with $\theta>0$ given independent of $N$ and, in the complex case, $W_N=(W_{i,j})_{1 \leq i,j \leq N}$ is a $N \times N$ Wigner Hermitian matrix with non necessarily Gaussian entries such that \\
(i) on the diagonal, the entries are real and the $\{ W_{i,i}, \: 1 \leq i \leq N \} \cup  \{ \Re
eW_{i,j} , \, \Im mW_{i,j}: \:  1 \leq i < j \leq N\}$ are real independent random
variables,\\
(ii) all these real variables have symmetric laws (as a consequence,
$\mathbb E[W_{i,j}^{2k+1}] =0$ for all $k \in \mathbb N ^*$),\\
(iii) $ \forall \, i <j, \, \mathbb E[(\Re e W_{i,j})^{2}]=\mathbb E[(\Im mW_{i,j})^{2}]=
\frac{\sigma^2}{2}$. The second moments of the diagonal elements
$W_{i,i}$ are assumed to be uniformly bounded,\\
(iv) all their other moments are assumed to be sub-Gaussian i.e. there exists a constant $\beta>0$
such that uniformly in $i,j$ and $k$,
$$\mathbb E[(\Re e W_{i,j})^{2k}] , \: \mathbb E[(\Im mW_{i,j})^{2k}]\leq {(\beta \, k )}^k.$$
In the real setting, $W_N=(W_{i,j})_{ 1 \leq i , j \leq N}$ is a $N\times N$ (non necessarily
Gaussian) Wigner symmetric matrix
which satisfies the following conditions\\
(i') the $\{ W_{i,j}, \: 1 \leq i \leq j \leq N\}$ are independent random variables,\\
(ii') the laws of the $W_{i,j}$ are symmetric (in particular,
$\mathbb E[W_{i,j}^{2k+1}] =0$),\\
(iii') for all $i <j$, $\mathbb E [W_{i,j} ^2]=\sigma^2$. The second moments of the $W_{i,i}$ are
assumed to be
uniformly bounded,\\
(iv') all the other moments of the $W_{i,j}$ grow not faster than the Gaussian ones. This means
that there is a constant $\beta>0$ such that, uniformly in $i, \, j$
and $k$, $\, \mathbb E[W_{i,j}^{2k}] \leq {(\beta \, k )}^k$.\\
When the entries of $W_N$ are further assumed to be Gaussian (with, on the diagonal, $W_{i,i} \sim
\mathcal N (0, \sigma ^2)$) that means in the complex (resp. real) setting when $W_N$ is element
of the so-called GUE (resp. GOE), we will denote by $M_N^G$ the corresponding Deformed
model.\\

Let $\lambda_1 \geq \cdots \geq \lambda_N$ be the ordered eigenvalues of $M_N$. At this point,
$M_N$ may be real or complex. If $\theta=0$, one recovers the classical Wigner Ensembles whose
spectrum is quite well-known. Our interest is to study the influence of the parameter $\theta$ on
this spectrum and mainly on the largest eigenvalues. Some answerings have yet been
obtained. \\
First, on a global setting, the classical Wigner Theorem is still satisfied whatever the parameter
$\theta \geq 0$ is (this is for example a consequence of Lemma 2.2 of \cite{Bai}). Thus, the
limiting behaviour of the empirical spectral measure $\mu_N=\frac{1}{N} \sum_{i=1}^N
\delta_{\lambda_i}$ of any ensemble of type $(i)-(iv)$ (or $(i')-(iv')$) is the semicircle law
$\mu_{\sigma}$ whose density is given by \be \frac{d \mu_{\sigma}}{dx}(x)=  \frac{1}{2 \pi
\sigma^2} \sqrt{4\sigma^2 - x^2} \, 1_{[-2\sigma,2 \sigma]}(x). \ee On the other hand, the
parameter $\theta$ may affect the limiting behavior of the largest eigenvalues. Let us recall the
results obtained for classical Wigner Ensembles. We denote by $\tilde \lambda _1 \geq \tilde
\lambda _2 \geq \cdots \geq \tilde \lambda _N$ the eigenvalues of $(\frac{1}{\sqrt{N}} W_N)_N$. It
is a fundamental result due to \cite{Geman} that the largest eigenvalue $\tilde \lambda _1$
converges almost surely to the right endpoint $2\sigma$ of the semicircle support. Then it was
established in \cite{TracyWidom} and for the GUE (resp. GOE) that, for all real $t$,
$$\displaystyle{\lim_{N \rightarrow \infty}\mathbb P \Big [ N^{2/3} ( \tilde
  \lambda _1 - 2 \sigma) \leq t \Big ]}=F_2^{TW}(t) \: (\text{resp. }F_1^{TW}(t)),\text{  where
$F_2^{TW}(t)$ (resp. $F_1^{TW}(t)$) }$$is  the well-known GUE (resp. GOE) Tracy-Widom distribution
(see \cite{TracyWidom} for precise definitions). A. Soshnikov later extended in \cite{Sos} these
results to arbitrary complex (resp. real) non-Gaussian Wigner matrices $(\frac{1}{\sqrt{N}}W_N)_N$
of type $(i)-(iv)$
(resp. $(i')-(iv')$). \\
Recently, the behavior of the largest eigenvalues of complex Deformed GUE was investigated in
details in $\cite{Pe}$ (see also $\cite{Fu-Ko}$). It is proved therein that the fluctuations of the
largest eigenvalue $\lambda _1 ^G$ of $(M_N^G)_N$ exhibit a phase transition according to the
value of $\theta$. Define \begin{equation}\label{rhoet sigma}\rho_{\theta}=\theta +
\frac{\sigma^2}{\theta} \text{ and } \sigma_{\theta}=\sigma \sqrt{\frac{\theta^2-\sigma^2}{\theta^2}}.
\end{equation}

\bt \label{ThPe} \cite{Pe} For any real $t$,
\begin{itemize}
\item if $\theta>\sigma,$ then
$\displaystyle{\lim_{N \rightarrow \infty} \mathbb P \left [N^{1/2}\left (
  \lambda_{1}^G-\rho_{\theta}\right )\leq t\right ]=\frac{1}{\sqrt {2\pi
}\sigma_{\theta}}\int_{-\infty }^t e^{\{-\frac{u^2}{2\sigma_{\theta}^2}\}} \, du,}$
\item if $\theta<\sigma,$ then
 $\displaystyle{\lim_{N \rightarrow \infty} \mathbb P \left[N^{2/3}\left ( \lambda_{1}^G-2\sigma\right ) \leq t\right ]= F^{TW}_2(t),}$
 \item if $\theta=\sigma,$ then
$\displaystyle{\lim_{N \rightarrow \infty} \mathbb P \left [N^{2/3}\left (
      \lambda_{1}^G-2\sigma\right ) \leq t\right ]= F^{TW}_{3}(t),\text{ where}}$
$F^{TW}_{3}$ is some generalized Tracy-Widom distribution (see \cite{Pe}, p. 2-4 and \cite{BBP} Subsection 3.3 for precise definitions).
\end{itemize}
\et 
The extension to the real case has not been obtained yet. Nevertheless, it can be inferred from
the results of \cite{Debashis} and communications with J. Baik (forthcoming paper \cite{BS}). In
particular, denoting by $F^{TW}_1$ the GOE Tracy-Widom distribution, one should obtain the
following result.

\bconj \label{conjectPecastheta><sigma}
Let $\lambda_{1}^G$ be the largest eigenvalue of the Deformed GOE. For all real $t$, \\
(i) If $ \theta> \sigma$, then $\displaystyle{\lim_{N \rightarrow \infty} \mathbb P \left
[N^{1/2}\left (
  \lambda_{1}^G-\rho_{\theta}\right )\leq t \right ]=  \frac{1}{\sqrt
    {4\pi }\sigma_{\theta}}\int_{-\infty }^t
  e^{\{-\frac{u^2}{4\sigma_{\theta}^2}\}} \, du.}$ \\
(ii) If $ \theta< \sigma$, then $\displaystyle{\lim_{N \rightarrow \infty} \mathbb P \left
[N^{2/3}\left (
  \lambda_{1}^G-2 \sigma \right )\leq t \right ]= F^{TW}_1(t)}$.
\econj

Some generalizations of Theorem \ref{ThPe} have been obtained. In $\cite{Fe},$ the almost sure
limit of the first largest eigenvalues of any complex or real Deformed Wigner model $(M_N)_N$ is
investigated. It is proved therein that for any complex or real Deformed Wigner matrix $(M_N)_N$,
the largest eigenvalue $\lambda _1$ a.s. jumps outside the support $[-2 \sigma, 2 \sigma]$ of the
semicircle law to the value $\rho _{\theta}$ as soon as $\theta > \sigma$. If $0 \leq \theta \leq
\sigma$, $\lambda _1$ still tends to the right edge
$2 \sigma$.\\

Our paper is mainly devoted to the study of fluctuations of the largest eigenvalue of non
necessarily Gaussian complex Deformed Wigner Ensembles of type $(i)-(iv)$ and of parameter
$\theta$. Our main result is that the universality holds for any $\theta>0$. Our investigation
also concerns non necessarily Gaussian real Deformed Wigner Ensembles of type $(i')-(iv')$ and
yields the proof of the second point of
Conjecture \ref{conjectPecastheta><sigma}.\\

We first prove the following universality result. \bt \label{Uni1Pe} Theorem \ref{ThPe} is true
for the largest eigenvalue $\lambda _1$ of any complex Deformed Wigner Ensembles of type
$(i)-(iv).$ \et

When $\theta<\sigma$, we can state a stronger result namely that the parameter $\theta$ does not
affect the asymptotic behavior of the distribution of the $k$ first largest eigenvalues of any
complex Deformed Wigner Ensemble of type $(i)-(iv)$, for any fixed integer $k\geq 1$. Hence, all
the asymptotic results established in \cite{Sos}  for general non-Gaussian Hermitian Wigner
Ensembles ($\theta=0$) extend to the case where $\theta<\sigma$.

\bt \label{UnikPe} Assume that $\theta<\sigma$. Let $k\geq 1$. Let $\lambda_{i}$ denote the
$i^{th}$ largest eigenvalue of a complex Deformed Wigner Ensemble of type $(i)-(iv)$. Then, for
all $(t_1, \cdots , t_k) \in \mathbb R ^k$,
$${\lim_{N \rightarrow \infty} \mathbb P \left[N^{2/3}\left (
    \lambda_{1}-2\sigma\right ) \leq t_1 , \cdots , N^{2/3}\left (
    \lambda_{k}-2\sigma\right ) \leq t_k\right ]= F^{TW}_{2,k}(t_1,
    \cdots, t_k),\text{ where}}$$
 $F^{TW}_{2,1}=F^{TW}_{2}$ and $F^{TW}_{2,k}$ is the Tracy-Widom limiting joint distribution of the $k$ first
    eigenvalues of the GUE (given e.g. in \cite{Sos}).
\et

All our middle results being also true in the real setting (with small modifications), we also
consider throughout this paper the real model. Actually, once the whole real version of Theorem
$\ref{ThPe}$ will be proven, our main Theorem $\ref{Uni1Pe}$ can readily be extended to the real
framework (see the next Section \ref{core} for a justification). Yet we prove an analog of Theorem
\ref{UnikPe} in the real framework which in particular gives the last point of Theorem \ref{ThPe}
in the real case.

\bt \label{theo: Unireal}Assume that $\theta< \sigma$.  Let $k\geq 1$. Let $\lambda_{i}$ denote the
$i^{th}$ largest eigenvalue of any real Deformed Wigner Ensemble of type $(i')-(iv')$. Then, for
all $(t_1, \cdots , t_k) \in \mathbb R ^k$,
$${\lim_{N \rightarrow \infty} \mathbb P \left[N^{2/3}\left (
    \lambda_{1}-2\sigma\right ) \leq t_1 , \cdots , N^{2/3}\left (
    \lambda_{k}-2\sigma\right ) \leq t_k\right ]= F^{TW}_{1,k}(t_1,
    \cdots, t_k),\text{ where}}$$
 $F^{TW}_{1,1}=F^{TW}_{1}$ and $F^{TW}_{1,k}$ is the Tracy-Widom limiting joint distribution of the $k$ first
    eigenvalues of the GOE (given e.g. in \cite{Sos}).
\et

\brem At this point we would like to point out the fact that results of Theorem \ref{ThPe} have
been proved for more complex Deformed GUE models. On the one hand, because of the rotational
invariance of the GUE distribution, Theorem \ref{ThPe} holds for arbitrary deterministic matrix
$A_N$ of rank one and of eigenvalue $\theta.$ On the other hand, the results of \cite{Pe} are
stated for any deterministic deformations $A_N$ of fixed rank $k\geq 1$. The natural problem of
the universality of the fluctuations arises for such deformations but is beyond the scope of this
paper. In a forthcoming paper, we will prove that universality does not hold for instance if one
chooses the diagonal matrix $A_N=\text{diag}(\theta, 0\ldots , 0).$ We will also investigate
deformations of fixed rank $k\geq 1$. \erem

The derivation of our results uses ideas and combinatorial techniques similar to those used by Y.
Sinai and A. Soshnikov in \cite{Si-So1}, \cite{Si-So2} and \cite{Sos}. Following especially the
approach developed in \cite{Sos}, we compute the limiting behavior of the expectation of traces of
high moments of $M_N$ defined by \be \label{traceg} \mathbb{E}[tr \, M_N^{2s(+1)}]=\sum_{1 \leq
i_0,i_1, \cdots , i_{2s-1(+1)}\leq N}
  \mathbb{E}[M_{i_0,i_1} \, \cdots \, M_{i_{2s-1(+1)},i_0}].
\ee for some powers $s=s_N$ such that $\lim_{N \to \infty}s_N=\infty$. In particular, we study in
details the contribution to (\ref{traceg}) from the closed paths $\mathcal{P}=\{i_0,i_1\, \cdots
\,i_{2s-1(+1)},i_0\}$ of length $2s(+1)$ on the set of vertices $\{1, \cdots, N\}$. The strategy
is to show that the leading term in the asymptotic expansion of (\ref{traceg}), for specific
exponent $s_N$, comes from the paths whose expectation only depends on $\theta$ and $\sigma$. This
implies that, up to a negligible error, (\ref{traceg}) has the same limiting behavior as in the
case the matrices $W_N$ are of the GUE. In other words we show that $\displaystyle{\mathbb{E}[tr
\, M_N^{2s(+1)}]=\mathbb{E}[tr \,\left ( M_N^G\right
  )^{2s(+1)}](1+o(1))}$. This strategy can be deepened to derive similar results for all
higher moments (see Section 6).\\
Our paper is organized as follows. In Section 2, we present the main elements of the proof of our
theorems and explain why they follow from universal limiting behavior of moments of high traces of
$M_N$. We then recall in Section 3 the needed specific terminology introduced in \cite{Sos}.
Sections 4, 5 and 6 are devoted to the complete proof of case $\theta >\sigma$. We next consider
the case where $\theta=\sigma$ in Section 7. At last, we justify the case where $0 < \theta <
\sigma$ in Remark $\ref{casetheta<sigma}$.\\
From a notationnal point of view, throughout this paper, the notations $C, C_i,1\leq i\leq 6, C'$ and $\tilde C$ will be used for different
positive constants.
\section{Core of the proof} \label{core}
Here, we first mainly concentrate on the case where $\theta>\sigma$. We show how universality of
the fluctuations of the largest eigenvalue of general complex Deformed Wigner Ensembles of type
$(i)-(iv)$ can be derived from the computation of limiting moments of traces of high powers of
$M_N$. This is inspired from the approach of $\cite{Sos}$. At the end of this section, we point
out the main modifications needed in both the cases where $\theta=\sigma$ and $\theta<\sigma$ and
also quickly discuss on
the real setting.\\

In the case where $\theta>\sigma$, we shall handle with powers $s_N$ of the order $\sqrt N$. It is
indeed expected, from Theorem \ref{ThPe}, that the largest eigenvalue $\lambda_1$ exhibits
Gaussian fluctuations around $\rho_{\theta}$ in the scale $\frac{1}{\sqrt N}.$ In particular, we
prove (in Section 5) the crucial fact that, for $s_N=[t\sqrt N]$ with $t>0,$
\begin{eqnarray} \label{limitTrace}
\lim_{N\rightarrow \infty} \left (\mathbb E\Big [Tr \left (\frac{M_N}{\rho_{\theta}}\right
)^{2s_N} + Tr \left (\frac{M_N}{\rho_{\theta}}\right )^{2s_N+1} \Big ] -\mathbb E\Big [Tr \left
(\frac{M_N^G}{\rho_{\theta}}\right )^{2s_N} + Tr \left (\frac{M_N^G}{\rho_{\theta}}\right
)^{2s_N+1} \Big ]\right)=0.
\end{eqnarray}
Basically, one intends to prove that only the largest eigenvalue $\lambda_1$ (resp. $\lambda_1^G$)
contributes to the first (resp. second) expectation in the l.h.s of $(\ref{limitTrace})$. Let
$\lambda_1^G \geq \cdots \geq \lambda_N^G $ be the ordered eigenvalues of $(M_N ^G)_N$ (of the
Deformed GUE). We decompose the eigenvalues of $M_N$ and $M_N^G$ as follows
\begin{eqnarray*}
&&\lambda_j=\rho_{\theta}\left (1+\frac{\xi_j}{2\sqrt N}\right ), \text{ if }\lambda_j>0 \: \text{
and }\lambda_j=-2\sigma+\frac{\tau_j}{N^{2/3}}, \text{ if }\lambda_j<0,\crcr
&&\lambda_j^G=\rho_{\theta}\left(1+\frac{\xi_j^G}{2\sqrt N}\right ), \text{ if }\lambda_j^G>0 \:
\text{ and }\lambda_j^G=-2\sigma+\frac{\tau_j^G}{N^{2/3}}, \text{ if }\lambda_j^G<0.
\end{eqnarray*}

The strategy to derive universality for the fluctuations of the largest eigenvalue from (\ref{limitTrace}) can be summarized in three steps.\\
The first step shows that for both $M_N$ and $M_N^G$ and for all $t>0$, the random variable \be
\label{eqbase} {\epsilon}_{N,t}^{(G)}=\frac{1}{2}\Big [Tr \left
(\frac{M_N^{(G)}}{\rho_{\theta}}\right )^{2[t\sqrt N]} + Tr \left
(\frac{M_N^{(G)}}{\rho_{\theta}}\right )^{2[t\sqrt N]+1} \Big ]- \sum_{|\xi_j^{(G)}|\leq
N^{1/6}}e^{t\xi_j^{(G)}} \ee converges to $0$ a.e. as in the $L^k$ norm, for all fixed $k\geq 1.$
Formula (\ref{eqbase}) will be proved below. In the case where $\theta>\sigma$ (only), the a.e.
and $L^1$ norm convergence are actually enough to derive the announced
fluctuations for the largest eigenvalue (see also the end of the section for the other cases). \\

The second step follows from results of \cite{Pe} which ensure that $\exists \: \delta>0$ such that for $|t|\leq \delta,$ \be \label{tdeL} {\lim_{N
\rightarrow \infty} \mathbb E \Big [
  \sum_{|\xi_j^G|\leq N^{1/6}}e^{t\xi_j^G} \Big ]=L_{\theta}(t),
  \quad }
\text{where $L_{{\theta}}$ is the Laplace transform of the law $\mathcal N (0, {\sigma
_{\theta}}^2)$}. \ee Let then $\mathbb P_N(\lambda_1, \ldots, \lambda_N)$ be the symmetrized joint
eigenvalue distribution on $\mathbb R ^N$ induced  by any Deformed Ensemble $M_N$ and, for all $1
\leq m \leq N$, denote by $\mathbb P_N^m$ one of its $m$-dimensional marginal. Define then the
associated $m$-point correlation function $R_m$ of $\mathbb P_N$ by \be R_m=\frac{N!}{(N-m)!}
\mathbb P _N ^m .\label{defrm}\ee Note that $R_m$ is a distribution in general. From (\ref{tdeL}),
the above results combined with the machinery developed in \cite{Sos} yield that the rescaled one
point correlation function $\rho_1(x)=\frac{1}{\sqrt
  N}R_1(\rho_{\theta}+\frac{x}{\sqrt N})$ of $M_N$ satisfies, for all $t \in \mathbb R$,
$\lim_{N \rightarrow
  \infty}\int_{t}^{\infty}\rho_1(x) \, 1 _{x\leq N^{1/6}} \,
dx=\int_{t}^{\infty}\frac{1}{\sqrt{2\pi
    \sigma^2_{\theta}}}\exp{\{-\frac{u^2}{2\sigma^2_{\theta}}\}} \,
du.$ Actually, the following stronger result holds \be \label{restotal}\lim_{N \rightarrow
  \infty}\int_{t}^{\infty}\rho_1(x) \,
dx=\int_{t}^{\infty}\frac{1}{\sqrt{2\pi
    \sigma^2_{\theta}}}\exp{\{-\frac{u^2}{2\sigma^2_{\theta}}\}} \,
du. \ee Indeed, in Lemma \ref{Prop: fondbis} given below, we will establish that there exist two
constants $C_3,C_4>0$ such that $\displaystyle{\mathbb P \Big [\sharp \{i: \, \lambda_i >
\rho_{\theta}(1 + \frac{1}{N^{1/3}})\} >0 \Big ] \leq C_3  \exp {\{-C_4 N^{1/6}\}.}}$ This yields
$(\ref{restotal})$ since
\begin{eqnarray}{\label{invrestronque}}
\int_{t}^{\infty}\rho_1(x) \, 1_{x> N^{1/6}} \, dx & \sim &
\int_{]\rho_{\theta}(1+\frac{1}{N^{1/3}}), \infty[}R_1(y) \, dy \leq NC_3
  \exp {\{-C_4 N^{1/6}\}} \to 0  \, \text{ as $N \to \infty$}.
\end{eqnarray}

At last, the third step is based on Lemma \ref{lemm: nombrevpenro} proven below which states that
only the largest eigenvalue of $M_N$ separates from the bulk and that it is close to
$\rho_{\theta}.$ To be more precise, this lemma implies that there exists $C_1, C_2 >0$ such that,
given $R<0$,  $\forall t\geq R$, \be\label{majod_n2}\mathbb P \left [\exists \text{ at least two
eigenvalues of }M_N \text{ in } ]\rho_{\theta}+\frac{t}{\sqrt N}, \infty[\right]\leq
C_1\exp{\{-C_2s_N\}},\ee for $N$ large enough. In this way, we claim that the largest eigenvalue
$\lambda _1$ of any complex Deformed Ensembles of type $(i)-(iv)$ has the same limiting behavior
as that of the Deformed GUE. To see it, let $D_N$ denote the random number of eigenvalues of $M_N$
in the interval $I=]\rho_{\theta}+\frac{t}{\sqrt N}, \infty[.$ Then, with $R_m(x_1,\ldots, x_m)$
given as in (\ref{defrm}), one has
\begin{eqnarray} {\label{P(DN>1)}}
\mathbb P[ \sqrt{N}( \lambda _1 - \rho _{\theta}) \leq t] & = & \mathbb P (D_N=0) \crcr &=
&1-\mathbb{E}\left (\sum_{i=1}^N 1_I(\lambda_i)\right)+\sum_{m\geq
  2}\frac{(-1)^m}{m!}\int_{I^m}R_m(x_1,\ldots, x_m) \,\prod_{i=1}^m dx_i.
\end{eqnarray}

Then, using (\ref{majod_n2}),
\begin{eqnarray*}
&& |\sum_{m\geq 2}\frac{(-1)^m}{m!}\int_{I^m}R_m(x_1,\ldots, x_m)dx_1\cdots dx_m \Big|=\Big|
\mathbb P (D_N=0)-1+\mathbb{E}\left (\sum_{i=1}^N 1_I(\lambda_i)\right)\Big |\crcr &&=\Big
|\mathbb{E}\left (\sum_{i=1}^N 1_I(\lambda_i)\right)-\mathbb P\left
  (D_N \geq 1 \right)\Big |
=\Big | \mathbb{E}\left (\sum_{i=1}^N 1_I(\lambda_i)1_{D_N\geq
    2}\right)-\mathbb P\left (D_N\geq 2\right)\Big |\crcr
&&\leq (N+1) \, \mathbb P (D_N\geq 2) \leq(N+1) C_1\exp{\{-C_2s_N\}}\rightarrow 0\text{ as }N \to
\infty.
\end{eqnarray*}

Noticing that $\displaystyle{\mathbb{E}\left (\sum_{i=1}^N
  1_I(\lambda_i)\right)=\int_{I}R_{1}(y)dy=\int_{t}^{\infty}\rho_1(x)  \,
dx}$, we derive Theorem \ref{Uni1Pe} for $\theta > \sigma$.

\paragraph{}Let us now return to formula (\ref{eqbase}) and prove the
announced convergence. We first show that the negative eigenvalues do not contribute to
$(\ref{eqbase})$. Given a real $c>0$, we define for $i=1,2$
\begin{eqnarray} \label{defri}
r_i={\sum}^{i} \Big
  (1/2\left(\frac{\lambda_j}{\rho_{\theta}}\right)^{2s_N}+1/2\left(\frac{\lambda_j}{\rho_{\theta}}\right)^{2s_N+1}\Big )
\end{eqnarray}
where $\sum ^{1}$ (resp. $\sum ^{2}$) corresponds to the summation over $\{j: \,
-2\sigma-c/N^{2/3}<\lambda_j<0 \}$ (resp.  $\{j: \,\lambda_j<-2\sigma-c/N^{2/3} \}$). As $\rho
_{\theta} > 2 \sigma$, it is an easy fact that $\displaystyle{|r_1|\leq \exp{\{-C s_N\}},}$ for
$N$ large enough (for a constant $C>0$ ). Considering $r_2$, we have that
$$\displaystyle{|r_2| \leq \frac{1+\rho_{\theta}/2\sigma}{2}
  \sum_{\lambda_j<-2\sigma-c/N^{2/3}} \left ( - \frac{\lambda
  _j}{\rho_{\theta}} \right )^{2 s_N +1} }.$$
But, denoting as above by $\tilde \lambda_i$ the $i^{th}$ largest eigenvalue of the corresponding
rescaled Wigner matrix $\frac{1}{\sqrt N}W_N$, the so-called interlacing property of eigenvalues
states that \be \label{interlac} \lambda_1\geq \tilde \lambda_1\geq \lambda_2\geq \tilde
\lambda_2\geq \cdots \geq \lambda_N\geq \tilde \lambda_N. \ee This implies that
$\displaystyle{|r_2| \leq
  \frac{1+\rho_{\theta}/2\sigma}{2} \times \frac{\rho_{\theta}}{2\sigma}
  \times Tr
  \left(\frac{W_N/\sqrt N}{\rho_{\theta}} \right) ^{2 s_N +2}.} $
By the investigations of $\cite{Si-So1}$ recalled in Theorem \ref{theo : soshni1} below, we then
deduce that all the moments of $r_2$
  vanish as $N$ goes to $\infty$.\\

Now, we examine the contribution of the positive eigenvalues which can be expressed as the sum
$r_3+r_4+r_5$ where for $i=3,4$ and $5$, $r_i$ is given by $(\ref{defri})$, the summation $\sum
^{3}$ is over $\{j: \, \lambda_j>0, \, |\xi_j|\leq N^{1/6} \}$,  $\sum ^{4}$ over $\{j: \,
\lambda_j>0, \, \xi_j\leq -N^{1/6}
  \}$ and $\sum ^{5}$ over $\{j: \, \lambda_j>0, \, \xi_j\geq N^{1/6}
  \}$. \\
First, it is an easy fact that
 $r_3=\sum_{|\xi_j|\leq N^{1/6}}e^{t\xi_j}
 (1+O(N^{-1/6}))$. We then show that the other terms lead to a negligible contribution. First, one readily has that $\displaystyle{|r_4|\leq N\exp{\{-C N^{1/6}\}}.}$
The analysis of the term $r_5$ leans on the following lemma.

\bl \label{Prop: fond} If $\theta>\sigma$, for all $k$ in $\mathbb N^* $, for any $t_i, 1\leq
i\leq k,$ in a compact subset $K$ of $\mathbb R ^{+*}$,
$$\exists\: C=C(K)>0, \:\mathbb{E}\left ( \prod_{i=1}^kTr \left
(\frac{M_N}{\rho_{\theta}}\right )^{2[t_i\sqrt N](+1)}\right )\leq C^k.$$ \el For the proof, we
refer to Sections 5 and 6. Thanks to this result, we can estimate the contribution of $r_5$. Note
that
\begin{eqnarray*}
&|r_5|&\leq \sum_{\xi_j\geq N^{1/6}}\left(\frac{\lambda_j}{\rho_{\theta}}\right)^{2s_N+1} \leq
\sum_{\xi_j\geq
N^{1/6}}\left(\frac{\lambda_j}{\rho_{\theta}}\right)^{2s_N+1}\left(\frac{\lambda_j}{\rho_{\theta}}\right)^{2s_N-1}\left
( \frac{1}{1+\frac{N^{1/6}}{2\sqrt N}}\right)^{2s_N-1}\crcr &&\leq Tr\left
(\frac{M_N}{\rho_{\theta}}\right)^{4s_N} \exp{\{-C N^{1/6}\}}.
\end{eqnarray*}
According to Lemma $\ref{Prop: fond}$, we trivially deduce that
all the moments of $r_5$ tend to zero as $N \to \infty$. This finishes the proof of formula $(\ref{eqbase})$, yielding the first step.\\

Lemma $\ref{Prop: fond}$ also ensures that the positive eigenvalues are not too large. This is
stated in the next lemma which completes the proof of $(\ref{invrestronque})$.

\bl \label{Prop: fondbis} There exist two positive constants $C_3$ and $C_4$ such that, for $N$
large enough,
$$\mathbb P \Big [\sharp \{i: \, \lambda_i > \rho_{\theta}(1 +
\frac{1}{N^{1/3}})\} >0 \Big ] \leq C_3  \exp (-C_4 N^{1/6}).$$ \el

\paragraph{Proof of Lemma \ref{Prop: fondbis}:}
From the Chebytchef inequality, we readily have that
$$\mathbb P \Big [\sharp \{i: \, \lambda_i > \rho_{\theta}(1 +
\frac{1}{N^{1/3}})\} >0 \Big ] \leq  \Big (\rho_{\theta} + \frac{\rho_{\theta}}{N^{1/3}} \Big )
^{- 2 [\delta N^{1/2}]}  \mathbb E\left ( Tr M_N ^{2 [\delta N^{1/2}]} \right )$$ where $\delta$
is a real $>0$. For $N$ large enough, we derive from Lemma \ref{Prop: fond} that $\exists\: C_3>0,
\tilde C_4>0$ such that $\displaystyle{ \mathbb P \Big [\sharp \{i: \, \lambda_i > \rho_{\theta}(1
+
\frac{1}{N^{1/3}})\} >0 \Big ] \leq  C_3 e^{-\tilde C_4 \delta N^{1/6}}. \, \square}$\\

At this stage, it remains to prove the following fundamental lemma which readily gives
$(\ref{majod_n2})$. \bl \label{lemm: nombrevpenro} There exists a positive constant $C_2$ such that
$$\mathbb P \Big [\sharp \{i \geq 2: \, \lambda_i > 2\sigma
  +(\rho_{\theta}-2\sigma)/2\} >0 \Big ]\leq \exp{\{-C_2 s_N\}}.$$
\el
\paragraph{Proof of Lemma \ref{lemm: nombrevpenro}:} By the interlacing
property of eigenvalues, it is clear that only one eigenvalue of $M_N$ is close to $\rho_{\theta}$
since
\begin{eqnarray*}
\mathbb P \Big [\sharp \{i \geq 2: \, \lambda_i > 2\sigma
  +(\rho_{\theta}-2\sigma)/2\} >0 \Big ] & \leq &
\mathbb P \Big [\sharp \{i \geq 2: \, \tilde \lambda_i > 2\sigma
  +(\rho_{\theta}-2\sigma)/2\} >0 \Big ] \cr
& \leq & \frac{\mathbb E[ Tr (W_N/\sqrt N) ^{2s_N}]}{\Big
  (2\sigma +(\rho_{\theta}-2\sigma)/2 \Big ) ^{2s_N}}
\leq  \exp{\{-C_2 s_N\}}
\end{eqnarray*}
for some positive constant $C_2.$
The last inequality follows from Theorem \ref{theo : soshni1}. $\square$\\

Thus we get the statement of Theorem \ref{Uni1Pe}  in the case $\theta>\sigma$. In the real
setting, one can expect the same proof with $\mathcal N (0, \sigma _{\theta} ^2)$ instead of the
law $\mathcal N (0, 2\sigma _{\theta} ^2)$ (recall Conjecture
\ref{conjectPecastheta><sigma}).\\

In both other cases where $\theta=\sigma$ and $\theta<\sigma$, the fluctuations of the largest
eigenvalue are expected to occur in the scale $N^{-2/3}$ around the edge $2\sigma$. This is
exactly as for classical Wigner Ensembles ($\theta=0)$ except the derivated limiting distribution.
The scheme to state the complete universality follows the same steps as in the case $\theta >
\sigma$ (with $2 \sigma$ instead of $\rho _{\theta}$ and replacing the law $\mathcal N (0, \sigma
_{\theta} ^2)$ with $F^{TW}_3$ if $\theta=\sigma$ and with $F^{TW}_2$ if $\theta<\sigma$). The asymptotics of correlation functions of the Deformed GUE required to establish the second step are straightforward from Propositions 2.1 and 2.2 in \cite{Pe} and Subsection 3.3 in \cite{BBP}.  The proof then
mainly boils down to universality of the limiting expectation of Traces of exponent of type
$o(N^{2/3})$ and $O(N^{2/3})$. Nevertheless the derivation of the result requires more complex
considerations than the previous analysis. Indeed, here, the largest eigenvalue does not separate
from the "bulk" and the whole spectrum lies in $[-2\sigma - \frac{1}{2 \sqrt{N}}, 2\sigma+
\frac{1}{2
    \sqrt{N}} ]$. In fact, the reasoning is very close
to that done by A. Soshnikov for general Wigner Ensembles and we refer to Sections 1,2 and 5 of
$\cite{Sos}$ for details. In particular, universality of all higher moments of the traces is
required. Note that in the case where $\theta<\sigma$, we actually prove (\ref{limitTrace}) and
(\ref{eqbase}) but also that the same formulaes hold with $M_N^G$ replaced with $\frac{1}{\sqrt N}
W_N^G$. Moreover, convergence of (\ref{eqbase}) in the $L^k$ norm for any fixed $k\geq 1$ ensures
universality of the limiting joint distribution of the $k$ first largest eigenvalues of any
Deformed Wigner Ensemble and that this limit is the one of the GUE. A detailed proof of this fact
is presented in \cite{Sos}. In the real setting (and again under $\theta<\sigma$), the same
reasoning shows that the fluctuations of eigenvalues of $M_N$ are compared to those, well known,
of the largest
eigenvalues of the GOE instead of $M_N^G$. \\

The rest of our paper is hence mainly devoted to the analysis of $(\ref{traceg})$ for some powers
$L_N=2 s_N(+1)$. This is based on the combinatorial machinery developed in $\cite{Si-So1}$,
$\cite{Si-So2}$ and $\cite{Sos}$. Before we proceed, we recall the main definitions needed in this
paper and introduced in \cite{Si-So1}-\cite{Sos}.

\section{Terminology: classification of instants and vertices}
To each term in the expectation (\ref{traceg}), we associate a path $\mathcal{P}= \{i_0,i_1, \cdots
, i_{L-1},i_{L}=i_o \}$  of length $L \geq 1$ where $i_j \in \mathbb N ^*$ (in this paper, we
restrict to vertices in $\{1, \cdots, N \}$). Note that loops are allowed i.e. it may happen that
$i_{j+1}=i_j$. To explain our counting strategy, we need to recall some definitions given in
\cite{Sos}.

\bdefi The instant $j=1, \ldots, L$ is said to be marked for the closed path $\mathcal{P}$ if the
unordered edge $(i_{j-1},i_j)$ occurs an odd number of times up to time $j$ (included). The other
instants are said to be unmarked. \edefi

Throughout this paper, we denote by $\mathcal{P}_{m,l}$ the set of paths $\mathcal{P}$ of length
$L=l+2m$ having $l+m$ marked instants and $m$ unmarked instants. In particular,
$\mathcal{P}_{m,0}$ corresponds to the classical even closed paths used in the framework of the
classical Wigner Ensembles. We associate to each path $\mathcal{P}$ a trajectory $x=\{x(t), \,
0\leq t \leq L\}$ of a simple random walk on the positive half-lattice such that
\begin{eqnarray*}
&&x(0)=0, \: x(L)=l; \: \forall t\in [0, L], \: x(t)\geq 0,\cr &&x(t)-x(t-1)=1 \: (\text{resp.
}-1) \quad \text{if $t\in \mathbb{N}^*$ is marked (resp. unmarked)} .
\end{eqnarray*}
Thus, the associated trajectory $x$ of a path $\mathcal{P}$ of $\mathcal{P}_{m,l}$ is such that
$l+m=\# \{t, \, x(t)-x(t-1)=1\}$ (up steps) and $m=\# \{t, \,
  x(t)-x(t-1)=-1\}$ (down steps). We define by $\mathcal{T}_{m,l}$ the set of such trajectories $x$
  and we let $T_{m,l}= \#\mathcal{T}_{m,l}$. The elements of $\mathcal{T}_{m,0}$
  are often called Dyck paths.

\bp \label{Prop: Tml}For $L=l+2m$, one has
$\displaystyle{T_{m,l}=C_{L}^{l+m}-C_{L}^{m-1}=\frac{L!}{(l+m+1)!m!}(l+1).}$ \ep

\brem Proposition \ref{Prop: Tml} is a straightforward consequence of the symmetry principle used
in the historical proof of the Wigner Theorem (c.f. $\cite{Bai}$ for example). One can also notice
that, amongst the paths of $\mathcal{T}_{m,l}$, exactly $T_{m,l-1}$ (resp. $T_{m-1,l+1}$) have a
last step up (resp. down). \erem

We also need to refine our classification of vertices of a path $\mathcal{P}$ of
$\mathcal{P}_{m,l}$.

\bdefi A marked instant $j$ is called an instant of self intersection of $\mathcal{P}$ if there
exists a marked instant $j'<j$ such that $i_{j'}=i_j.$ \edefi \bdefi A vertex $i$ is said to be a
vertex of simple (resp. $k$-fold) intersection of $\mathcal{P}$ if there exist exactly two (resp.
$k$) marked
instants such that $i_j=i.$\\
A path without self-intersection will be called a simple path. \edefi

We can now split the vertices of $\mathcal{P}\in \mathcal{P}_{m,l}$ into $l+m+1$ disjoint subsets
such that \be \{1,\ldots, N\}=\bigcup_{k=0}^{l+m}\mathcal{N}_k, \ee where $\mathcal{N}_k$ is the
subset of vertices of $k-$fold self intersection. Taking $N_k= \# \mathcal{N}_k$, such a path
$\mathcal{P}$ will be said of type $(N_o,N_1, \cdots, N_{l+m})$. In particular, the simple paths of
$\mathcal{P}_{m,l}$ are of type $(N-(l+m),l+m,0, \cdots, 0)$. It may happen that the origin $i_o$
is unmarked and then it is in $\mathcal{N}_o$. Otherwise $i_o$ is marked and $\mathcal{N}_o$
contains only vertices not belonging to ${\cal P}$.

\begin{center}
\begin{tabular}{cp{0.5cm}c}
\includegraphics[height=4cm,angle=0]{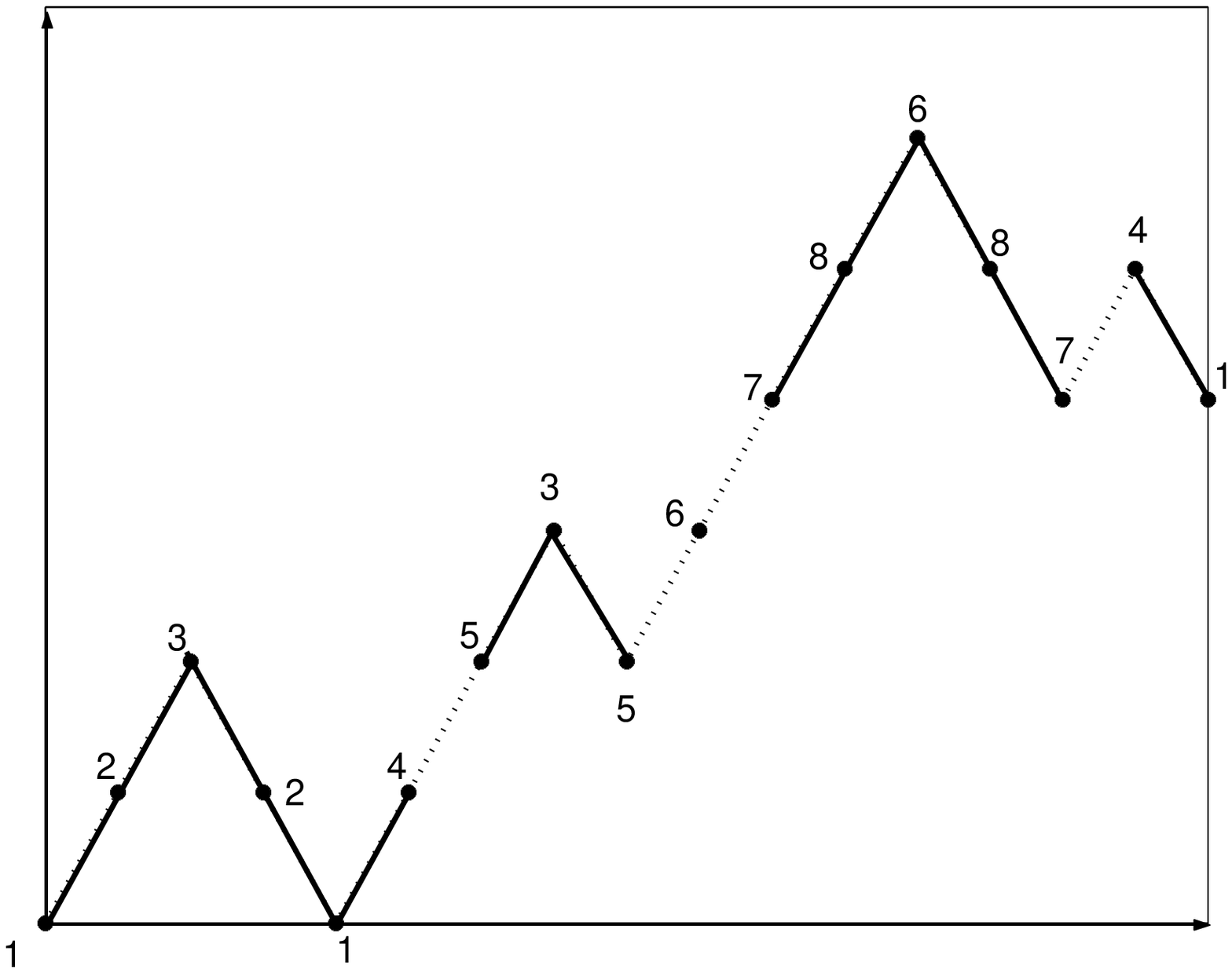} & &
\includegraphics[height=4cm,angle=0]{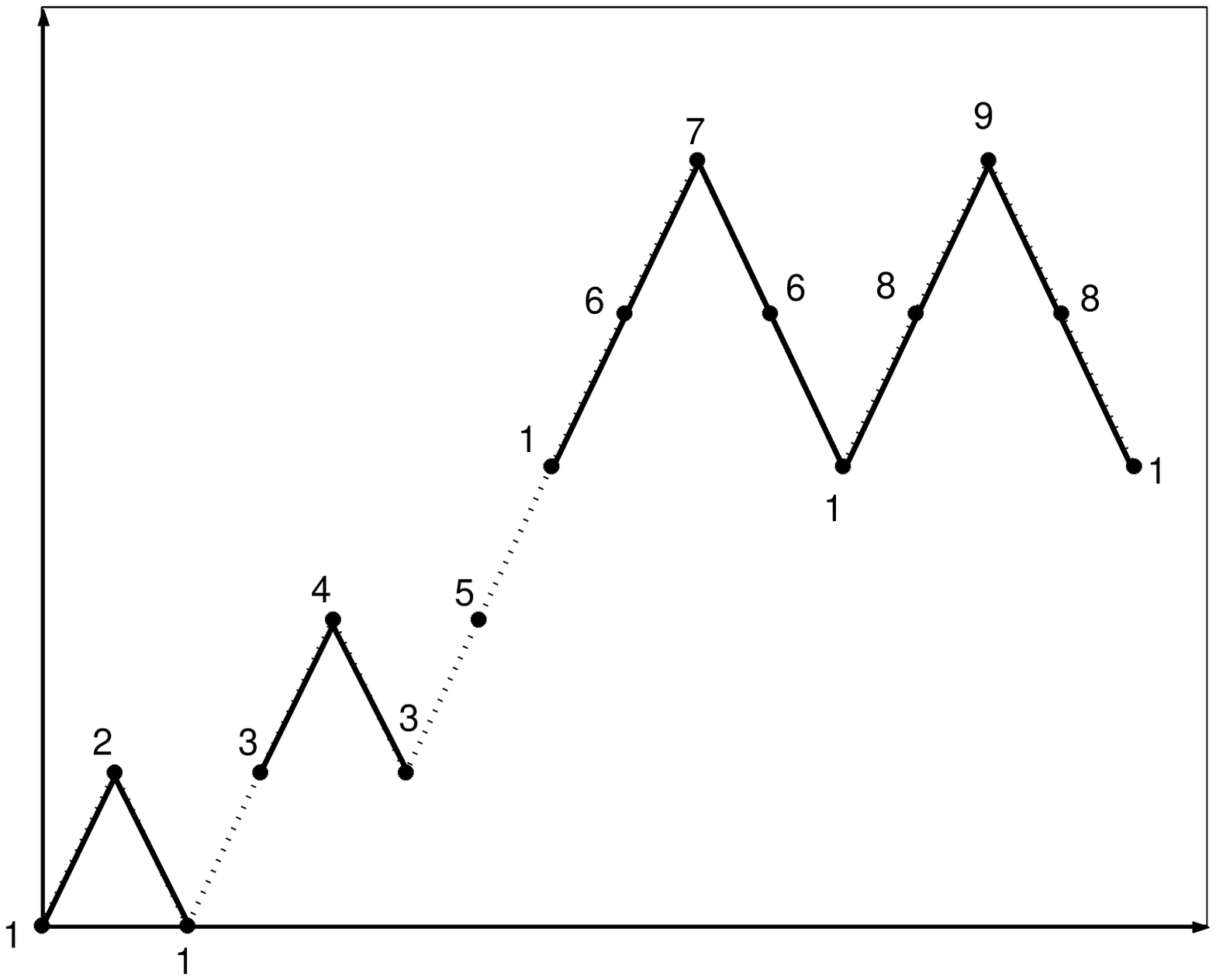} \\
\multicolumn{3}{p{13cm}}{{\small{Figure 1. \textit{Left fig.: }A path $\mathcal P \in
\mathcal{P}_{6,4}$ with origin $1$ unmarked. $\mathcal{N}_1=\{2,5,7,8\}\text{ and}\hspace*{.5cm}$
$\mathcal{N}_2=\{3,4,6\}$. \textit{Right fig.: }A simple path $\mathcal P$ of ${\mathcal P}_{6,3}$
where the origin $1$ is marked.}}}
\end{tabular}
\end{center}

\section{Asymptotics of $\mathbb E[Tr M_N^{2s_N (+1)}]$ for $1<<s_N<<\sqrt N
  $ and $\theta > \sigma$ } \label{Sec : s_n<<sqrtn}
The aim of this section is to prove the following theorem for any general complex or real Deformed
Wigner Ensemble of type $(i)-(iv)$ or $(i')-(iv')$.

\bt \label{trace} Assume $\theta > \sigma $. Then, for any sequence $1 << s_N << \sqrt{N}$ and for
$L_N=2s_N(+1)$ \be  \label{trws} \mathbb E[tr M_N^{L_N} ]= \rho_{\theta}^{L_N}(1+o(1)). \ee \et
For the sake of clarity, in the whole paper, we only consider traces of even powers of $M_N$, since
the reasoning is exactly the same for odd exponents. Theorem \ref{trace} must be compared with the
analogous result established in \cite{Si-So1} for classical Wigner Ensembles (i.e. $\theta=0$).

\bt \cite{Si-So1}\label{theo : soshni1} For any sequence $s_N$ such that $1 << s_N << \sqrt{N}$,
\be  \label{trw1} \mathbb E\big [tr \left (\frac{W_N}{\sqrt N}\right)^{2s_N+1} \big]= 0 \quad
\text{and} \quad \mathbb E \big [tr \left (\frac{W_N}{\sqrt N}\right)^{2s_N} \big ]=
\frac{N}{\sqrt{\pi} s_N^{3/2}} {(2\sigma)}^{2s_N}(1+o(1)). \ee \et

\brem \label{Rem : cases theta leq sigma}The strategy used in the sequel also leads to the
following universal estimates, very close to those of Theorem \ref{theo : soshni1}, in the case
where $\theta\leq \sigma.$ One has
\begin{eqnarray*}
&&\mathbb E\big[tr \left(\frac{M_N}{2\sigma}\right)^{2s_N+1}\big ]= o(1) \quad \text{and} \quad
\mathbb E\big[tr \left(\frac{M_N}{2\sigma}\right)^{2s_N}\big ]= \frac{N}{\sqrt{\pi}
s_N^{3/2}}(1+o(1)),  \text{ if }\theta<\sigma,\quad \quad\crcr && \mathbb{E}\big[tr
\left(\frac{M_N}{2\sigma}\right)^{2s_N+1}\big ] =(1/2+o(1)) \quad \text{and } \quad \mathbb
E\big[tr \left(\frac{M_N}{2\sigma}\right)^{2s_N}\big ]= \frac{N}{\sqrt{\pi} s_N^{3/2}}(1+o(1)),
\text{ if }\theta=\sigma.
\end{eqnarray*}
\erem Our proof of Theorem \ref{trace} often refers to that of Theorem \ref{theo : soshni1} which
we briefly recall. The first estimate of (\ref{trw1}) follows from the symmetry and independence
assumptions on the Wigner matrix entries. It is proved in \cite{Si-So1} that the main contribution
to the second expectation in (\ref{trw1}) comes from even simple paths with an unmarked origin, as
in the proof of the classical Wigner Theorem (c.f. \cite{Bai}). As each of these paths has exactly
$s_N$ edges passed twice, once in one direction and once in the reverse direction, it is uniquely
determined by a trajectory $x$ of $\mathcal{T}_{s_N,0}$, the gift of the origin and of the
vertices at marked instants. One then readily deduces that their total contribution is of the order of
$\displaystyle{T_{s_N,o} \: N^{s_N+1} \: (\frac{\sigma ^2}{N})^{s_N}=\frac{(2s_N)!}{s_N!
(s_N+1)!} \: N
 \: \sigma ^{2s_N}.}$
Stirling's formula yields then the result.

For the proof of our Theorem \ref{trace},  we shall examine paths of the whole $\{ \mathcal P
_{m,l}, \, l+2m=2s_N, \, l \geq 0 \}$. Theorem \ref{theo : soshni1} actually gives the contribution
of the even paths ($l=0$) of any Deformed model $M_N$ . Indeed, one has that
$\mathbb{E}(|M_{i,j}|^2)=\frac{1}{N}(\sigma^2+\frac{\theta^2}{N}), \forall i,j \leq N$. But, in
the whole paper, we will replace $\sigma^2+\frac{\theta^2}{N}$ with $\sigma^2$ (for any $\theta
>0$), since, in no cases, the error made will affect the final result. Then, as $\rho _{\theta} >
2 \sigma$ if $\theta > \sigma$, the contribution of even paths is negligible (compare with
$(\ref{trws})$). The investigation of paths with at least one unreturned edge (i.e. $l >0$) is
quite different. In particular, the origin can here be marked in typical paths (i.e. those giving
the main contribution to the expectation). In fact, to obtain the precise estimate of Theorem
\ref{trace}, we need to refine the counting procedure of \cite{Si-So1} (because of the ``$l$ not
returned edges''). We shall also consider separately the cases where the origin $i_o$ is marked or
not. We prove that, when the origin $i_o$ is marked, simple paths of $\{\mathcal{P}_{m,l}, \, l
\geq 1 \}$ are typical (c.f. Subsection \ref{subsec: iomarque}). In the case where $i_o$ is
unmarked, we establish that paths with only one simple self-intersection are typical. This later
result requires a finer study which will allow us to boil down to paths with a marked origin (c.f.
Subsection \ref{sub: corresp}).

\subsection{Paths with marked origin $i_o$ \label{subsec: iomarque}}
Let us first compute the contribution of simple paths. Consider such a simple path $\mathcal{P}$
of length $2s_N$ belonging to some $\mathcal{P}_{m,l}$, with $l \geq 1$ and $l+2m=2s_N$. Since
$\mathcal{P}$ belongs to $\mathcal P_{m,l}$, it has exactly  $l+m$ marked instants and $m$
unmarked instants. As $\mathcal P$ is simple, it has exactly $l$ edges that appear once and $m$
edges that appear twice, once in one direction and once in the other direction. Thus the
contribution of $\mathcal P$ to the expectation $\mathbb E[tr M_N^{2s_N}]$ is equal (at the leading
order) to $\frac{\sigma^{2m} \times \theta^l}{N^{l+m}}$. The last point is that a simple path of
$\mathcal P_{m,l}$ with a marked origin is uniquely determined (see Remark
\ref{rempositionorigsimplepath} below) by a trajectory $x$ of $\mathcal T_{m,l}$ and the $l+m$
distinct values at its marked instants. From this we deduce (as $\theta > \sigma$) that the total
contribution of simple paths with marked origin is $\displaystyle{\mathbb E[ \text{ simple paths
with a marked origin}] = \sum_{l>0,\: l \text{ even}}T_{m,l}\theta^l
\sigma^{2m}e^{-\frac{(l+m)^2}{2N}}=
  \rho_{\theta}^{2s_N}(1-\frac{\sigma^2}{\theta^2})(1+o(1)).}$
  \brem \label{justif cases theta leq sigma} One has that $\sum_{m=0}^{s_N} T_{m,l} \theta^l \, \sigma^{2m} = o((
2\sigma)^{2s_N})$ if $\theta < \sigma$ and $\sum_{m=0}^{s_N} T_{m,l} \theta^l \, \sigma^{2m} = (2
\sigma) ^{2s_N}/2 $ if $\theta = \sigma$. This fact combined with the following analysis then
justifies Remark \ref{Rem : cases theta leq sigma}. \erem

\brem \label{rempositionorigsimplepath} If $\mathcal{P}$ is a simple path of length $2s_N$ with
marked origin, it is not hard to see that the instant $T_o$ of the marked occurence of the origin
is uniquely determined. If the last step of $\mathcal{P}$ is up, then $T_o= 2s_N$. Otherwise
$T_o=\inf \{ t > 0, \, x(t)=l \, \text{ and }
\forall t'\geq t, \, x(t') \geq l \}$ (see Figure 1 above).\\
\erem

The following theorem shows that typical paths with marked origin are simple if $s_N =o(\sqrt{N})$.

\bt \label{tracewsi} Assume $\theta > \sigma $. Then, for any sequence $1 << s_N << \sqrt{N}$ and
for $L_N=2s_N(+1)$ \be  \label{trwsi}
\mathbb E[ \text{ Paths with a marked origin}]= \rho_{\theta}^{L_N}(1-\frac{\sigma^2}{\theta^2})(1+o(1)).\\
\ee \et

We shall now show that, amongst paths with a marked origin, paths with some multiple self
intersections give a negligible contribution to the expectation. For such paths, one can observe
that the instant of the marked occurrence of the origin is not determined any more and there are
multiple ways to close edges opened previously. Throughout the rest of this subsection,  we only
consider paths with a last step up. The case of paths with marked origin and a last step down will
be studied at the end of the Subsection 4.2. Assuming the last step is up avoids technicalities
(see Subsection 4.2). Indeed, when the last step is up, the origin is well defined, once the
vertices at marked instants have been chosen. Thus, with respect to the analysis made in
\cite{Si-So1}, one essentially only has to pay attention to the unreturned edges of the path. Such
edges have to be taken into account in the estimation of the closing of the path as we explain now.

\subsubsection{Closing of the path when the last step is up.}
Before proving Theorem \ref{tracewsi}, we give an important technical result which is analogous to
Lemma 1 of \cite{Si-So2}. Let $\mathcal{P}$ be a path in $\mathcal P_{m,l}$ of type $(N_o,N_1,
\cdots, N_{l+m})$ and whose last step is up. By definition of the $N_k$ (see Section 3), one can
readily verify the following relations
 \begin{eqnarray}{\label{type}}
&&\sum_{k=0}^{l+m} N_k=N, \quad \sum_{k=1}^{l+m} kN_k=l+m.
      \end{eqnarray}
Assume that we have chosen the distinct vertices occuring at the marked instants of $\mathcal{P}$.
Call then $\Omega _{m}$ the number of ways to fill in the blanks of $\mathcal{P}$ at the unmarked
instants (notice that in references \cite{Si-So1}, \cite{Si-So2} and \cite{Sos}, the authors use
the notation $W_m$). We then set $\Omega _m \mathbb E _{max}:=\underset{\mathcal{P} \text{ of type
}(N_o,N_1, \cdots, N_{l+m})}{\max}\Omega _m |\mathbb{E} \Big [\prod_{j =0}^{2s_N-1} M_{i_ji_{j+1}}
\Big ]|$.

\bp \label{Prop: majoWmE} There exists $C>0$ independent of $N, s_N$ and $m$ such that
\begin{eqnarray}\label{majoWmE}
\Omega _m \mathbb E _{max}  & \leq & \prod_{k= 2}^{10}\left
(Ck\right)^{kN_k}\prod_{k=11}^{l+m}\left (Ck \right ) ^{4kN_k/3} \, \frac{\theta^l
\sigma^{2m}}{N^{l+m}}
\end{eqnarray}
\ep
\paragraph{Proof of Proposition \ref{Prop: majoWmE}:} The main difference from the proof of Lemma 1 in
\cite{Si-So2} follows from the existence of the $l$ odd edges. One can readily check that $\Omega_m\leq \prod_{k\geq 2}(2k)^{kN_k}.$ Then, using the same arguments as in
\cite{Si-So2}, we obtain that
\begin{eqnarray}
&& \Omega _m |\mathbb{E}\prod_{j =0}^{2s_N} M_{i_ji_{j+1}}|\cr
&&\leq\frac{\sigma^{2m}\theta^l}{N^{l+m}} \prod_{(i'j')\text{ even, } l(i'j')\geq 2}\frac{(4\beta
l(i'j'))^{l(i'j')}}{l(i'j')!}
\prod_{(ij)\text{ odd, }l(ij)\geq 2} \frac{(4\beta l(ij))^{l(ij)+1} }{l(ij)!} \prod_{k\geq 2}4^{kN_k}\prod_{k\geq 2}(2k)^{kN_k}\cr &&\label{contributionnonref}\\
&&\leq \frac{\sigma^{2m}\theta^l}{N^{l+m}}\prod_{k=2}^s\left (C_1k\right )^{kN_k}\prod_{
(ij) \text{ odd}, \:l(ij)\geq 2}l_{ij}.\nonumber
\end{eqnarray}
In (\ref{contributionnonref}), we have used that if an edge $(i'j')$ (resp. $(ij)$) occurs an even (resp. odd ) number of times $2l(i'j')$ (resp. $2l(ij)+1$), then the path is closed $l(i'j')$ (resp. $l(ij)$) times along the same edge. The term $\prod_{k\geq 2}4^{kN_k}$ comes from edges read three times. For such edges $(ij)$, one can check that $i$ or $j$ is necessarily a self intersection and that $\mathbb{E}(M_{ij}^3)\leq 4\sigma^{2}\theta/N^2$ (and $|\mathbb{E}(\overline{M_{ij}}^2M_{ij})|\leq 4\sigma^{2}\theta/N^2$).\\
Now, let $(ij)$ be an odd edge of $\mathcal P$ for which $l(ij)\geq 2.$ Then
$l(ij)=k(ij)+k(ji)-1$, where $k(ij)$ denotes the number of times $i$ is marked in the edge $(ij).$
Then $l(ij)\leq 2\max \{k(ij), k(ji)\}$ and denoting by $m(ij):=\max \{k(ij), k(ji)\}$ we have that
\begin{eqnarray}&\displaystyle{\prod_{(ij) \text{ odd, } l(ij)\geq 2}l(ij)}&\leq \prod_{(ij) \text{ odd},\: m(ij)\geq 2}m(ij)\prod_{k\geq 2}2^{kN_k}\cr
&&\leq \prod_{k\geq 2}2^{kN_k}\prod_{(ij) \text{ odd},\: 2\leq m(ij)\leq
3}m(ij)\prod_{(ij) \text{ odd},\: 3\leq m(ij)}m(ij). \label{majolij}
\end{eqnarray}
Now, for a vertex of self intersection $i$, we denote by $p(i)$ the number of odd edges $(ij_l), l=1, \ldots, p(i),$ such
that
$$k(ij_l)=\max\{k(ij_l), k(j_li)\} \text{ and } k(ij_l)\geq 3.$$
Then, denoting by $k_i$ the number of times $i$ is marked in $\mathcal P$,
$\displaystyle{\prod_{l=1}^{p(i)} k(ij_l)\leq k_i^{p(i)}\leq k_i^{k_i/3},}$ since necessarily $k_i\geq 3p(i).$
Inserting this in (\ref{majolij}), we obtain that
 $\displaystyle{\prod_{\text{edges }(ij) \text{ odd}}l(ij)\leq \prod_{k\geq 2}6^{kN_k}\prod_{k\geq 3}k^{kN_k/3}.}$
Thus \\
$\displaystyle{\Omega _m |\mathbb{E}\prod_{j =0}^{2s_N} M_{i_ji_{j+1}}|  \leq 
 \prod_{k\geq 2}\left (6C_1k\right)^{kN_k}\prod_{k\geq 3}
 k^{kN_k/3}\frac{\theta^l \sigma^{2m}}{N^{l+m}}.}$
We then readily deduce $(\ref{majoWmE})$. $\square$
\brem \label{rem: prem avril}One can note that Proposition \ref{Prop: majoWmE} also holds, up to minor modifications, for paths with last step down if, in this case, $\Omega_m$ denotes the number of ways to fill in the blanks of $\mathcal{P}$ once vertices at the origin and marked instants are given.
\erem

\subsubsection{Contribution of paths with last step up and self intersections }
We now come back to the proof of Theorem \ref{tracewsi} for paths with last step up and prove that
paths with self intersections give a
  negligible contribution to the expectation. Consider $(l,m)$ such that
  $l+2m=2s_N$. Given $(N_o,N_1, \cdots,
N_{l+m})$ satisfying $(\ref{type})$, the number of ways to distribute the $l+m$ marked instants of
such a path is $\frac{(l+m)!}{\prod_{k\geq 1} (k!)^{N_k}}$ and the number of ways to affect the
vertices is $\frac{N!}{{N}_o! {N}_1! \cdots  N_{l+m}! }.$ Once the marked vertices and origin are
chosen, the contribution to the expectation of paths of $\mathcal P_{m,l}$ of type $(N_o,N_1,
\cdots, N_{l+m})$ with last step up can be bounded from above by the r.h.s. of $(\ref{majoWmE})$.
Accordingly, the added contribution of paths of $\mathcal P_{m,l}$  of type $(N_o,N_1,\ldots,
N_{l+m})$ and having a last step up
 can be estimated from above by
\begin{eqnarray}
&&\frac{N!}{N_o!N_1!\ldots N_{l+m}!} \:T_{m,l-1}\:
  \frac{(l+m)!}{\prod_{k=2}^{l+m}(k!)^{N_k}} \Omega _m \mathbb{E}_{max}\cr
&&\leq  T_{m,l-1}\theta^l \sigma^{2m}
  N^{N-N_o-(l+m)}\prod_{k=2}^{10}\frac{\Bigl(C(m+l)\Bigr)^{kN_k}}{N_k!}
  \prod_{k=11}^{l+m}\frac{\Bigl((C(m+l))^{4/3}\Bigr)^{kN_k}}{N_k!}. \label{majofact}
\end{eqnarray}
 The last
inequality follows from Proposition \ref{Prop: majoWmE} (and $T_{m,l-1}=\# \{x \in
\mathcal T _{m,l-1} \, {\text{with a last step up}} \}$). As before, $C$ denotes a positive
constant whose value may change from line to line.\\
Moreover, as $N-N_o-(l+m)=- \sum_{k=2}^{l+m} (k-1)  N _k $ (by $(\ref{type})$), one obtains that\be
\label{majofactbis} (\ref{majofact}) \leq   T_{m,l-1}\theta^l \sigma^{2m} \,
\prod_{k=2}^{10}\frac{1}{N_k!} \Bigl
  (\frac{(C(m+l))^k}{N^{k-1}}\Bigr)^{N_k}\prod_{k=11}^{l+m}\frac{1}{N_k!} \Bigl
  (\frac{(C(m+l))^{4k/3}}{N^{k-1}}\Bigr)^{N_k}.
\ee Set now $\sum_{k =2}^{l+m}N_k=M_1+M_2 \text{ with
  } \,  M_1=\sum_{k \geq 11 }N_k \text{ and } \, M_2=\sum_{k=2}^{10}N_k.$
We have, as long as $s_N=O(N^{2/3})$ (which will be the greatest scale of use in this paper) that
\begin{eqnarray} \label{majofactM1}
&&\sum_{(N_k, k\geq 11) /  \, \sum_{k \geq 11 }N_k=M_1} \:\prod_{k\geq 11}\frac{1}{N_k!} \Bigl
  (\frac{(C(m+l))^{4k/3}}{N^{k-1}}\Bigr)^{N_k} \crcr
  && =   \frac{1}{M_1!}\sum_{(N_k, k\geq 11) / \: \sum_{k \geq 11 }N_k=M_1} \:\frac{M_1!}{\prod_{k\geq 11}N_k!} \prod_{k\geq 11}\Bigl
  (\frac{(C(m+l))^{4k/3}}{N^{k-1}}\Bigr)^{N_k}
  \leq    \frac{(2C N ^{-2/9})^{M_1}}{M_1!} .
\end{eqnarray}
Similarly, one obtains that \be \label{majofactM2} \sum_{ (N_k, k=2, \ldots, 10) / \, \sum_{2 \leq
k \leq 10 }N_k=M_2} \:\frac{1}{N_k!} \Bigl
  (\frac{(C(m+l))^{k}}{N^{k-1}}\Bigr)^{N_k}
 \leq  \frac{1}{M_2!} \left ( \frac{C s_N^2}{N}\right)^{M_2}.
\ee Recalling that paths with multiple self-intersections are such that $0 < \sum_{k=2}^{l+m}k
N_k$, the summation of (\ref{majofactbis}) over all the $(N_0, \cdots, N_{l+m})$ such that $M_1>0$
(and any $M_2 \geq 0$) is not greater than

\be \label{majoMavec expinit} T_{m,l-1} \: \theta^l \, \sigma^{2m} \: \times \exp \{ \frac{(C
(l+m))^2}{N} \}  \times \left
       [\exp{(C N^{-2/9})} -1\right ] .
\ee As $\exp{(C N^{-2/9})} -1=O(N^{-2/9})=o(1)$, this ensures that the contribution of paths for
which there exists a vertex of $k-$fold self-intersection such that $k\geq 11$ is negligible. At
last, we shall consider paths with multiple self intersection of type smaller than $10$ which
amounts to consider the summation of (\ref{majofactbis}) over all the $(N_0, \cdots, N_{l+m})$
such that $M_2>0$. It is not greater than
\be \label{majoMavec exp} T_{m,l-1} \: \theta^l \, \sigma^{2m} \: \times \left [\exp \{\frac{(C
(l+m))^2}{N} \}-1 \right ] . \ee 
As $s_N=o(N^{1/2})$, one has that $\exp \{\frac{(C
(l+m))^2}{N} \}-1=o(1)$ uniformly in $l$. As a consequence (and noticing that $\sum _{l \geq
2}T_{m,l-1} \theta^l \sigma^{2m}=(\theta - \sigma ^2 /\theta) \rho _{\theta} ^{2s_N-1} (1+o(1))$),
we deduce that the summation of $(\ref{majoMavec exp})$ over all $l$ is negligible with respect to
$(\ref{trwsi})$. Theorem $\ref{tracewsi}$ for paths with last step up is established. $\square$


\subsection{\label{sub: corresp} Paths with a last step down }

In this section, we still consider paths where at least one edge is passed an odd number of times
(i.e. $l>0$). Indeed, even paths ($l=0$) are considered in Theorem $\ref{theo : soshni1}$ and give
a negligible contribution (as $\theta > \sigma$). We first investigate paths whose origin is
unmarked and establish the following result.

\bt \label{theo : iononmaruqesimin} Assume $\theta > \sigma $. For all $1 << s_N << \sqrt{N}$ and
for $L_N=2s_N(+1)$ \be  {\label{traceavec1si}} \mathbb E[ \text{ Paths with an unmarked origin}]=
\rho_{\theta}^{L_N}\frac{\sigma^2}{\theta^2}(1+o(1)). \ee \et

In the last part of this subsection, we will consider paths with last step down but a marked
origin and hence finish the proof of Theorem \ref{tracewsi}.

\paragraph{}
Before we proceed the proof of Theorem \ref{theo : iononmaruqesimin}, we give a sketch of our
reasoning. A path ${\cal P}\in {\mathcal P}_{m,l}$ with $l\geq 1$, $i_o$ unmarked and a last step
down has necessarily a self-intersection. But, as soon as there exist multiple self-intersections,
the geometry of the path becomes complex (compare with simple paths handling in the previous
subsection). For instance, a path of type $(N-(l+m), l+m, 1, 0, \cdots,0)$ with $i_o \in
\mathcal{N}_o$ (and $l>0$), is defined by a trajectory $x$ of $\mathcal{T} _{m,l}$, $l+m$ distinct
vertices and the instant of the second ``marked'' occurence of the vertex of simple
self-intersection. Analyzing in details the geometry of paths with last step down, we build a
correspondence, fundamental throughout the paper, between paths with last step down and unmarked
origin (resp. marked origin) and paths with marked origin and a last step up. In particular, the
correspondence is such that paths of type  $(N-(l+m), l+m, 1, 0, \cdots,0)$ with $i_o \in
\mathcal{N}_o$ are ``in bijection'' with simple paths having a last step up whose contribution was
precisely estimated in Subsection $4.1$. Thus paths with several self-intersections and  $i_o \in
\mathcal{N}_o$ are associated to paths with a marked origin and at least one multiple
self-intersection whose contribution is negligible according to
the reasoning of the previous subsection. Finally, the correspondence is used to consider paths with last step down and marked origin.\\

\noindent \textbf{Proof of Theorem \ref{theo : iononmaruqesimin}:}\\
\indent \textit{$1^{st}$ step: Construction of the correspondence. }Consider a term in the trace
$Tr M_N ^{2s}$ (for any $s$) \be M_{i_oi_1}M_{i_1i_2}M_{i_2i_3}\cdots M_{i_{2s-1}i_o}
\label{lepathP}\ee such that the corresponding path $\mathcal P \in \mathcal P _{m,l}$ has an
unmarked origin. Assume that the first odd edge of this path is read for the first time at instant
$2k+1$. Here we mean that the left endpoint of this edge occurs at time $2k$. Denote by $(vw)$
this oriented edge (i.e $v=i_{2k}$ and $w=i_{2k+1}$). Here we assume that the first odd edge is
read at an odd instant, but the reasoning is similar in the case where it is even. Consider then
the path $\mathcal P'$ defined by the term

\be M_{i_{2k+1}i_{2k+2}}\cdots M_{i_{2s-1}i_o}. M_{i_oi_1}M_{i_1i_2}M_{i_2i_3}\cdots
M_{i_{2k-1}i_{2k}}M_{i_{2k}i_{2k+1}}\label{lepathP'}, \ee which has the same edges, visited with
the same multiplicity as in $\mathcal P$ (that means that $\mathcal P$ and $\mathcal P'$ have the
same weight). Note that in $\mathcal P'$, the origin $w=i_{2k+1}$ is marked and determined by the
way we distribute the marked instants, since the last step of $\mathcal P'$ is up. Furthermore,
all the edges $(i_oi_1),(i_1i_2), \ldots, (i_{2k-1}i_{2k})$ are even in $\mathcal P$ as in
$\mathcal P'.$ Given $k$ and $\mathcal P'$, we can identify $\mathcal P.$
We simply move one by one the $2k+1$ last edges of $\mathcal P'$ to the beginning, reversing the operation leading from (\ref{lepathP}) to (\ref{lepathP'}), until we obtain a path with unmarked origin.\\
Let then $p$ be the number of edges opened but not closed before time $2k$ in $\mathcal P.$ Note
that $p\geq 1.$ Now we shall estimate the number of uplets $(\mathcal P', k)$ given the level
$1\leq p\leq m.$ Since the weight of the path $\mathcal P$ is equal to the one of the
corresponding path $\mathcal P'$, we only need to count the number of the possible underlying
trajectories $x'$ of $\mathcal P'$. It is not straightforward to make such a numbering and we
shall build a new transformation on the trajectory $x'$. To do this, we observe that in $\mathcal
P'$, the origin $i_o$ of $\mathcal P$ occurs at the instant $T:=2s-2k-1$. Then after time $T$, one
closes exactly $p$ edges opened before $T$ and the last edge is an up step. Thus after time $T$
where $i_o$ is at a level $l+p-1$, the trajectory of the path $\mathcal P'$ remains above the
level $l-1$ until time $2s-1$ and makes a final up step to reach the level $l.$ Let us call $x_o'$
be the subtrajectory of $\mathcal P'$ of length $2k+1$ in between the instants $[T,2s]$. We then
define $\tilde x'_o$ to be the trajectory $x'_o$ read in the reverse direction. Then $\tilde x'_o$
starts from $0$, makes a first down step and reaches level $p-1$ remaining above the level $-1$.
Finally, define $x''_o$ to be $\tilde x'_o$ where the first step is replaced with an up step. Thus
$x''_o$ is a path of $\mathcal T_{2k+1}^{p+1}$ which does not go below the level $1$ after the
first step.

Define now the trajectory $x''$ as follows: from time $0$ to time $T$, $x''$ coincides with $x'$
and then one reads the subtrajectory $x''_o.$
Then the following holds.\\
$\bullet$ $x''$ is a path of $\mathcal T_{m-p, l+2p}.$\\
$\bullet$ $T+1$ is the first time where $x''$ reaches the level $l+p$ without
  going below afterwards i.e.
$$T+1=\inf \{t, x''(t)=l+p, x''(t')\geq x''(t), \forall t'\geq t\}.$$
Thus, by the above procedure, we have, given the level $p$, built a bijection
between the uplets $(x', k)$ and the trajectories $x''$ of $\mathcal T_{m-p, l+2p}.$\\
\begin{center}
\begin{tabular}{cp{0.5cm}c}
\includegraphics[height=4cm,angle=0]{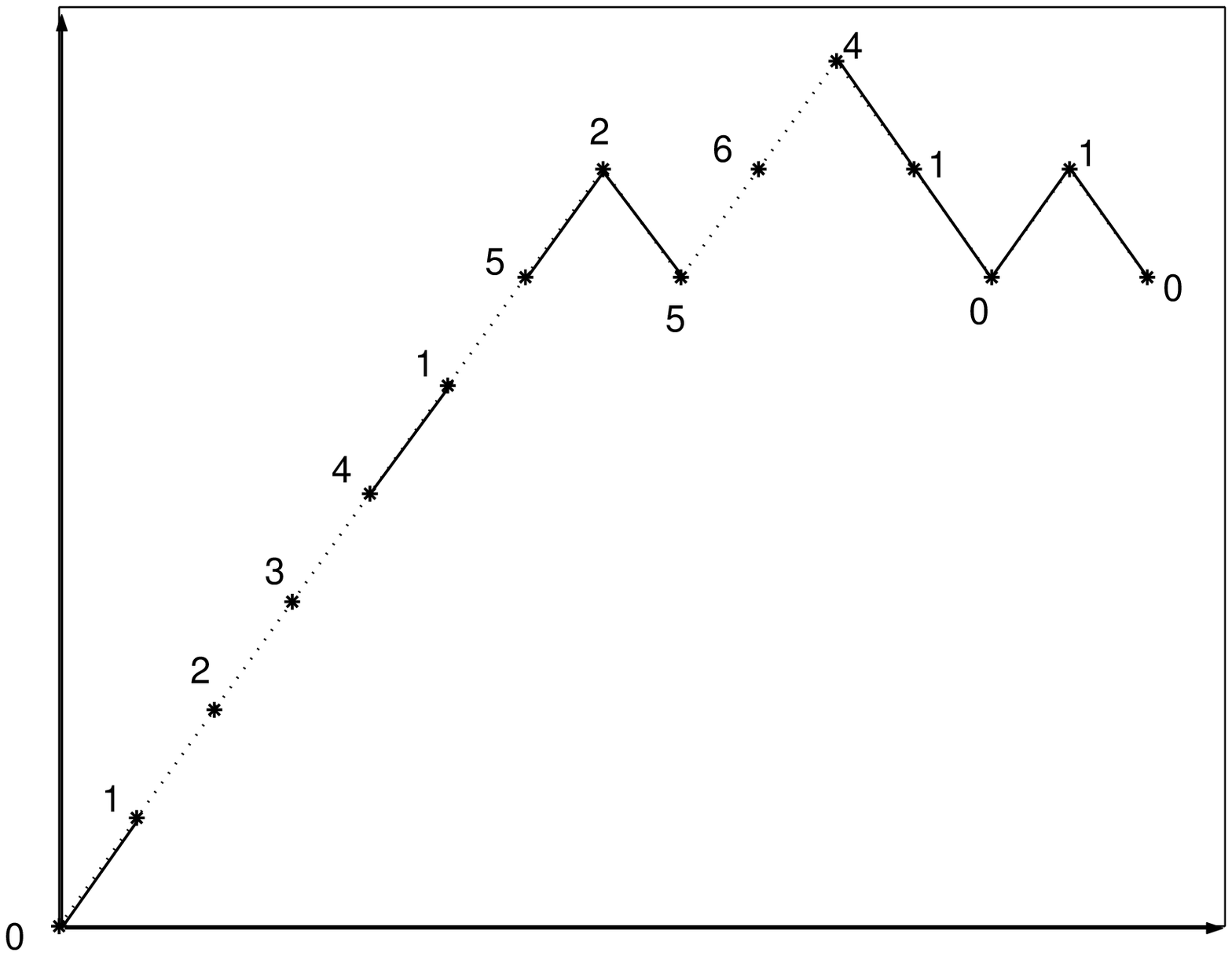} & &
\includegraphics[height=4cm,angle=0]{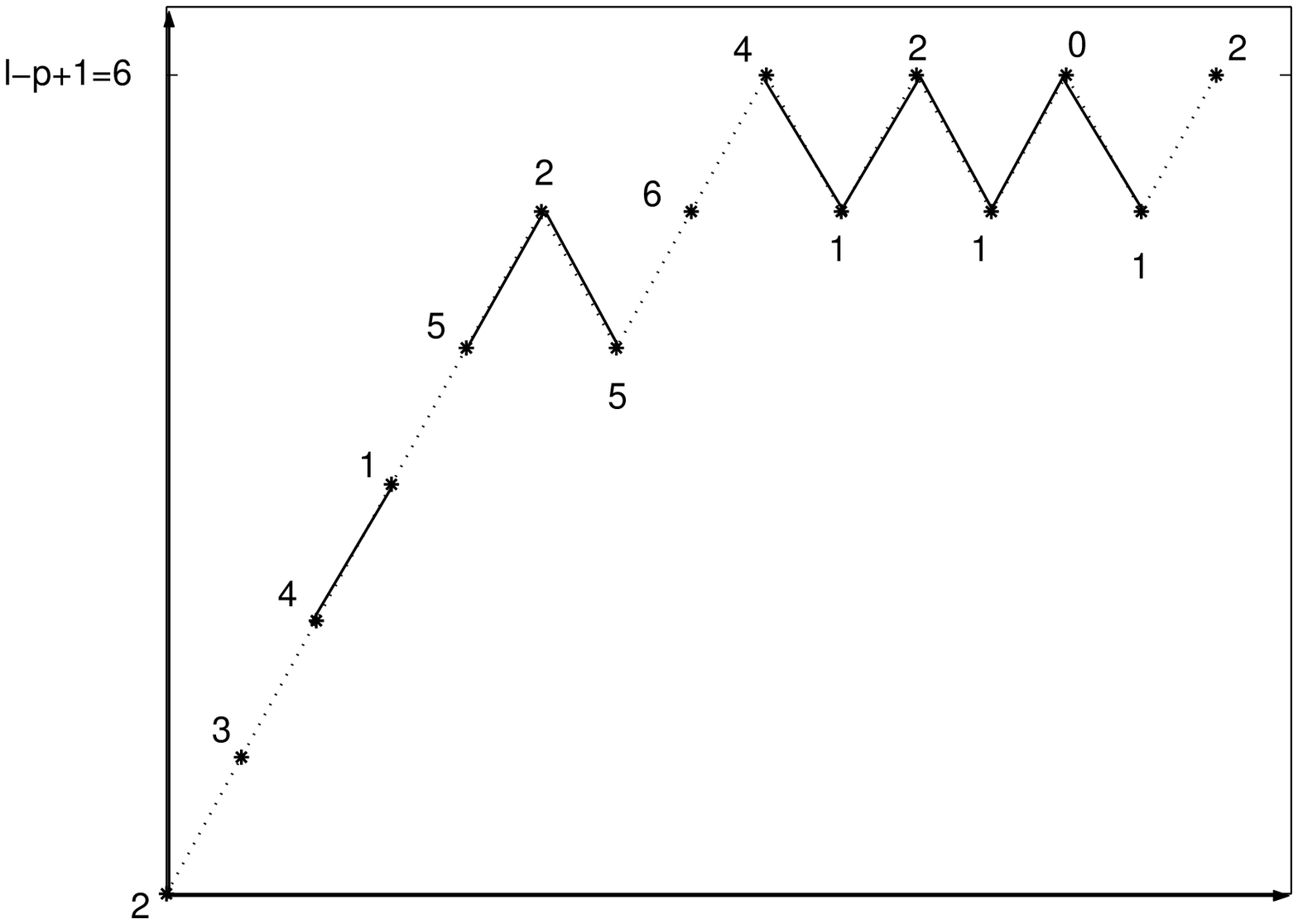} \\
\multicolumn{3}{p{14cm}}{{\small{Figure 2. \textit{Left:} A path $\mathcal P$ of $\mathcal P
_{4,6}$, with unmarked origin $0$, $\mathcal N _1=\{ 3,5,6\}$, $\mathcal N _2=\{2 ,4\}$, $\mathcal
N_3=\{1\}$. The second edge (12) is the first odd edge,
  $v=1$ and $w=2$. \textit{Right:} This path $\mathcal P'\in \mathcal P
  _{4,6}$ is in correspondence with $\mathcal P$: $i_o=2$, $p=1$, $\mathcal N _1=\{0,1,3,5,6\}$, $\mathcal N _2=\{4,2\}$. }}}
\end{tabular}
\end{center}

\textit{$2^{nd}$ step: Contribution of paths with origin unmarked. } As a path $\mathcal P$ and
its corresponding path $\mathcal P'$ have the same weight, we deduce (using computations as in
Subsection 4.1 with $M_1= \sum _{k \geq 11} N_k$ and $M_2=\sum _{k=2}^{10} N_k$) that the
contribution of paths with unmarked origin is at most
\begin{eqnarray}
&&\sum_{m=1}^{2s_N}  \, \sum_{1\leq p\leq m} T_{m-p,l+2p}\sum_{N_0, N_1, \ldots N_{l+m}}
\frac{N!}{N_o!N_1!\ldots N_{l+m}!} \:  \frac{(l+m)!}{\prod_{k=2}^{l+m}(k!)^{N_k}} \Omega _m
  \mathbb{E}_{max}\crcr
  &&\leq \sum_{m=1}^{2s_N} \sum_{M_1, M_2} \, \sum_{1\leq p\leq m} T_{m-p,l+2p}\theta^l\sigma^{2m}
  \frac{1}{M_2!}\left( \frac{Cs_N^2}{N}\right)^{M_2}\frac{1}{M_1!}\left (2CN^{-2/9}\right)^{M_1}.
\end{eqnarray}
So, reasoning as in the proof of Theorem \ref{tracewsi} (and using formula $(\ref{sumTml})$), it
is not hard to see that, as $s_N=o(\sqrt N)$, the paths $\mathcal P$ which correspond to paths
$\mathcal P'$ for which there exists at least a vertex of self intersection give a negligible
contribution to the expectation of the trace. Note that, in general, we are not able to say
something on the type of $\mathcal P'$ with respect to that of $\mathcal P$. Nevertheless, in the
scale $s_N=o(\sqrt N)$, it will appear that paths $\mathcal P$ having at least two
self-intersections or at least one self intersection of type larger than $2$ are negligible.
Indeed, the following step shows that our correspondence establishes a bijection between the set
of paths $\mathcal P$ with an unmarked origin and a sole simple self-intersection and the set of
simple paths
$\mathcal P'$ with last step up.\\

\textit{$3^{rd}$ step: Contribution of paths with
  unmarked origin and a single simple self-intersection. }
Here, we show that our previous correspondence is such that paths with an unmarked origin and a
single simple self-intersection are in bijection with some paths with a marked origin and without
self-intersection. Then, using the fact that these later paths are
  readily to count, we will deduce that their contribution to the
expectation is the rhs of $(\ref{traceavec1si})$. \\
Let $\mathcal P=\{i_o, i_1, \cdots, i_{2s-1}, i_o\}$ in $\mathcal P _{m,l}$ having an unmarked
origin with a sole self-intersection $v$. The typical geometry of such a path can be precisely
described. There exist six distinguished instants $0 \leq t_o<t_1 \leq t'_1 <t_2 \leq t_2' <t_3
\leq 2s$ defined as follows:
\begin{itemize}
\item $t_o=\max \{t \geq 0, \, x(t)=0\}$ and $t_3=\inf\{t\geq 0, x(t)=l \text{ and } x(t')\geq l\: \forall\: t'\geq t\}$. On $[0,t_o]$ (resp. $[t_3, 2s]$), the path $\mathcal P$ describes a simple sub-Dyck
  path $\mathcal P _o$ (resp. $\mathcal P' _o$) with origin $i_o$.
\item $t_1$ is the first marked occurence of $v$. Call $\mathcal P _i$ the part of the path $\mathcal
  P$ in between
  $[t_o,t_1]$: all its edges are even and there is $p \geq 1$ edges
  $\{e_j, 1 \leq j \leq p \}$ opened before $t_1$ and closed after $t_1$.
\item $t_2$ is the second marked occurence of $v$. The
  instant $t_1'$ and $t'_2$ are given by\\
$t_1'= \inf \{t \geq t_1, \, x(t_1')=x(t_1) \quad \text{and} \quad x(t') > x(t_1'), \, \forall
  t'>t'_1 \}$\\
$t_2'= \inf \{t \geq t_2, \, x(t_2')=x(t_2) \quad \text{and} \quad x(t') < x(t_2'), \, \forall
  t'_2 <t' \leq t_3 \}$\\
and such that on $[t_j,t_j']$ ($j=1,2$), the path $\mathcal P$ describes a sub-simple Dyck
  path $\mathcal P _v^j$ with origin $v$.
Note that it may happen that $t_1'=t_1$ (and thus $\mathcal P ^1_v$ is empty) or $t_2'=t_2$.
\item Let us denote by $w$ the
  vertex occuring at the instant $t_1'+1$: $(vw)$ is the first odd (simple) edge
  of $\mathcal P$. On $[t'_1+1,t_2]$, the
  subpath $\mathcal P _i$ begins in $w$, ends at $v$ both by an up step and remains above the level $x(t'_1+1)$.
\item Consider now the rest of the path. Just before $t_2'$, one closes the edge $e_p$ in the reverse
  sense. Next we successively return in the reverse direction
the other edges $e_{p-1}, \cdots, e_1$ and thus reach $i_o$ closing
  $e_1$ at time $t_3$ ($t_2' < t_3 \leq 2s$). These returns can be interspersed with
sub-simple Dyck paths.
\end{itemize}

The edge $(vw)$ is the first odd (simple) edge of $\mathcal P$: this is the distinguished
unreturned edge defining the origin $w$ of the new path $\mathcal P'$ obtained from our
correspondence. It is easy to see that $\mathcal P'$ is simple with a marked origin $w$ well
determined since the last step is up, the vertex $i_o$ is also marked and well defined. Moreover,
the vertex $v$ is now of type 1 in $\mathcal P'$: its marked occurence at time $t_2$ in $\mathcal
P$ is its sole marked occurence in $\mathcal P'$ whereas its marked occurence at time $t_1$ in
$\mathcal P$ is changed by an unmarked occurence in $\mathcal P'$. Thus, according to the
\textit{$1^{rst}$ step} of this proof, the numbering of such simple paths $\mathcal P'$ is of \be
\label{sumTml} \sum_{p=1}^m T_{m-p, l+2p}=C_{2s_N}^{m-1}=C_{2s_N}^m-T_{m,l}. \ee Then, one easily
finds the r.h.s. of $(\ref{traceavec1si})$. This finishes the
proof of Theorem \ref{theo : iononmaruqesimin}. $\square$\\

Let us now illustrate our correspondence on a path whose origin is unmarked and which has only one
simple self-intersection.
\begin{center}
\begin{tabular}{cp{1cm}c}
\includegraphics[height=4.5cm,angle=0]{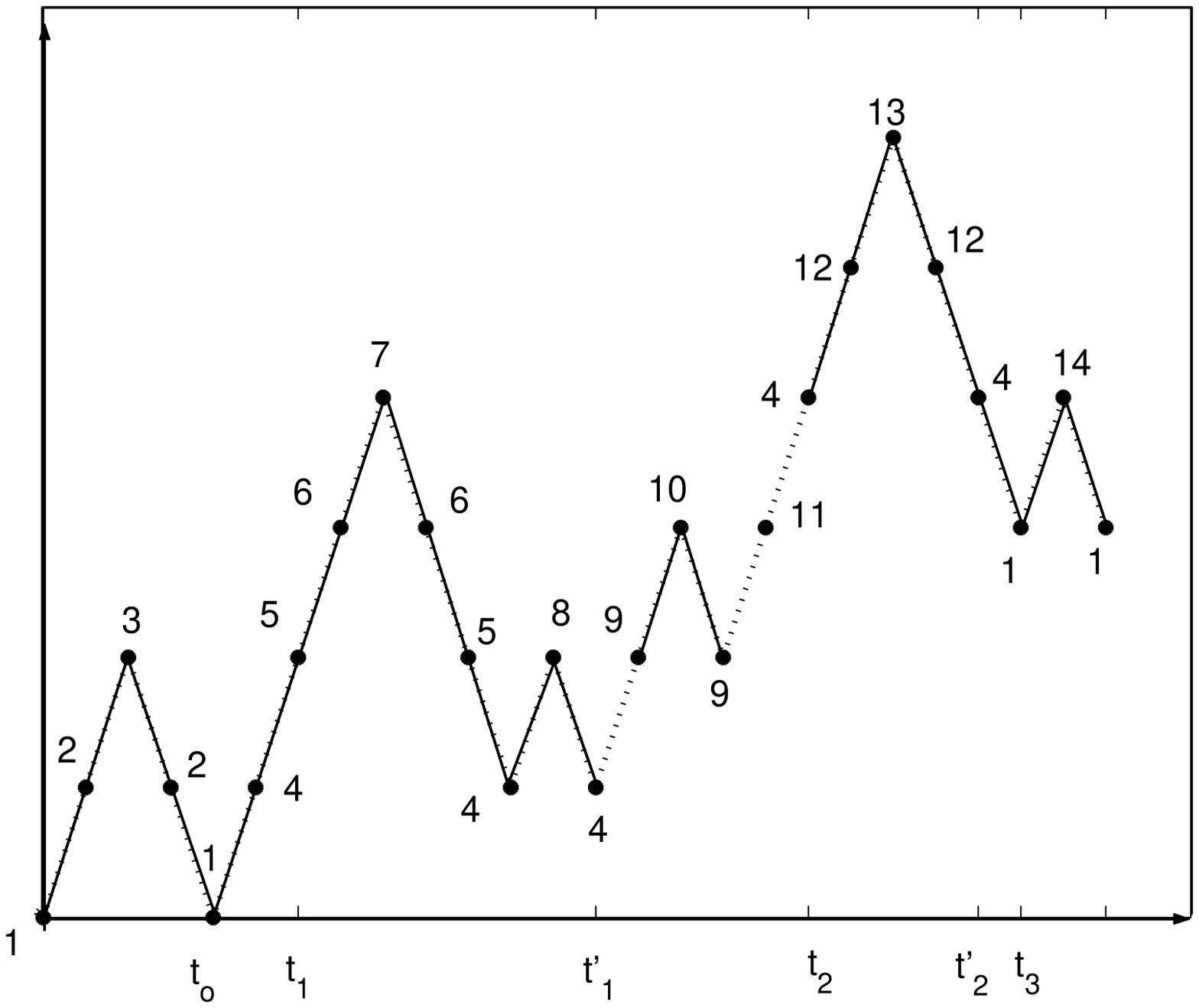} & &
\includegraphics[height=4.5cm,angle=0]{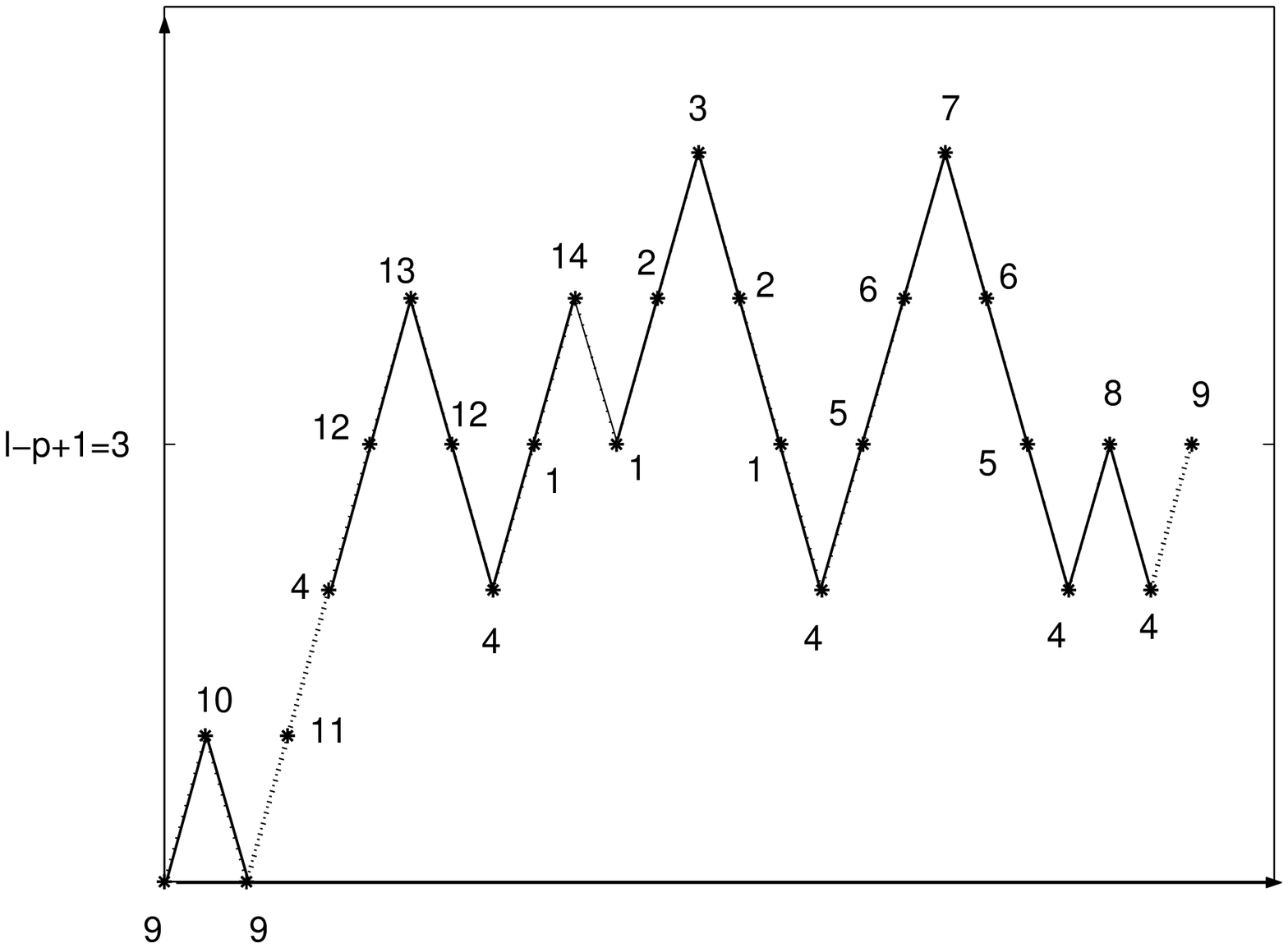} \\
\multicolumn{3}{p{13cm}}{{\small{Figure 3. \textit{Left: }$\mathcal P$ is in $\mathcal P _{11,3}$.
Its sole
  self-intersection is the vertex $4$ which is simple. $v=4$ and $w=9$. \textit{Right: }This is the simple path $\mathcal P'$  corresponding to $\mathcal P$. It belongs to $\mathcal P
  _{11,3}$ and $p=1$.  }}}
\end{tabular}
\end{center}

We now complete the proof of Theorem \ref{tracewsi} (and Theorem \ref{trace}). We shall then consider non-simple paths with marked origin and a last
step down. To this aim, it is
enough to notice that our correspondence still works for any path with marked origin and ending with a down step. The sole difference from the case where the
origin is unmarked is that the level $p$ (introduced in the $1^{st}$ step of the previous proof)
of the first odd edge can now be equal to $0.$ In this way, one can note that a path $\mathcal P$ which is not simple is associated to a non simple path $\mathcal P'.$ It is then easy to see (referring to the previous $2^{nd}$ step) that the non-simple   paths with last step down and marked origin give a negligible contribution to the expectation. We do not explain more.

\section{Computations of $\mathbb E[Tr M_N^{2s_N (+1)}]$ for $s_N=O(\sqrt N)$ if $\theta > \sigma$ \label{Sec : fluctuations}}
Here, we shall prove that, in the scale $s_N=O(\sqrt N)$ and as $N \to \infty$, the behavior of
the expectation of the Trace is the same for any Deformed Wigner Ensemble of type $(i)-(iv)$
(resp. $(i')-(iv')$). We use as before $M_N^G$ to denote the corresponding Deformed GUE (resp.
GOE) model.
\bt \label{Theo: expectation } Let $(L_N)$ be a sequence such that $\exists c>0, \, \lim_{N \to
\infty}\frac{L_N}{\sqrt N}=c$. Let $M_N$ be a Deformed Wigner matrix of type $(i)\!-(iv)$ (resp. $(i')\!-(iv')$). Then $\exists C'>0$ such that, for $N$ large enough,
$$\mathbb{E} \left [Tr M_N^{L_N}\right ]\leq C' \rho_{\theta}^{L_N} \, \text{ and
} \, \mathbb{E} \left [Tr M_N^{L_N}\right ]=\mathbb{E} \left [Tr \left (M_N^G
  \right )^{L_N}\right ] (1+o(1)).$$
 \et
To be more precise, in the complex setting, one has $\mathbb{E} \left [Tr M_N^{L_N}\right
  ]=\rho_{\theta}^{L_N}\exp{\{\frac{L_N^2}{2N}
(\frac{\sigma_{\theta}}{\rho_{\theta}})^2\}}(1+o(1)).$
This can trivially be deduced from the result of Theorem $\ref{ThPe}$  combined with some considerations of Section 2. Note that
similar exact estimates, with $\sigma_{\theta}$ replaced by $\sqrt{2} \sigma_{\theta}$, can be
expected for the real model (see Conjecture
\ref{conjectPecastheta><sigma}).\\
We only consider even powers $L_N=2s_N$ since the proof is similar for odd powers. The main part
of this section is dedicated to paths with a last step up. We show that the typical paths with a
last step up have at most simple self intersections, no loops and edges passed at most twice. The
last fact ensures in particular that the expectation of the Trace is the same for any Deformed
Wigner Ensemble. In Subsection \ref{sub : downstep}, thanks to the fundamental
correspondence built in Subsection 4.2, we translate this analysis to paths with a last step down
and show that universality holds too.

\subsection{Paths with last step up}{\label{sub: downstep}}
Throughout this section, we only consider paths with a last step up. Their contribution is at
least of the order of $\rho_{\theta}^{2s_N}$ since it can easily be seen from the preceding
section that
    the contribution of such simple paths is of the order of $\sum_{l\geq 2}T_{m,l-1}\theta^{l}\sigma^{2m}e^{\{-\frac{(l+m)^2}{2N}\}}\geq (1-(\frac{\sigma}{\theta})^2)e^{\{-\frac{2s_N^2}{N}\}} \rho_{\theta}^{2s_N}(1/3+o(1)).$ 
The following Proposition shows that there are at most simple self intersections in the typical
paths, i.e. those contributing to the Trace in a non negligible way. \bp \label{Prop :
notriple}Typical paths with last step up have no self intersection of multiplicity $k\geq 3$. \ep
\paragraph{Proof of Proposition \ref{Prop : notriple} :}Let $Z_o(m)$ denote the contribution of paths of $\mathcal P_{m,l}$ of type \\$(N_o, N_1, \ldots, N_{l+m})$
such that $\sum_{k\geq 3}N_k\geq 1.$ For such paths, set $M_1=\sum_{k>10} N_k$ and
$M'_2=\sum_{k=3}^{10} N_k$. Then, from (\ref{majofactbis}), one has that \be Z_o(m)\leq \sum_{N_2,
M_1, M'_2,M_1+M'_2>0 }T_{m,l-1} \, \theta ^{l} \sigma ^{2m} \frac{1}{N_{2}!}\left
(\frac{C(l+m)^2}{N}\right )^{N_{2}} \frac{1}{M'_2!} \Bigl(\frac{(C''(m+l))^3}{N^{2}}\Bigr)^{M'_2}
\frac{\left (CN^{-2}\right )^{M_1}}{M_1!}. \label{majofactbisbis} \ee The summation of
(\ref{majofactbisbis}) over all the $M_1,N_2, M'_2 $ such that $M'_2>0$ is not greater than \be
\label{majoMavec expbis} T_{m,l-1} \: \theta^l \, \sigma^{2m} \exp{\{\frac{(C (l+m))^2}{N}\}} \,
\Big [\exp{\{\frac{C''(l+m)^3}{N^2}\}}-1 \Big ]=T_{m,l-1} \: \theta^l \, \sigma^{2m} \times
O(\frac{(l+m)^3}{N^2}) \ee with $O(\frac{(l+m)^3}{N^2})= o(1)$ as $s_N=O(\sqrt{N}).$ Similarly,
the summation over all the $M_1,N_2, M'_2 $ such that $M_1>0$ is at most of the order $T_{m,l-1}
\: \theta^l \, \sigma^{2m} \times O(N^{-2}).$ Thus
$Z_o:=\sum_{m=0}^{s_N-1}Z_o(m)=O(\rho_{\theta}^{2s_N}\frac{1}{\sqrt N}),$
which is negligible w.r.t. the contribution of simple paths. $\, \square$\\

Given a path $\mathcal P$ of $\mathcal P _{m,l}$ (with $l+2m=2s_N$) with last step up and of type
$(N_o,N_1, \cdots, N_{l+m})$, we define $M:=
  \sum_{2}^{l+m} (k-1) N_k$ to be the number of its
  self-intersections. This quantity will be important in the following
  and is the object of the next proposition.

\bp \label{Prop : majoreM}The number of self intersections of typical paths satisfies
$M\leq s_N^{\epsilon},$ for any $0<\epsilon <1$. \ep

\paragraph{Proof of Proposition \ref{Prop : majoreM}:} By Proposition \ref{Prop : notriple} and (\ref{majofactbisbis}), it is clear that adding the contribution
of paths where ${\sum_{k\geq 2}k N_k\geq s_N^{\epsilon}}$ (for any $\epsilon >0$) gives a final
contribution which is
$o(\rho_{\theta}^{2s_N})$. $\square$\\

In the following, we investigate in details paths with a last step up and that have only simple
self intersections. Note that for such
paths, each edge is passed at most four times. We first discuss on those having edges read at most
twice. Then, we show that those admitting at least one edge passed three or four times and those
with at least one loop can be neglected.

\subsubsection{Paths with only simple self-intersections, edges read at most twice and last step up}
Our goal is here to prove that for a path with only simple self intersections, there exists
different ways of closing the path given the vertices at marked instant. In the denomination of
\cite{Sos}, this means that there are "non closed vertices" in typical paths. The definition will
be recalled later. This explains that the expectation of the Trace differs in the real and the
complex setting (see the comments just before Definition \ref{defnonclosed} below).

\paragraph{} Define $Z_1(M,m)$ to be the contribution of paths of type $(N_o, N_1, M,0, \ldots, 0)$ with a marked
origin, last step up and edges passed at most twice. Denote also by $Z_1(m)=\sum _M Z_1(M,m)$ the
total contribution of such paths. We want to establish that there exists a constant $D$
independent of $N$ such that, for $N$ large enough,
\begin{eqnarray}
Z_1:=\sum_{m\leq s_N-1}Z_1(m)& \leq & \exp{\{D\frac{s_N^2}{N}\}}\, \rho_{\theta}^{2s_N}.
\end{eqnarray}
Consider a path $\mathcal P$ contributing to $Z_1(M,m).$ By Propositions \ref{Prop : notriple} and
\ref{Prop :
  majoreM}, one can assume that $M=N_2\leq s_N^{\epsilon}$ for some arbitrary
$0<\epsilon<1$. Let then $t_{j_1}< t_{j_2}<  \cdots < t_{j_M} $ be the instants of
self-intersection of ${\mathcal P}$. We now choose the vertices occuring at the marked instants,
and thus fix the origin of the path. First, there are $\prod_{j=0}^{l+m-M-1}(N-j) \sim N^{l+m-M}
e^{\{-\frac{(l+m)^2}{2N}\}}$ different ways to choose the distinct vertices occuring in the path
in the order of their appearance. If a vertex of self intersection occurs at some instant
$t_{j_i}$, there are $j_i -i$ possible choices for such a vertex. It is indeed chosen amongst the
marked vertices, which have already occured in the path but have not yet been repeated. Note that
if $i_o$ is a vertex of self intersection, then $j_M=l+m$,
 $t_{j_M}=l+2m$ and there are at most $l+m-M$ choices for the vertex $i_o$.

\begin{remark}\label{Rem: nbclsoed}If there is no choice for closing edges at unmarked instants, then the number of such paths is at most $N^{l+m}$ since
 $$\sum_{M \geq 0}\sum_{j_i\leq l+m}\prod_{j=0}^{l+m-M-1}(N-j)
 \prod_{k=1}^M(j_k-k)=N^{l+m}e^{\{-(l+m)^2/2N\}}\sum_{M \geq 0}\frac{1}{M!}\left (
 \frac{(l+m)^2}{2N}\right )^M(1+o(1)).$$
Here the $o$ is uniform due
 to the fact that $M\leq s_N^{\epsilon}$ for some $0<\epsilon <1/32.$
 In the general case, there are many choices for closing edges from a
 vertex of self intersection and the number of paths of type $(N_o,N_1,N_2,0,\ldots,0)$ can then be of the order
 $N^{l+m}e^{C s_N^2/N}.$
 \erem
We now count the number of ways to close the path at unmarked instants. One can close an edge
starting from a vertex belonging to $\mathcal N _1$. In this case there is no choice for closing
it. We can also close an edge starting from a vertex in ${\mathcal
  N}_2.$ Then, we can close it in at most 3 ways: along the edge used to arrive at this
vertex for the first or second time, or along the edge used to leave it for the first time. Such
consideration leads to the notion of non-closed vertex. \bdefi {\label{defnonclosed}} A vertex of
self-intersection is said to be non closed if there are several possibilities of return from this
vertex at an unmarked instant. \edefi

For example, in the left path of Figure 1, the
vertex $3$ is closed whereas $4$ and $6$ are non closed.\\

Here we show that paths of type $(N_o,N_1, \cdots, N_{l+m})$ with non-closed vertices contribute
in a non negligible way to the expectation of the Trace, if $l>0$. The fact that typical paths
admit non closed vertices explains that the expectation of the Trace (and thus the limiting
distribution of $\lambda_1$) differs between the real and complex case. Indeed, assuming edges
appear at most twice, an oriented edge repeated with the same orientation has the weight
$\frac{\theta^2}{N^2}<<\frac{1}{N}$ in the complex case instead of $\frac{\sigma^2}{N}$ in the real case. \\
  Assume first that $i_o$ is of type one, so that the vertex at the origin is defined by the gift of
the distinct vertices occuring in the path. Then, by the definition of the associated trajectory
$x \in \mathcal T_{m,l}$, there are at most $x(t)$ ways of choosing a non-closed vertex of
self-intersection appearing at some instant $t$. Then the number of ways to choose the vertices
occuring at remaining marked instants, once the distinct vertices occuring in the path have been
chosen, is bounded from above by \be \label{nombreunclosed}\sum_{r=0}^M\sum_{1\leq
j_1<\cdots<j_{M}< l+m}\sum_{1\leq l_1<\cdots <l_r\leq M}(j_1-1)\overline{(j_{l_1}-l_1)} \cdots
\overline{(j_{l_r}-l_r)} \cdots (j_{M}-M)\prod_{i=1}^rx(t_{j_{l_i}}).\ee Here the overlining means
that the term does not appear in the expression. One then has that
$$(\ref{nombreunclosed})\leq \sum_{r=0}^M\frac{C_M^r}{M!}\left
  (\sum_{j=1}^{l+m}(j-1) \right )^{M-r} \left ( \sum_{j=1}^{l+m}
  x(t)\right)^r\leq \frac{1}{M!}\left ( \frac{(l+m)^2}{2}\right )^{M} \sum_{r=0}^M\frac{1}{r!}\left ( \frac{2M}{l+m}\max_{t}x(t)\right)^r.$$
Then, paths where $i_o$ is of type $2$ are negligible. Indeed, once the distinct vertices occuring
in the path have been chosen, one first chooses the vertex $i_o$ occuring at the marked instant
$t_{j_M}=l+2m$ and thus the origin of the path. Then one chooses the vertices occuring at the
remaining marked instants. The total number of ways to do so is at most
\begin{eqnarray} &&\label{nombreunclosed2}\sum_{r=0}^M\sum_{1\leq j_1<\cdots<\cdots <j_{M-1}< j_M=(l+m)}\sum_{1\leq l_1<\cdots <l_r< M}\prod_{i=1}^{M-1}(j_i-i)\prod_{i=1}^r\frac{x(t_{j_{l_i}})}{(j_{l_i}-l_i)}\: (l+m-M)\crcr
&&\leq \sum_{r=0}^M\frac{2M}{l+m}\frac{1}{M!}\left (
  \frac{(l+m)^2}{2}\right )^{M} \frac{1}{r!}\left (
  \frac{2M}{l+m}\max_{0\leq t\leq
    2s_N}x(t)\right)^r=(\ref{nombreunclosed}) \times o\left(\frac{1}{(l+m)^{1-\epsilon}}\right).
\end{eqnarray}

We now obtain an upper bound for $Z_1$. One has, in the general case, that $\max_{0\leq t\leq
2s_N}x(t)\leq l+m.$ This estimate implies that one can not neglect in (\ref{nombreunclosed}) the
different ways of closing the edges. Moreover, the contribution of a path with vertices of type at
most two, $r$ non-closed vertices and edges passed at most twice is bounded above by $
\frac{1}{N^{l+m}} \sigma ^{2m} \, \theta^l \times 3^r.$ Thus, using summation,
\begin{eqnarray}
& Z_1&\leq  \sum_{m=0}^{s_N-1}T_{m,l-1} \theta^{l}\sigma^{2m}  \sum_{M \leq
  s_N^{\epsilon}} \, \sum_{r\leq M}e^{\{-\frac{(l+m)^2}{2N}\}}\frac{\left (
    \frac{(l+m)^2}{2N}\right )^M}{M!}\frac{( 6M)^r}{r!}(1+o(\frac{1}{(l+m)^{1-\epsilon}})) \crcr
    &&\leq \theta\rho_{\theta}^{2s_N-1} \, \exp \Big (\frac{e^6 s_N^2}{N} \Big ). \label{estimZ_1}\end{eqnarray}

This proves that paths with edges passed at most twice and having possibly non-closed vertices
give a non-negligible but uniformly bounded contribution to the expectation $\mathbb{E}Tr\left (
\frac{M_N}{\rho _{\theta}} \right)^{2s_N}$. And $i_o$ is of type one in typical paths. $\square$

\brem The fact that there are non closed vertices in typical paths is a main difference from the
case where $l=0$ (see \cite{Si-So2}). Indeed, in that case and loosely speaking, ``$\max_t
x(t)=O(\sqrt {s_N})$'' implying that typical paths of length of order $\sqrt N$ do not have any
non-closed vertices (see Section 7 for some comments about this maximum). \erem

\subsubsection{\label{subsub: edges 3or4}Paths with edges passed four
  or three times and last step up}
Call $Z_2(M, m)$ the contribution of such paths $\mathcal P$ of $\mathcal P _{m,l}$ having only
simple self intersections and exactly $M$ self-intersections and where at least one edge is passed
three or four times. Let $Z_2=\sum_{M \leq s_N^{\epsilon}} \sum_{m\leq s_N-1} Z_2(M,m)$. We shall
now prove that $Z_2=o(Z_1)$. The considerations here are very close to those of the appendix of
\cite{Si-So2}. Nevertheless, in order to have a paper self-contained, we give the main steps of the
demonstration. \\
Assume that the distinct vertices occuring in $\mathcal P$ are known and that the origin of the
path is chosen (if $i_o$ is of type 2). Consider an unoriented edge $e=(vw)$ which is read (at
least) three times in the path. Two situations must be examined according to the directions of the
two up occurences of $e$.
\\The first one is when one reads two up oriented edges $(vw)$ (for instance) that is twice in the same direction. In this case, $w$ is a vertex of self intersection and is the rightendpoint of two edges started at $v$.
Assume first that $w\not=i_o$. We then need to introduce another characteristic of the path
$\mathcal P$, $\nu _N(\mathcal P)$, which is the maximum number of vertices that can be visited at
marked instants from a given vertex. If we denote by $t$ the (marked) instant of the second up
occurence of the edge $e$, then there are at most $ \nu _N(\mathcal P)$ choices for the vertex $w$
(since we shall look amongst the vertices already occured in the path and being in an up edge with
left endpoint $v$).
If $w=i_o$, then $t=l+2m$ and the vertex $i_o$ is then of type $2$ and occurs in an edge passed more than twice.\\
We now analyze the case where one first reads the up oriented edge $(vw)$ and then the second up
occurence of the edge $e$ in the reverse direction $(wv)$ (at the instant $t$). Then either $v$ is
a vertex of self intersection or $v=i_o$ is of type one and $t=l+2m$. In the case where $v$ is a
vertex of self intersection and $v\not=i_o$, the number of possible choices for the vertex of self
intersection at the marked instant $t$ (that is $v$) is at most the type of $w$, which is smaller
than 2 here. If $v=i_o$, then $i_o$ is of type $2$ and $t=l+2m.$ There remains to investigate the
case where $v=i_o$ but $i_o$ is of type one. Then the edge $wi_o$ is opened for the second time at
$t=2s_N$. In this case, the edge is read three times.
Note that this can happen only once and one can check that $\mathbb{E}\left (M_{ij}\right )^3$ (or $\mathbb{E}\left (|M_{ij}|^2M_{ij}\right )$) depends only on $\theta$ and $\sigma^2$ for any Deformed Wigner Ensemble. This possible event will thus not affect universality of the expectation and we will not consider it any more.\\
As before, paths where $i_o$ is of type 2 can be shown to be negligible with respect to those
where $i_o$ is of type one as before.
Thus it is now enough to obtain a bound for $\nu _N(\mathcal P)$.\\

It was shown in \cite{Si-So2} that, for even path $\mathcal P$ of length $2s_N$, $\nu _N(\mathcal
P)$ grows not faster than $s_N^{ \gamma}$, whatever $0 \leq \gamma <1$ is (at least for the
corresponding typical paths). The arguments of \cite{Si-So2} (also used in \cite{Sos}) can readily
be extended to our setting as explained in Remark \ref{rem: nun} below.  Thus, for the rest of
this section, it is enough to consider paths with only simple self intersections and such that
$\nu_N(\mathcal P) \leq (l+m)^{1/2-2\epsilon},$ for any value of $m$ and where $\epsilon<1/32.$
Let as above $t_{j_1}<t_{j_2}<\cdots <t_{j_M}$ be the instants of self intersection. Amongst these
instants, let $t_{u_1}<t_{u_2}<\cdots <t_{u_r}<l+2m$ and $t_{v_1}<t_{v_2}<\cdots<t_{v_q}<l+2m$ be
the instants respectively of nonclosed vertices and vertices belonging to edges passed three or
four times. Then, the contribution $\sum_m Z_2(M,m)$ of such paths is at most (overlining the
terms that do not appear in the product)
\begin{eqnarray*}
&&\sum_{m=0}^{s_N-1} \, \sum_{r= 0}^M \, \sum_{q =
1}^MT_{m,l-1}\theta^{l}\sigma^{2(m-q)}\sum_{1\leq j_1<j_2<\cdots <j_M \leq
l+m}\:\sum_{u_1<u_2<\cdots <u_r}\:\sum_{v_1<v_2<\cdots <v_q}\cr &&\times
\prod_{y_=0}^{l+m-M-1}(N-y)\Big ((j_{1}-1)\cdots \prod_{d=1}^r\overline
  {(j_{u_d}-u_d)}\prod_{i=1}^q\overline{(j_{v_i}-v_i)}\cdots (j_{M}-
  M)\Big )\crcr
  &&\times 3^r\prod_{1}^r x(t_{j_{u_d}}) \Big( C(l+m)^{1/2-2\epsilon}
  \Big)^q \frac{1}{N^{l+m}}\left(1+O((l+m)^{\epsilon-1})\right)\\\label{**}
& \leq & 2\sum_{m=0}^{s_N-1} \, T_{m,l-1}\theta^{l}\sigma^{2m}
e^{\{-\frac{(l+m)^2}{2N}\}}\frac{\left (
    \frac{(l+m)^2}{2N}\right )^M}{M!}\sum_{r\leq M}\frac{( 6M)^r}{r!}
 \Big( \exp ( \frac{C (l+m)^{3/2-\epsilon}}{N})-1 \Big).
\end{eqnarray*}
The term $\left(1+O((l+m)^{\epsilon-1})\right)$ comes from the fact that the origin $i_o$ can be
of type $2$, which gives a negligible contribution to the expectation as in
(\ref{nombreunclosed2}). Now, by summing over all $M$, we get (it is exactly as for the term
$Z_1$) $Z_2=o(Z_1)$.

\brem \label{rem: nun} We quickly explain why one can assume $\nu_N (\mathcal P)\leq
s_N^{1/2-2\epsilon}$. Consider a vertex $i$ which is the starting point of $\nu_N(\mathcal P)$ up
edges. Define $K=r+\sum_{k\geq 3}kN_k+1$. One can then split the interval $[0, 2s_N]$ into $K$
subintervals such that, inside such subintervals, edges can only be closed in the reverse
direction that the one used to open the edge. Furthermore, there exists such a subinterval in which
$i$ is the starting point of $\nu_N(\mathcal P)/K$ edges. In this interval, there are no choices
for closing the edges. Thus the trajectory falls $\nu_N(\mathcal P)/K$ times at the same level
(that of the first occurence, marked or not, of $i$ in this interval) without going below. It can
then be shown (see \cite{Sos} for more details, p. 41) that the number of such trajectories is not
greater than $\displaystyle{4 s_N^2\exp{\{-C_1 \nu_N(\mathcal P)/K\}}}T_{m,l-1},$ for some constant $C_1$ independent of $N, m$ and $s_N$. By Proposition
\ref{Prop : majoreM}, we can assume that $K\leq s_N^{\epsilon}<< s_N^{1/2-4\epsilon} .$ The
contribution of paths with $\nu_N(\mathcal P)\geq s_N^{1/2-2\epsilon}$, with $\epsilon<1/32$, is
then of the order
$$8 \sum_{m=0}^{s_N-1} \, T_{m,l-1}\theta^{l}\sigma^{2m}s_N^2\exp{\{ \frac{Cs_N^2}{N}+Cs_N^{\epsilon}\}}\sum_{\nu_N\geq s_N^{1/2-2\epsilon}}\exp{\{-C_1 \nu_N s_N^{-\epsilon}\}}=o(1)Z_1.$$
\erem

\subsubsection{\label{subsub: no loops} Paths with loops and last step up}
At this stage, we know that in the scale $s_N=O(\sqrt N)$, typical paths with last step up are of type $(N_o, N_1, M, 0, \cdots, 0)$ with $M=N_2\leq s_N^{\epsilon}$ (with $0
< \epsilon <1$) and have edges passed at most twice. To conclude to universality, we shall
state that those having at least one loop are negligible.\\
For this, given a path $\mathcal P=\{i_o, \cdots,i_{2s _N-1}, i_o\}$ of $\mathcal P _{m,l}$ with
$i_o$ marked and at most self-intersections of type $2$, we define $d$ to be the number of its
loops. We also introduce $u_1< \cdots < u_d$ the instants of self-intersection which lead to a
loop. That means that at the instant $u_j$, one opens a loop and $i_{u_j-1}=i_{u_j}$. So, there is
no choice for the vertices occuring at the $d$ instants $u_j$. The total contribution $Z'_1$ of
such paths is at most (it is similar to the approach leading to (\ref{estimZ_1})) 
\begin{eqnarray*}
&Z'_1 &\leq \sum_{m=0}^{s_N-1}T_{m,l-1}\theta^{l}\sigma^{2m}\sum_{M \leq
s_N^{\epsilon}}\sum_{d=1}^M  \sum_{r=0}^M \:\:\sum_{1\leq j_1<\cdots<j_{M}\leq l+m}\:\:\sum_{1\leq
  l_1<\cdots <l_r\leq M}\:\:\ \sum_{1\leq u_1<\cdots <u_d\leq M}\cr
  &&\frac{\prod_{i=1}^M (j_i-i)} {\prod_{i=1}^d (j_{u_i}-u_i) \times \prod_{i=1}^r
  (j_{l_i}-l_i)} 3^r \prod_{i=1}^rx(t_{j_{l_i}}) \frac{1}{N^M}\left(1+O((l+m)^{\epsilon-1})\right).
\end{eqnarray*}
Therefore, as $M=O(\frac{s_N^{2+\epsilon}}{N})$, we obtain that $\displaystyle{ Z'_1 \leq Z_1
\times \Big ( \exp {\{\frac{2s_N^{\epsilon}}{s_N}\}}-1 \Big )=o(Z_1).}$
\brem \label{casloopN3} In the whole generality, given a type of paths and assuming the origin is
marked (or unmarked), paths with loops have a negligible contribution with respect to those
without loop. \erem
\subsection{Paths with a last step down}\label{sub : downstep}
We now justify that universality still holds for paths with a last step down. This follows from
the correspondence established in Subsection \ref{sub: corresp}. We will not be able to
characterize the type of the typical paths. Nevertheless, by the correspondence, we only have to
consider the paths ${\mathcal P} \in
 \mathcal P _{m,l}$ ($l>0$) with
 last step down for which the associated path ${\mathcal P}'$ are typical.\\
By Subsection 4.2, to any path $\mathcal P\in  \mathcal P _{m,l}$ with a last step down and a
marked (resp. unmarked) origin, we can associate a path $\mathcal P' \in  \mathcal P _{m,l}$ with
the same weight. In each case, the number of trajectories and preimages of paths $\mathcal P'$ is at
most $\sum_{m=0}^{s_N-1}\sum_{p=0}^m T_{m-p,
  l+2p}=\sum_{m=0}^{s_N-1}C_{2s_N}^{m}.$ So, the contribution of such paths,
which we note $Z_4$, is then (with computations as in the last subsection) at most
\begin{eqnarray}
Z_4 & \leq & 2 \sum_{m=0}^{s_N-1}  \, C_{2s_N}^{m}\theta^{l}\sigma^{2m}
\exp{\{-\frac{(l+m)^2}{2N}\}}\crcr && \sum_{M\leq s_N^{\epsilon}}\frac{1}{M!}\left (
    \frac{(l+m)^2}{2N}\right )^M\sum_{r\leq M}\frac{( 6M)^r}{r!}\sum_{q\geq 0}
    \frac{1}{q!}\left (Cs_N^{-1/2+\epsilon}\right)^q\sum_{d\geq 0}\frac{1}{d!}\left (Cs_N^{-1+\epsilon}\right)^d\crcr
    &\leq & 2 \rho_{\theta}^{2s_N}\exp{\{e^6 s_N^2/N\}} \sum_{q\geq 0} \frac{1}{q!}
    \left (Cs_N^{-1/2+\epsilon}\right)^q\sum_{d\geq 0}\frac{1}{d!}\left (Cs_N^{-1+\epsilon}\right)^d.\label{iomarquelastdown}
\end{eqnarray}
We then readily see that we can neglect paths $\mathcal P$ with edges passed at least three times
($q\geq1$) or/and loops ($d\geq 1$) and that $Z_4$ is of the order of $Z_1$. This is sufficient to
 conclude to the announced universal behavior.
\section{Higher moments of $Tr M_N ^{2s_N (+1)}$ for
  $s_N=O(\sqrt N)$, $\theta > \sigma$}
In this section, we show that, for $s_N=O(\sqrt N)$ and in the large limit $N$, universality still
holds for any moment of the traces $Tr M_N ^{2s_N (+1)}$. As yet mentioned in Section 2, such a
result is not required for the proof of Theorem \ref{Uni1Pe} in the case where $\theta>\sigma$. Yet
it will be needed in the case where $\theta\leq \sigma$ and the analysis is similar in both cases.
We have the following result. \\
Let $M_N^G$ be of the deformed GUE (resp. GOE) of parameter $\theta> \sigma$ and $A>0, \epsilon>0$
be fixed.
\bp \label{Prop: variance } Let $M_N$ be a Deformed Wigner matrix of type of type $(i)-(iv)$
(resp. $(i')-(iv')$) and of parameter $\theta >\sigma.$ For any $k \geq 1$ and for any
$\epsilon<c_i<A, i=1, \ldots, k$, there exists a constant $C'(A,\epsilon, k)>0$ such that, for $N$
large enough,
\be \mathbb{E}\prod _{i=1}^k \Big [Tr (\frac{M_N}{\rho_{\theta}}) ^{[c_i \sqrt
N]}\Big]\leq {C'(A,\epsilon, k)} \, \text{ and } \, \mathbb{E}\prod _{i=1}^k \Big [Tr
(\frac{M_N}{\rho_{\theta}}) ^{[c_i \sqrt N]}\Big] = \mathbb{E}\prod _{i=1}^k \Big [Tr
(\frac{M_N^G}{\rho_{\theta}})^{[c_i \sqrt N]}\Big]\, (1+o(1)). \label{expdiff}\ee \ep
In the complex setting, universality of the variance (combining Theorem $\ref{ThPe}$ and some
considerations of Section 2) leads to $\mathbb{E}\left[ Tr M_N^{[c \sqrt N]}\right
]^2=\rho_{\theta}^{2[c \sqrt N]}\exp{ \{2(\frac{c \sigma_{\theta}}{\rho_{\theta}})^2\}}(1+o(1))$
for any $c>0$.\\

For the proof of Proposition \ref{Prop: variance }, we first focus on the variance of $Tr
M_N^{L_N}$, in the case where $L_N/ \sqrt N \to c \geq 0$. Then, we indicate the modifications
needed to consider higher moments. We use a method once more inspired by \cite{Si-So1} (and also
used in \cite{Si-So2} and \cite{Sos}). It calls on a natural extension of the \textit{construction
procedure} introduced in Section 3 of \cite{Si-So1}.
\subsection{The variance.\label{subsec: variance}}
The variance of $TrM_N^L$, for arbitrary $L \geq 1$, can be written, noting $M_N=M$ for short,
\be \label{variancetrace}
\sum_{\mathcal P _1= \{i_o,\ldots , i_{L-1}, i_o\} }\sum_{\mathcal
  P _2= \{i'_o,\ldots , i'_{L-1}, i'_o\} } \mathbb{E}\left (\prod_{j=1}^L M_{i_{j-1}i_j}\prod_{j'=1}^L M_{i'_{j'-1}i'_{j'}}\right )\!-\!\! \mathbb{E}\left( \prod_{j=1}^L M_{i_{j-1}i_j}\right )\! \mathbb{E}\left (\prod_{j'=1}^L M_{i'_{j'-1}i'_{j'}}  \!\right ).
\ee
Now, by the independence assumptions on the entries of $M_N$, the only non zero terms in the above sum
  correspond to the pairs of paths ${\mathcal P}_1$ and ${\mathcal P}_2$ having
  a common edge.
\bdefi We call correlated any pair of paths with a common edge. It is simply correlated if each
edge appears at most twice in the union of the two paths. \edefi
Assume given two correlated paths ${\mathcal P}_1$ of ${\mathcal P_{m_1,l_1}}$ and
  ${\mathcal P}_2$ of ${\mathcal P_{m_2,l_2}}$ (with $2m_i+l_i=L$ for $i=1,2$) such that the first
  common edge they share is $(vw)$. We then glue these two paths,
  erasing the edge $(vw)$, by the so-called \emph{construction
  procedure} as follows. We first read ${\mathcal P}_1$ until meeting the
  left endpoint $v$ then jump to ${\mathcal P}_2$ and follow
  ${\mathcal P}_2$ (in the reverse orientation of ${\mathcal P}_2$ if
  the edge has the same orientation in the two paths). Then, when
  meeting $w$ after $L-1$ steps, we jump back to ${\mathcal P}_1$. In this way, we obtain a path
  ${\mathcal P}$ denoted by $\mathcal P= \mathcal P_1 \vee
  \mathcal P_2$ with $2L-2$ steps and a certain amount of up steps. \\
  We now explain how to invert the procedure.
  Assume that $\mathcal P$ belongs $\mathcal P_{m,l}$ where $l+2m=2L-2$. We shall describe the structure of preimages of $\mathcal
P$, in order to obtain a bound for their number. Call $x$ its underlying trajectory. Let $\tau$ be the instant of the left endpoint of first common edge shared by the paths glued. Then either $x(\tau)=0$ or $x(\tau)>0.$  In the latter case, by the definition of
$\tau$ , there exists an interval of time $I$ containing at least $[\tau, \tau+L-1]$ such that
$x(t)\geq x(\tau), \forall t\in I.$ This trivially holds if $x(\tau)=0$. Then, once the instant $\tau$ is
chosen, the ``end'' of $\mathcal P_2$ occurs $L-1$ steps after the instant $\tau$. Now there remains to
fix its origin and orientation: there are at most $2L$ ways of doing so. From this we deduce that
the number of preimages of a path $\mathcal P$ of length $L_o$ (here $L_o=2L-2$) is at most
$$2L K_N(x), \text{ where }K_N(x)= \sum_{\tau\leq L_o-L+1}\: \prod_{s \in [\tau, \tau+L-1]}1_{x(s)\geq x(\tau)}.$$
It is enough for our purpose to show, if $\mathbb{E}_{m,l}$ denotes the uniform distribution on the set
$\mathcal T_{m,l}$, that \be \label{espkn}\mathbb{E}_{m,l}(K_N(x))\leq C'l+ C \sqrt{L_o}, \ee
 for some constants $C, C'.$
 If $l=0$, such a result has been established in \cite{Si-So1}. To prove (\ref{espkn}) in the general case, and loosely speaking, we show that
 $K_N(x)$ does not grow faster than
 $ K_N(y)+lC'$, for some Dyck path $y,$ introduced in the following.
\paragraph{}To this aim, given $l$ and $m$, we associate bijectively to any trajectory
$x$ of $\mathcal T_{m,l}$ a sequence of $p'+1$ sub-Dyck paths as follows. Let $t_i$ be the last
instant where $x(t)=0.$ In between $[0, t_i]$, the $x$ trajectory defines a sub Dyck path called
$y_0$ of length $2m_0.$ Let then $l_1$ be defined by\\
$\displaystyle{l_1=\min \{x(t): \: x(t)-x(t-1)=-1 \text{ and } t\geq t_i\}\text{ or } l_1=l\text{ if }\{t\geq t_i, x(t)-x(t-1)=-1 \}=\emptyset.}$ \\
In other words, $l_1$ is the maximum level under which the path never falls after reaching it. Let
$t_1^1$, $t_1^2$ be the instants $t_1^1=\min \{t\geq t_i, x(t)=l_1\}$ and $t_1^2=\max \{t\geq t_i,
x(t)=l_1\}.$ Then, in between $t_1^1$ and $t_1^2$, $x(t)-l_1$ defines a Dyck path $y_1$ of length
$2m_1=t_1^2-t_1^1.$ If after $t_1^2$ no down step occurs, we set $l_2=l-l_1$. Otherwise we then define $l_2=\min \{x(t): \: x(t)-x(t-1)=-1 \text{ and } t\geq t_1^2\}$
and the associated times $t_2^1$ and $t_2^2$ in the same way as before. Then the path $\mathcal P$
makes $l_1$ up steps, then follows the Dyck path $y_1$ of length $2m_1$. Afterwards it makes $l_2$
up steps and follows the second Dyck path $y_2$ and so on. Thus, we obtain a sequence $l_i\geq 1$,
$i=1,\ldots, p'$, such that $\sum_{i=1}^{p'} l_i=l$.  Now it is easy to see that each path $x$ of $\mathcal
T_{m,l}$ can be uniquely defined by a number $1\leq p'\leq l$, a sequence of $p'$ positive integers
$l_i, i=1,\ldots, p'$ such that $\sum_{i=1}^{p'}
  l_i=l$ and a sequence of $p'+1$ sub-Dyck paths $y_i, \, 0 \leq i\leq p'$ of
positive length $2m_i, i=0,\ldots ,p'$ (except $m_o, m_{p'} \geq 0$) satisfying $\sum_{i=0}^{p'}
2m_i=2m$. We then say that the trajectory $x$ is of type $(p', l_1, \ldots, l_{p'}, 2m_0, \ldots,
2m_{p'}).$ We denote $y$ the Dyck-path of length $2m$ made of the succession of the $p'+1$
trajectories $y_0, \ldots , y_{p'}.$ \\

The proof of (\ref{espkn}) is now based on the following Lemma. In the sequel, given $l\leq
L_o$, $\mathcal T_{L_o}^{l}$ is the set of paths with $L_o$ steps and ending at level $l$. And we let $T_{L_o}^{l}$ be its cardinal.
Given a trajectory $x$ in $ T_{L_o}^l$, we denote by $y_0, \ldots, y_{p'}, p'\leq l$ its sub Dyck
paths
and  by $\nu(y_i), i\leq p',$ the number of returns to $0$ of the Dyck path $y_i$. Let then $k_o$ be a given integer.%
 \bl \label{Lem: KN} Assume that $L_o\leq k_o (L-1)$ and let $x\in \mathcal T_{L_o}^{l}$.
Then, \be \label{forbase}\hspace*{-3.5cm}(i) \: K_N(x)=\sum_{\tau\leq L_o-L+1}\prod_{s\in [\tau, \tau+L-1]} 1_{x(s)\geq
x(\tau)}\leq l+\sum_{i=0}^{p'}1_{m_i\geq L}K_N(y_i)+ \sum_{i=0}^{p'}\nu(y_i),\ee 
$(ii)$ for any $k$, there
exists a polynomial $Q_k \in \mathbb{R}(X,Y)$ of total degree $k$ such that \be \mathbb{E}_{m,l}(K_N(x))^k \leq
Q_k(\sqrt L, l).\label{borne sur lesmoments} \ee
\el
\paragraph{Proof of Lemma \ref{Lem: KN} :} Point $(i)$ follows from the fact that either $\tau$ is chosen amongst the ``rises $l_i,i\leq l,$'' or
$\tau\in \cup y_i.$ Assume in the latter case that $\tau \in y_i$ and $y_i(\tau)>0.$ Necessarily $
\tau+L-1\in y_i$. Thus $\tau$ can be chosen amongst vertices of sub Dyck paths of length greater than
$L-1$. The last case is when $\tau$ is chosen amongst the
returns
to $0$ of the sub Dyck paths. This yields (\ref{forbase}).\\
For $(ii)$, it can easily be inferred from \cite{Si-So1} that there exists $C>0$ 
such that for any integer $k$, $\mathbb{E}_{m_i,0}(K_N(y_i))^k\leq (C\sqrt L)^k,$ for any $m_i\in \{L-1,\ldots, k_o(L-1)\}.$ Similarly there
exist $ C'>0$ independent of $l, L_o,p',m_i$ such that for any $k\geq 1,$
 $\mathbb{E}_{m_i,0}(\nu (y_i))^k\leq {C'}^k$. 
The above estimates then give (\ref{borne sur lesmoments}) since $p'\leq l$ and $ \#\{i \leq p': \, {m_i\geq
     L} \} \leq k_o$. $\square$\\

We now examine in more details the contribution of correlated paths ${\mathcal
  P}_1$ and ${\mathcal P}_2$ to the variance of the Trace in case $\theta >
\sigma$ and $L=L_N$ with $L_N=o(\sqrt N)$. Following Section 4, one knows that the corresponding
path ${\mathcal
  P}_1 \vee{\mathcal P}_2$ has a non-null contribution if it is without self-intersection as soon as its
origin is marked and has a unique self-intersection otherwise. Assume, for ease, that the path
${\mathcal P} \in \mathcal P _{m,l}$ (such that $2m+l=2(L_N-1)$) obtained by the gluing has no
self intersection and last step up. This implies that all the edges are passed at most twice
except the common edge. So, this common edge may appear three or four times at most. But, the
probability of such an event vanishes as $N \to \infty$. This is based on the following argument
already given in \cite{Si-So1}: it can be easily shown that the ratio of the number of simple
paths of length $2(L_N-1)$ that have an edge $(vw)$ to the whole number of simple paths of length
$2(L_N-1)$ tends to zero when $N \to \infty$. As a result, one can assume in the limit that the
common edge appears exactly one time in each of the subpaths ${\mathcal
  P}_1$ and ${\mathcal P}_2$ which give a weight of at most $\theta^l \sigma^{2m}
\frac{\sigma ^2}{N}.$ Thus, using $(\ref{espkn})$, the contribution of all the paths gluing to the
path $\mathcal P$ is at most $C{\frac{L_N^{3/2}+ L_Nl}{N} \, \theta^{l}\sigma^{2m} \sigma ^2}=
o(1) \theta^{l}\sigma^{2m} \sigma ^2$ (as $L_N=o(\sqrt N)$). Finally, in the case
 where ${\mathcal P}$
has a last step down, it is easy to check, using our fundamental correspondence of Subsection
\ref{sub: corresp}, that $K_N(x')\leq K_N(x'').$ This is enough to ensure that the above analysis
extends readily to paths with last step down.
 We do not explain
more. Consequently, in the limit, the variance $var(Tr
M_N ^{L_N})$ is universal and $var(Tr
M_N ^{L_N})=o(\rho _{\theta} ^{2L_N})$.\\
In the scale $L_N=O(\sqrt N)$, the approach is similar and relies on results of Section 5. In
particular, paths ${\mathcal P}$ with simple self-intersections (including the case where the
glued edge is a multiple edge) have to be taken into
 account. The details of the proof that the variance only depends, at the
 leading order, on the two
 first moments of the entries $M_{i,j}$ are easy and left.
This finishes the proof of Proposition \ref{Prop: variance } for the variance.

\subsection{Higher moments. }To investigate higher moments, we refer the reader to \cite{Si-So1} for the extension of the
construction procedure to glue more than two paths. We here indicate the main changes to be done
to adapt the construction procedure to our case. 
For ease of explanatory (minor modifications should be done to consider expectations as in (\ref{expdiff})), we consider expectations of the type \be \label{exphigh} \mathbb{E} \left ( Tr M_N^L
-\mathbb{E}(Tr M_N^L) \right )^{k'}, \: k'>2. \ee

It is clear that, due to the independence
assumptions on the entries of $M_N$, the expectation (\ref{exphigh}) splits into a product of expectations
over different \emph{clusters } defined as follows.%
 \bdefi A set of paths $\mathcal P_i, i=1, \ldots ,j$ of length $L$ is called a cluster if\\
 $\bullet$ for
 any pair $1\leq i<i'\leq j$ one can find a chain of paths from $\mathcal P_i$ to $\mathcal P_{i'}$ such
 that any two neighbor paths in the chain share a common edge.\\
 $\bullet$ the subset $\{ \mathcal P_i, i=1, \ldots ,j\}$ can not be enlarged with the preservation
 of the preceding condition.
 \edefi
It is now enough to consider the contribution of a cluster of correlated paths to the expectation
 (\ref{exphigh}). Let then $\mathcal C=\{ \mathcal P_i, i=1, \ldots ,k\}, k\leq k'$ be a given cluster of correlated
 paths.
We define the so-called \emph{modified construction procedure} as follows.\\
We find the first edge along $\mathcal P_1$ shared with
some other path of $\mathcal C$ and call this edge $(vw).$ Among the paths of $\mathcal C$ sharing
the edge $(vw)$, denote by $\mathcal P_{R,2}$ the path in which $(vw)$ has the largest multiplicity. If there are choices amongst such paths  we choose the one having the smallest number of edges common
with $\mathcal P_1$ (and of lower index). We then form $\mathcal P_1 \vee \mathcal P_{R,2}.$ \\
\textbf{Regular step : }If
$\mathcal P_1 \vee \mathcal P_{R,2}$ and the $k-2$ remaining paths still form a cluster, replace $\mathcal P_1$ with $\mathcal P_1 \vee \mathcal P_{R,2}$ and go on the procedure with the remaining paths to be glued. This is as in \cite{Si-So1}.\\
\textbf{Modified step : }If $\mathcal P_1 \vee \mathcal P_{R,2}$ is disconnected of the other paths, we denote by $n_1$
the number of paths of $\mathcal C$ containing $(vw)$ and the corresponding paths by $\mathcal P_1, \mathcal P_{\delta_1}, \ldots, \mathcal P_{\delta_{n_1-1}}$. Then $(vw)$ is a simple edge of these $n_1$ paths. If $n_1$ is even, we glue these $n_1$ paths
as in \cite{Si-So1}. We read $\mathcal P_1$ until meeting $v$, then switch to $\mathcal P_{R,2}$
and make $L-1$ steps until meeting $w$, then read successively the $n_1-2$ others paths and finally
read the end of $\mathcal P_1.$ We denote by $\mathcal P_1 \vee \mathcal P_{\delta_1}\vee\ldots \vee \mathcal P_{\delta_{n_1-1}}$ the path obtained. This is called the \emph{standard modified step}. We then replace $\mathcal P_1$ with $\mathcal P_1 \vee \mathcal P_{\delta_1}\ldots \vee\mathcal P_{\delta_{n_1-1}}$ and go on the procedure.
The last case, which is not
covered in \cite{Si-So1}, is when $n_1$ is odd as the edge $(vw)$ appears only once in each of the
$n_1$ paths. As $\mathcal P_1 \vee \mathcal P_{R,2}$ is disconnected of the rest of the cluster,
$\mathcal P_{R,2}$ has no other edge than $(vw)$ in common with any of the other paths of
$\mathcal C \setminus{\{\mathcal P_1\}}$. And $\mathcal P_{R,2}$ has no other edge than $(vw)$ in
common with $\mathcal P_1$, otherwise $\mathcal P_1$ (and thus $\mathcal P_1 \vee \mathcal P_{R,2}$ ) shares
also an edge distinct of $(vw)$ with some path of $\mathcal C \setminus{\{\mathcal P_1, \mathcal
P_{R,2}\}}.$ Then, we take off $\mathcal P_{R,2}$ of the cluster. We then perform  the standard modified step (or the regular step if $n_1=3$) with $\{\mathcal P_1,\mathcal P_{\delta_1}, \ldots, \mathcal P_{\delta_{n_1-1}}\}\setminus \{\mathcal{P}_{R,2}\}$ . We then define the set of ``erased paths'' as $\mathcal E_1=\{P_{R,2}\}.$ We then go on the procedure with $\mathcal P_1$ replaced by $\mathcal P_1 \vee \mathcal P_{\delta_1}\overline {\vee  {\mathcal{P}_{R,2}}} \ldots \vee \mathcal P_{\delta_{n_1-1}},$ and the remaining paths of $\mathcal C\setminus{\mathcal E_1}$ to be glued.  This is the \emph{non standard Modified Step.}\\
The above \emph{modified construction procedure} ends with a set of erased paths $\mathcal E_{n_o}=\{\mathcal P_{e,i}, i=1, \ldots, n_o\}$ for some $0\leq n_o\leq k-1$ and a path $\mathcal P'$ of length $(k-n_o)L-q$, $(k-n_o)\leq q \leq 2(k-n_o),$  which has been obtained from $\mathcal C':=\mathcal C\setminus \mathcal E_{n_o} $ by the \emph{construction procedure} described above. The erased paths $\mathcal P_{e,i}, i\leq n_o,$ are disconnected of $\mathcal P'$ and pairwise disconnected. Each erased path shares at least one ''distinguished'' simple common edge with the set of erased edges (the one leading to the non standard modifed step). If a path $\mathcal P_{e,i}$ shares more than one edge with the set of erased edges, then the multiplicity of such an edge is at most $2k$ in $\mathcal P_{e_i}.$
\paragraph{}If $n_o\geq 1$, we show that the contribution of the cluster $\mathcal C$ is at most of the order of the contribution of a set of ``clusters'' $\mathcal D=\{\mathcal C', \mathcal P_{e,i}', i\leq n_o \}$, where the $\mathcal P_{e,i}', i\leq n_o,$ are paths of length $L$ pairwise disconnected and disconnected from $\mathcal C'.$ This follows from the fact that the edges common to $\mathcal P_{e,i}, i\leq n_o$ and  $\mathcal C'$ are only even erased edges of multiplicity at most $2k.$
The number of possible choices for $n_o$ and the corresponding set of erased paths is a function of $k$ only. Then, fixing the erased edges in $\mathcal C'$ and the $n_o$ ''distinguished'' edges determines the corresponding paths $\mathcal P_{e,i}.$ 
There are at most $C_{2k}^{k}k!$ different ways to choose the even ``distinguished'' erased edges corresponding to the erased paths. The number of ways to affect some of the remaining $q-n_o$ even erased edges to some paths $\mathcal P_{e_i}$ (if any) is a function of $k$ only. As $\mathbb{E}(M_{vw})^{i+j}\leq C_k\mathbb{E}(M_{vw})^{j}\mathbb{E}(M_{vw})^{i},\forall i,j\leq 2k$, we deduce that $\mathcal C$ and $\mathcal D$ have weights of the same order.  From this we deduce that it is enough to consider the case where $n_o=0$ only.\\
To prove Proposition \ref{Prop: variance }, it is now enough to estimate the number of preimages of a path obtained by the gluing of $k$ paths of length $L=L_N=O(\sqrt N)$, using only regular and standard modified steps. 
Assume that we have fixed the order of the paths read, the moments of time we use a modified step and how many paths we combine at each modified steps. Note that the number of such choices is a function of $k$ only. Denote by $g$ (resp. $g'$) the number of regular (resp. modified) steps. 
Denote by $2q_i, i\leq g'$ the length of the cluster at each of the modified steps.
Then the number of erased edges is $q=\sum 2q_i+2g$. Let then $\mathcal P$ be a path of length $kL-q$ ending at level $l$. Then it can be shown (see \cite{Sos}, p 42) that the number of preimages $\mathcal P$ divided by $N^{q/2}$ does not grow faster, for typical paths, than 
$ \frac{C_k(2L)^{k-1} (K_N(x))^{k-1}}{N^{k-1}}$ for some constant $C_k$ depending on $k$ only.
This follows from the following fact. Given $g, g',q_i, i=1,\ldots, g'$,
the number of preimages of the path is at most 
$C_k(2L)^{g+\sum_{i=1}^{g'}2q_i-g'}K_N(x)^{g+g'}/N^{\sum_{i=1}^{g'} q_i+g}.$
Conversely each time a modified step with $2q_i$ paths is performed, the contribution of the path is decreased of a factor $(s_N/N)^{2q_i-2}$ with that of a typical path of the same length. Indeed, the instant $\tau_i$ at which one starts the modified step determines $2(q_i-1)$ vertices of the path. Replacing these vertices with pairwise distinct vertices of $\mathcal N_o\setminus\{i_o\}$
 only decreases the weight of the path of at most $O(s_N^{2(q_i-1)})$ (without changing the constraint on the possible choices for $\tau_i$) since it can easily be shown that typical clusters are such that $\mathcal P'$ has no vertices of type greater than $C \sqrt {s_N}$ for some constant $C.$
Combining the whole, and using that $s_N=O(N^{2/3})$ in the whole paper, gives that the contribution of such paths is of the order of 
$C'_k(2L^{3/2}/N)^g  N^{-\sum (q_i-1)/3} (K_N(x)/\sqrt L)^{g+g'}.$ This is enough to ensure the above result. In this way, we also deduce that clusters for which $g'\geq 1$ have a negligible contribution.
Using (\ref{borne sur lesmoments}), one can then deduce
Proposition \ref{Prop: variance }, mimicking the arguments of Subsection \ref{subsec: variance}. $\square$

\section{The case where $\theta=\sigma$}
As explained in Section 2, in both cases where $\theta=\sigma$ and $\theta<\sigma$, the results of
Theorems \ref{Uni1Pe}, \ref{UnikPe} and \ref{theo: Unireal} follow from the universality of all
the limiting moments of Traces of exponent of type $o(N^{2/3})$ and $O(N^{2/3})$. This section is
devoted to the case where $\theta =\sigma$ and we discuss on the case where $0 < \theta <\sigma$
in the last Remark $\ref{casetheta<sigma}$. Here, we only establish the universality of limiting
expectation of traces of $M_N^{L_N}$ with $L_N=O(N^{2/3})$. It is then an easy task to translate
the computations to higher moments (c.f. Section 6) as well as to the case $L_N=o(N^{2/3})$ and to
see that universality holds in this case also (the detail is left). Here is the main result of
this section.

\bt \label{Theo : theta=sigma} Let $M_N$ be a Deformed Wigner matrix of type  $(i)-(iv)$
(resp. $(i')-(iv')$) and of parameter $\theta =\sigma$. Let $(L_N)_N$ be a sequence such that
$\exists c>0, \lim_{N \to \infty}  \dfrac{L_N}{N^{2/3}}= c$. Then, there exists $C'>0$ such that,
for $N$ large enough,
$$\mathbb{E} \left (Tr M_N^{L_N}\right )\leq C'(2\sigma)^{L_N}\text{ and } \,
\mathbb{E} \Big [Tr (\frac{M_N}{2\sigma})^{L_N}\Big ]=\mathbb{E}\Big[Tr
(\frac{M_N^G}{2\sigma})^{L_N}\Big] (1+o(1)).$$ \et

\brem \label{Rem : highermoments theta=sigma}Similar results to those of Theorem 3 in \cite{Sos}
can be stated for joint moments. For instance $\mathbb{E}\left (Tr M_N^{L_N}\right
)^k=(1+o(1))\mathbb{E}\left (Tr
  (M_N^G)^{L_N}\right )^k\leq C_k(2\sigma)^{kL_N}\exp{\{C'k L_N^3/N^2\}},$ $\forall k \geq 1$. \erem

As in the preceding sections, we only consider paths of even length $L_N=2s_N$. As
$\theta=\sigma$, $\rho _{\theta}=2 \sigma$ and even paths ($l=0$) have to be taken account. More
precisely, for $l=0$ and $s_N=O(N^{2/3})$, it is proved in Section 4 of \cite{Sos} that even
typical paths have unmarked origin, self-intersections of type $3$ at most (with only a finite number of type $3$) and edges passed only twice. Their total contribution is of the order of
\be \label{contrieven}
N T_{s_N,0}\sigma^{2s_N}\sim \frac{N(2\sigma)^{2s_N}}{\sqrt{\pi}s_N^{3/2}}.
\ee
To prove Theorem \ref{Theo : theta=sigma}, we shall now examine paths such that $l>0$. We first
show that those with last step up are negligible. This is intuitively clear since the analysis is
close to that of \cite{Sos}, where typical paths have a non marked origin. Then, we call on our
fundamental correspondence to investigate paths with a last step down. As in Section 5, we are not
able to characterize the type of typical paths. Nevertheless, we show that, in the set $\{
\mathcal{P} _{m,l}, \, l+2m=2s_N \text{ and } l>0\}$, typical paths with last step down
 have edges passed at most twice and are in correspondence with paths admitting self intersections of type 3 at most. Note that, as $\theta=\sigma$, every typical path has the weight $\sigma^{2s_N}$.\\

Before entering the details of the proof, we show some important technical results.

\bp  \label{lem : majol} There exists $C_{crit}>0$ independent of $N$ such that the paths (with
marked or unmarked origin) for which $l\geq  C_{crit}\sqrt N$ give a negligible contribution
w.r.t. $(2\sigma)^{2s_N}$ (and (\ref{contrieven})). \ep
This proposition is based on the following lemma.
\bl \label{lem : majoexpTml} There exists $C_o>0,$ independent of $l$ and $N$, such that
$\displaystyle{T_{m,l}\leq (l+1)e^{\{-C_o \frac{l^2}{s_N}\}}T_{s_N,0}.}$ \el
\paragraph{Proof of Lemma \ref{lem : majoexpTml}: }
Writing $s_N=m+q$ where $q=s_N-m=l/2$, we have that
\begin{eqnarray}
&&\dfrac{C_{2s_N}^{l+m}}{C_{2s_N}^{s_N}} =\exp{\{-\sum_{j=m+1}^{m+q}\log (1+\frac{q}{j})\}} \leq
\exp{\{-\frac{l}{2}\log (1+\frac{l}{2s_N})\}}\leq \exp{\{-C_o\frac{l^2}{s_N}\}}.\nonumber
\end{eqnarray}
Thus $\displaystyle{\dfrac{T_{m,l}}{T_{s_N,0}}=
  \frac{(l+1)(m+1)}{l+m+1}\dfrac{C_{2s_N}^{l+m}}{C_{2s_N}^{s_N}}\leq
  (l+1)\exp{\{-C_o \frac{l^2}{s_N}\}}.} \, \square$

\paragraph{Proof of Proposition \ref{lem : majol}:} Here we do not assume that the origin is marked. The contribution of paths with $l>0$ is then at most, from Remark \ref{rem: prem avril} and mimicking (\ref{majofact}) to (\ref{majofactM2}),
\begin{eqnarray} \label{sommesurl}
&&T_{m,l} N\: \sigma^{2s_N} \, \exp \{\sum_{k=2}^{10} \frac{(C \, (l+m))^k}{N^{k-1}} \} \,
\sum_{M_1=\sum_{i>10}N_i}\frac{1}{M_1!}\left (\frac{C}{N^{2/9}}\right)^{M_1}\cr &&\leq N T_{m,l}\:
\sigma^{2s_N}\:  \exp{\{C_{1} \frac{s_N^2}{N}\}}\leq N T_{m,l}\: \sigma^{2s_N}\: \exp{\{C_{2}
N^{1/3}\}},
\end{eqnarray}
for some constant $C_1, C_2>0$ independent of $N$. Here an extra $N$ has been added to include the
case where the last step is down. By Lemma \ref{lem : majoexpTml}, 
it is easy to see that there exists a
constant $C_{crit}>0$ such that, in the large $N$ limit, the subsum in (\ref{sommesurl}) over
$l\geq C_{crit} \sqrt N$  gives a contribution which is negligible with respect to
$(\ref{contrieven})$ (i.e. to that of even paths). $\square$

\brem \label{remls_N=o(N^{2/3})} For $s_N=o(N^{2/3})$, it is enough to consider the paths with
$l\leq C_{crit}s_N^{3/2}N^{-1/2}=o(N^{1/2})$ for some constant $C_{crit}.$ \erem

The arguments used to prove Theorem \ref{Theo : theta=sigma} are quite similar to those yielding
Lemmas 5,6 and 7 (Section 4) in \cite{Sos}. In particular, we show that typical paths may contain
some non-closed vertices. But, using Remark \ref{Rem: nbclsoed}, we need to have some control on
the number of possibilities of choosing a non-closed vertex. In Section \ref{Sec : fluctuations},
given a path $\mathcal P$ in some $\mathcal{P} _{m,l}$ (with $l+2m=2s_N$) and denoting by $x\in
\mathcal{T} _{m,l}$ the associated trajectory, the number of possible choices was estimated from
above  by ${\rm{\max}} _t x(t)\leq l+m$. Here, in the scale $O(N^{2/3})$, a better estimate is
required.
We now show that the number of choices of a non-closed vertex can be estimated through a quantity
as $\max_t y(t)$, which involves the particular Dyck path $y \in \mathcal{T} _{m,0}$ introduced in
Subsection \ref{subsec: variance}. Let then $z(t)$ denote the maximal number of ways to choose a non closed vertex at the
instant $t$. We now prove that \be \label{estz}\sum_{t=1}^{l+m}z(t)\leq \frac{l^2}{2}+ m \times
\max _t y(t) +lm.\ee Recall that $z(t)$ can be bounded from above by the number of marked vertices
opened but not closed before $t$. If at the instant $t$, one crosses the subpath $y_i$, then there
are at most $z(t) \leq \sum_{j \leq i} l_i+ y_i(t) $ such vertices. It may also happen that the
instant $t$ corresponds to, for example, the rise $``l_i''$ and then $z(t) \leq \sum_{j \leq i}
l_i$ (note also that in any case $\max z(t)\leq l+\max y(t)$). Therefore, by a straightforward
computation, one has $\sum_{t=1}^{l+m}z(t) \leq \frac{l^2}{2}+\sum_{i=1}^{p'}
\sum_{t=1}^{m_i} y_i(t) +\sum_{i=1}^{p'} l_i \sum_{j\geq i}m_j.$
This yields (\ref{estz}).  \\
The estimated needed on the quantity $\max_{t\leq 2m} y(t)$, with $y$ as before, is given by Lemma
\ref{Lem : max y} stated below. Assuming this lemma, we are in position to establish Theorem
\ref{Theo : theta=sigma}: we essentially mimic the ideas in Section 4 of \cite{Sos}, distinguishing
paths with a last step up or down.

\paragraph{Proof of Theorem \ref{Theo : theta=sigma}:}
From now on, one considers paths with $l\leq C_{crit}\sqrt N,$ where $C_{crit}$ is the constant
given by Proposition \ref{lem : majol}.

\textit{$1^{rst}$ case: Paths with last step up}. It can be easily inferred from the computations
of the preceding Sections and of \cite{Sos} (p. 41) that, up to a negligible error, for each $l$,
one can assume that
\begin{itemize}
\item[{($A_1$)}] the number of self intersections is smaller than $B\frac{s_N^2}{N},$ for some fixed $B>0$ (large enough).
\item[{($A_2$)}]  the maximal type of a vertex of $\mathcal P$, denoted $\nu'_N(\mathcal P)$, satisfies $\nu'_N(\mathcal P)\leq CN^{1/3}/\ln N$.
\item[{($A_3$)}] $i_o$ is of type one.
\end{itemize}
Assumptions $(A_1)$ is straightforward from Section 5 and \cite{Sos}. Assumption $(A_2)$ can be
proved using formula (\ref{majoMavec exp}). Assumption $(A_3)$ follows from the fact that the
contribution of paths where $i_o$ is of type at least $2$ is of the order $M/s_N=O(N^{-1/3})$ that
of paths for which $i_o$ is of type one. \brem Actually, we can prove that $\nu'_N(\mathcal P)\leq
4$ for paths giving a non negligible contribution to the expectation of the Trace. Yet this
estimate, needed to consider edges passed at least three times, requires some technical tools we
only develop in the sequel. \erem

We can now proceed to the estimation of the contribution of paths with last step up, under
assumptions $(A_1)$ to $(A_3)$. Consider paths of type $(N_o, N_1, \ldots, N_{l+m})$. Let $r$
denote the number of non closed vertices of type $2$ and $q$ be the number of vertices of ${\cal
N}_2$ for which the edge is read at least three times for any such path.  Then, given $(N_o, N_1,
\ldots, N_{l+m}), r\leq N_2$ and $q\leq N_2$,  it is easy to see that Proposition $\ref{Prop:
majoWmE}$ reads for such paths as
\begin{eqnarray*}
\Omega _m \mathbb{E}_{max} \leq \frac{1}{N^{s_N}} \: 3^r\: \prod_{k=3}^{10}\left (Ck\right
)^{kN_k}\prod_{k=11}^{s_N} \left (Ck\right )^{4kN_k/3} \left (\frac{\theta}{\sqrt N}\right )^l
\sigma^{2m-q}C^q.
\end{eqnarray*}
Denote by $u_1<u_2<\cdots<u_r$ the ranks of the instants of self intersections of type $2$ where
the non closed vertices are chosen. Let also $v_1<v_2<\cdots<v_q$ be the marked instants of ${\cal
N}_2$ associated to an edge read at least three times. We also set $\nu_N'':=\nu_N(\mathcal
P)+\nu'_N(\mathcal P).$ The contribution of paths of type $(N_o, N_1, \ldots, N_{l+m})$ with last step up is then at most
\begin{eqnarray}
&&  e^{\{-\frac{(l+m)^2}{2N}\}} e^{\frac{Bs_N^3}{N^2}}
  \frac{\sigma^{2s_N}}{N^{N_2}} \times \sum_{j_1<j_2<\cdots<j_{N_2}}\:\sum_{u_1<u_2<\cdots <u_r}\:\sum_{v_1<v_2<\cdots <v_q}\: \prod_{i=1}^{N_2}(j_i-i)\crcr
  &&
 \prod_{k=1}^r \frac{3z(j_{u_k})}{(j_{u_k}-u_k)}\prod_{s=1}^q \frac{\nu''_N}{{(j_{v_s}-v_s)}}  \prod_{k= 3}^{10}\frac{1}{N_k!}\left (\frac{C(l+m)^k}{N^{k-1}}\right)^{N_k}\prod_{k= 11}^{l+m}\frac{1}{N_k!}\left (\frac{C(l+m)^{4k/3}}{N^{k-1}}\right)^{N_k}\crcr
&&\leq e^{\{-\frac{(l+m)^2}{2N}+\frac{Bs_N^3}{N^2}\}} \sigma^{2s_N} \, \frac{1}{(N_2-r-q)!r!q!}
\left (\frac{(l+m)^2}{2N}\right )^{N_2-r-q}\left (\sum_{t=1}^{l+m}\frac{3z(t)}{N}\right )^r\left
(\frac{Cs_N \nu''_N}{N}\right)^q\cr && \prod_{k= 3}^{10}\frac{1}{N_k!}\left
(\frac{C(l+m)^k}{N^{k-1}}\right)^{N_k}\prod_{k= 11}^{l+m}\frac{1}{N_k!}\left
(\frac{C(l+m)^{4k/3}}{N^{k-1}}\right)^{N_k}\cr &&\leq
\exp{\{-\frac{(l+m)^2}{2N}+\frac{Bs_N^3}{N^2}\}} \, \sigma^{2s_N} \, \frac{1}{(N_2-r-q)!r!q!} \left
  (\frac{(l+m)^2}{2N}\right )^{N_2-r-q} \left (\frac{Cs_N
  \nu''_N}{N}\right)^q \cr
&&\left (\frac{\frac{3}{2}l^2+ 3m \max y(t)+ 3ml}{N}\right )^r \prod_{k=
3}^{10}\frac{1}{N_k!}\left (\frac{C(l+m)^k}{N^{k-1}}\right)^{N_k}\prod_{k=
11}^{l+m}\frac{1}{N_k!}\left (\frac{C(l+m)^{4k/3}}{N^{k-1}}\right)^{N_k}. \hspace{0.9cm}
\label{formulesansvn}
\end{eqnarray}

Let $Z_3(l)$ be the contribution of all the paths of $\mathcal P _{m,l}$ with last up step without
self-intersections of type greater than $4$ (included),  for which $q=0$ and with edges read at
most twice. \bl \label{lem: z3N23}One has that $\displaystyle{Z_3:=\sum_{l=2}^{C_{crit}\sqrt
N}Z_3(l)=O(\frac{ (2\sigma)^{2s_N}}{N^{1/3}}).}$ \el
\paragraph{Proof of Lemma \ref{lem: z3N23}:} We first need a few definitions.
Given $m\geq 1$, $1\leq n\leq m$, and a sequence of integers $s_i>0, i=1,\ldots, n,$  such that $\sum s_i=m$, we consider the set of trajectories of $\mathcal T _{m,0}$ made of the succession of $n$ sub-Dyck
paths of length $2s_i$, $i=1,\dots n$. We then say that such trajectories are of class $T=(n, s_1,\ldots, s_n)$ and  call $\mathbb{T}(m)$ the set of such classes. One should note that the type $(p',l_1,\ldots, l_{p'},m_o,\ldots, m_{p'-1},m_{p'}=0)$ of a trajectory $x \in \mathcal T_{m,l}$ with last step up naturally defines the class of the associated Dyck path $y$. We denote $T_{in}$ the induced class.
 Then, denoting by $\mathbb{E}_{Y,T}$ the expectation
with respect to the uniform distribution on the trajectories of $\mathcal T _{m,0}$ of class $T$ and using (\ref{formulesansvn}), one has that
\begin{eqnarray} &Z_3(l)&\leq T_{m,l-1}  \exp{\{\frac{3l^2/2+3ml }{N}+\frac{C s_N^3}{N^2}\}} \sigma^{2s_N}\cr
&&\times \frac{1}{T_{m,l-1}}\sum_{(p',l_1,\ldots, l_{p'},m_o,\ldots, m_{p'-1}) }\:\mathbb{E}_{Y,T_{in}} \left (\exp{\{\frac{3m^{3/2}}{N} \frac{\max
    y(t)}{\sqrt m}\}}\right)\: \:\prod_{0\leq i\leq p'-1/ m_i>0} T_{m_i,0}\nonumber\\
    &&\leq T_{m,l-1}  \exp{\{\frac{3ml }{N}+\frac{Cs_N^3}{N^2}\}}\sigma^{2s_N} \max_{T \in \mathbb{T}(m)}\mathbb{E}_{Y,T} \left (\exp{\{\frac{3m^{3/2}}{N} \frac{\max
    y(t)}{\sqrt m}\}}\right)
    \label{formuletheta=sigma}
\end{eqnarray}
since $\sum_{(p',l_1,\ldots, l_{p'},m_o, m_1,\ldots ,
    m_{p'-1})}\:\prod_{0\leq i\leq p'-1/ m_i>0} T_{m_i,0}=T_{m,l-1},$
and $3l^2\leq 3C_{crit}^2 N.$
In the following, $C,\tilde C,$ $ C_1, \ldots,
C_6$ denote some constants independent of $N$.
In Lemma \ref{Lem : max y}, whose proof is postponed to the end of
this section, we show that, given a constant $C>0$, there exists $\tilde C>0$ such that,
$\max_{T \in \mathbb{T}(m)}\mathbb{E}_{Y,T}\left [\exp{\{\frac{C\max y(t)}{\sqrt m}\}}\right ]\leq \tilde C,$ for any $m$.\\
Assuming this holds we finish the proof of Lemma \ref{lem: z3N23}. It is now enough to show that
\be \exists C_1>0, \: \sum_{l\leq C_{crit}\sqrt N}N^{1/3}T_{m,l-1}\exp{\{\frac{3ml}{N}\}} \leq
C_1N T_{s_N,0}.\label{lastetim}\ee First, by Lemma $\ref{lem : majoexpTml}$, one can find $K_o>0$
such that \be N^{1/3}\sum_{K_o N^{1/3}\leq l\leq C_{crit}\sqrt N
}T_{m,l-1}\exp{\{\frac{3ml}{N}\}}\leq C' NT_{s_N,0},\label{estimqueue} \ee for some constant $C'$.
Indeed, as $m \leq C_2 N^{2/3}$, one has that
\begin{eqnarray}
&\displaystyle{\sum_{l=K_o N^{1/3}}^{C_{crit}\sqrt N }}N^{1/3}T_{m,l-1}\exp{\{\frac{3ml}{N}\}}
&\leq C_2 NT_{s_N,0}\sum_{l\geq K_o N^{1/3} }\frac{
l}{N^{2/3}}\exp{\{l(\frac{3m}{N}-\frac{C_ol}{s_N})\}}\cr
&&\leq  \frac{C_{3}}{N^{1/3}}NT_{s_N,0} \sum_{l\geq K_o N^{1/3} }\frac{l^2}{4s_N}\exp{\{-\frac{C_ol^2}{2s_N}\}}\label{remplace mlNparl2n}\\
&&\leq NT_{s_N,0} \frac{C_4}{N^{1/3}}\sum_{l\geq K_o N^{1/3}}\exp{\{-l\frac{C_oK_o
N^{1/3}}{4s_N}\}}\label{sum k_o}
\end{eqnarray}
where in (\ref{remplace mlNparl2n}), we have chosen $K_o\geq \max\{\frac{12 s_N^2}{C_o N^{4/3}},
\frac{4s_N}{C_o N^{2/3}}\}.$ As $\sum_{l\geq K_o N^{1/3}}\exp{\{-l\frac{C_oK_o N^{1/3}}{4s}\}}
=O(N^{1/3}),$ we obtain that $\displaystyle{(\ref{sum k_o})\leq C_5 NT_{s,0}}$. This yields
$(\ref{estimqueue})$. It is also straightforward that
 \be
\label{estimqueue2}N^{1/3}\sum_{l\leq K_o N^{1/3}}T_{m,l-1}\exp{\{\frac{3ml}{N}\}} \leq C_62^{2s_N},
\ee 
since $\displaystyle{\sum_{2\leq l\leq K_o N^{1/3}}T_{m,l-1}\leq
  C_{2s_N}^{s_N+1}\sim  s_N T_{s_N, 0}.}$ Then, (\ref{estimqueue}) and (\ref{estimqueue2}) yield (\ref{lastetim}) and Lemma \ref{lem: z3N23}. $\square$
\paragraph{}We then denote by $Z_4(l)$ the subsum of (\ref{formulesansvn}) over paths for which
$\nu_N(\mathcal P)\leq s_N^{1/2-\epsilon},$ for some (small) $\epsilon>0$. We then show that for
such paths and any $l$, typical paths have edges passed at most twice. Note that in this case
$\nu''_N(\mathcal P)\leq 2CN^{1/3}/\ln N.$ Thus, it is not hard to see that the summation of
(\ref{formulesansvn}) over paths for which $\nu_N(\mathcal P) \leq s_N^{1/2-\epsilon}$, and
$\sum_{i\geq 4} N_i\geq 1$ or $N_3\geq B s_N^{3/2}/N$ , whatever $q$ is, is $o(Z_3(l))$ in the
large $N$ limit. Finally, assuming that there are no self intersection of type strictly greater
than $3$, it is also easy to see that the contribution of paths with  $q\geq 1$  gives a
contribution of the order $Z_3(l) s_N\nu''_N/N=o(Z_3(l)).$
Assuming then that $q+\sum_{k\geq 4}k N_k=0$, we can then proceed as above to show that no vertex of type $3$ is in an edge read at least three times and that there are no loops.\\

Let finally $Z_5(l)$ denote the contribution in (\ref{formulesansvn}) of paths for which
$\nu_N\geq s_N^{1/2-\epsilon},$ where $\epsilon<1/32.$   \bl \label{lem: z5juste}One has
$Z_5:=\sum_{l=2}^{C_{crit}\sqrt N}Z_5(l)=o(Z_3)$.\el
\paragraph{Proof of Lemma \ref{lem: z5juste} :}We first need to introduce a few notations.
Given a path of type $(N_o, N_1, \ldots, N_{l+m})$ with $r$ non closed vertices, we set
$K:=r+\sum_{i\geq 3}iN_i+1$, $K_o:= r+\sum_{i= 3}^{200}N_i$ and $K':= 200K_o +\sum_{i\geq
200}iN_i.$
Note that $K \leq K'+1.$ The choice of the constant $200$ is not optimal here but is enough for our next computations. Let also $\Gamma_j$ be the event \\
$\displaystyle{\Gamma_{j}:= \: \{ \exists\: 1\leq t_1<t_2<\cdots <t_{j}\leq
  2s_N: \: x(t_i)=x(t_1),\: \forall i\leq j\text{ and } x(t)\geq x(t_1), \, \forall t \in [t_1, t_j]\}.}$ Then \\$\Gamma_{j+1}\subset \Gamma_{j},$ $\forall \: j\geq 2.$ As explained in Remark \ref{rem: nun}, given $K'\geq 1$, if a vertex is the left endpoint of $\nu_N$ up edges in a path, then the associated trajectory $x$ necessarily belongs to $\Gamma_{\frac{\nu_N}{K'+1}}.$ Note also that these returns are necessarily made inside a sub-Dyck path $y_i, i\leq p'$ of the trajectory. It is then an easy fact that there exists $C_1>0$ such that, for any $j$, $l\geq 0$ and for any class $T \in \mathbb{T}(m)$,
\be \label{majonu} \mathbb{P}_{Y,T}\left (\Gamma_{j}\right) \leq 4s_N^2 \exp\{-C_1j\}. \ee Thus,
using (\ref{formulesansvn}) and arguments of \cite{Sos} (p. 41), there exists $C>0$ such that
\begin{eqnarray}
&Z_5(l)\leq &\sigma^{2s_N}T_{m,l-1}\exp{\{\frac{Cs_N^3}{N^2}\}}\sum_{N_i,i\geq
200}\prod\frac{1}{N_i!}\left ( \frac{(Cs_N)^{4i/3}}{N^{i-1}}\right)^{N_i}\sum_{\nu_N \geq
s_N^{1/2-\epsilon}}  \exp{\{\frac{Cs_N \nu_N}{N}\}}\crcr &&\!\!\! \times \max_{T\in \mathbb{T}(m)}
\mathbb{E}_{Y,T}\Big ( \sum_{K_o}\frac{1}{K_o!}\left ( \frac{3m \max y(t)+ 3ml+ C
s_N^{3/2}}{N}\right)^{K_o} 1_{\Gamma_{\frac{\nu_N}{K'+1}}}\Big ).\label{calculZ5l}
\end{eqnarray}

Let $Z'_5(l)$ be the subsum over $K_o$ and $ N_i, i\geq 200,$ such that $K'\leq
s_N^{1/2-2\epsilon}$. Then
\begin{eqnarray*}
&Z'_5(l)\leq &\sigma^{2s_N}T_{m,l-1}\exp{\{\frac{Cs_N^3}{N^2}\}}
\left(\max_{T \in \mathbb{T}(m)}\mathbb{E}_{Y,T}\Big( \exp{\{2.\frac{3m \max y(t)+ 3ml+ C
s_N^{3/2}}{N}\}}\Big)\right)^{1/2}\crcr &&\sum_{\nu_N \geq s_N^{1/2-\epsilon}} \exp{\{\frac{Cs_N
\nu_N}{N}\}}s_N^{1/2-2\epsilon}\left (\max_{T \in \mathbb{T}(m)}\mathbb{P}_{
Y,T}\Big(\Gamma_{\frac{\nu_N}{s_N^{1/2-2\epsilon}+1}}\Big)\right)^{1/2}.
\end{eqnarray*}
One readily deduces from (\ref{majonu}) that $Z'_5(l)=O(e^{s_N^{\epsilon}/4}s_N^{3})= o(Z_3(l))$
and thus $Z'_5=o(Z_3).$

Let then $Z''_5(l)$ be the subsum over $K'\geq s_N^{1/2-2\epsilon}$. From Lemma \ref{Lem : max y}
proven below, the proportion of paths for which $\max_{t\leq 2m} y(t)\geq AN^{1/6}s_N^{1/2}$
decreases as $C \exp{\{-c_oA^2N^{1/3}\}}.$ Choosing $A$ large enough then ensures that the
contribution of the sole paths for which $\max_{t\leq 2m} y(t)\leq AN^{1/6}s_N^{1/2}$ has to be
taken into account. We now restrict to such paths in $Z''_5(l).$ This implies in particular that
there exists $C_2>0$ such that $\displaystyle{\frac{3m \max y(t)+ 3ml+ C s_N^{3/2}}{N}\leq
C_2N^{1/6}}.$

Using (\ref{majonu}), we can find $K_1$ such that
$$\sum_{ s_N^{1/2-\epsilon}\leq \nu_N} \exp{\{\frac{Cs_N \nu_N}{N}\}}K_1N^{1/3}\max_{T \in \mathbb{T}(m)}\Big [\mathbb{P}_{ Y,T}\big (\Gamma_{\nu_N/(1+K_1N^{1/3})}\big)\Big]^{1/2}\leq
C_3 s_N^{3}.$$

Denote then $Z''_{5,1}(l)$ the subsum over $s_N^{1/2-2\epsilon}\leq K' \leq K_1 N^{1/3}.$ Then,
either $K_o\geq s_N^{1/2-2\epsilon}/400$ or $\sum_{i\geq 200}iN_i\geq s_N^{1/2-2\epsilon}/2$. In
the latter case, $\sum_{i\geq 200} N_i\geq 1.$ Thus, there exists $C'>0$ such that
\begin{eqnarray}
&Z''_{5,1}(l)\leq &C'\sigma^{2s_N}T_{m,l-1}\exp{\{\frac{Cs_N^3}{N^2}\}}C_3 s_N^{3}\crcr &&\times
\left ( \sum_{K_o\geq \frac{s_N^{1/2-2\epsilon}}{400}}\frac{1}{K_o!}\left
(C_2N^{1/6}\right)^{K_o}+ \sum_{N_i,i\geq 200, \sum N_i\geq 1}\prod_{i\geq 200}\frac{1}{N_i!}\left
( \frac{(Cs_N)^{4i/3}}{N^{i-1}}\right)^{N_i}\right)\crcr &&\leq
C''\sigma^{2s_N}T_{m,l-1}\exp{\{\frac{Cs_N^3}{N^2}\}}C_3 s_N^{3}\times  N^{-20}.
\end{eqnarray}
In the last line, we have used Stirling's formula and the multinomial identity. This yields that
$Z''_{5,1}(l)=o(Z_3(l)).$

Denote finally by $Z''_{5,2}(l)$ the subsum over $K' \geq K_1 N^{1/3}.$ Then there exists $C>0$
such that
\begin{eqnarray}
&Z''_{5,2}(l)\leq &\sigma^{2s_N}T_{m,l-1}\exp{\{\frac{Cs_N^3}{N^2}\}}\exp{\{C N^{1/3}\}}\crcr
&&\!\!\!\times \left ( \sum_{K_o\geq  \frac{K_1 N^{1/3}}{400}}\frac{\left
(C_2N^{1/6}\right)^{K_o}}{K_o!}+ \sum_{N_i,i\geq 200, \sum iN_i\geq \frac{K_1
N^{1/3}}{2}}\prod_{i\geq 200}\frac{1}{N_i!}\left (
\frac{(Cs_N)^{4i/3}}{N^{i-1}}\right)^{N_i}\right)\label{z''52}
\end{eqnarray}
To consider the above sum, we introduce $N_{max}$ which is the maximal type of a vertex in the
path (in particular $N_{max}\leq \nu'_N$) and $M_o:=\sum_{i\geq 200}iN_i$. As $\sum_{i\geq 200}
N_i\leq \frac{M_o}{200}$ and $\sum_{N_i, i\leq N_{max}, \sum iN_i=M_o}\leq
C_{M_o+N_{max}-1}^{N_{max}-1},$ one has
\begin{eqnarray} \label{sommeMo}
&\displaystyle{\sum_{N_i, \sum iN_i\geq \frac{K_1 N^{\frac{1}{3}}}{2}}\prod_{i\geq
200}\frac{1}{N_i!}\left ( \frac{(Cs_N)^{4i/3}}{N^{i-1}}\right)^{N_i}} &\leq\sum_{M_o\geq \frac{K_1
N^{1/3}}{2}}{C'}^{M_o}N^{-\frac{M_o}{18}}\sum_{N_{max}\leq M_o} C_{M_o+N_{max}-1}^{N_{max}-1} \crcr
&&\leq \sum_{M_o\geq \frac{K_1 N^{1/3}}{2}}{C'}^{M_o}N^{-\frac{M_o}{18}}2^{2M_o}\leq 2\left (
\frac{4C'}{N^{1/18}}\right)^{\frac{K_1 N^{1/3}}{2}}
\end{eqnarray}
Inserting this in (\ref{z''52}) then yields that $Z''_{5,2}(l)=o(Z_3(l)).$ This finishes the proof
that $Z_5<<Z_3$ and finishes the proof of Lemma \ref{lem: z5juste}.$\square$ \brem In the case
$s_N=o(N^{2/3})$, one can see that the paths with non closed
vertices can be neglected since $l<<\sqrt N$ (c.f. Remark \ref{remls_N=o(N^{2/3})}).
\erem

\textit{$2^{nd}$ case: Paths with last step down.} Here we show that the contribution of paths
with last step down is of the
order of $(2\sigma)^{2s_N}$ and that (contrary to paths with last step up) they contribute in a non negligible way to the expectation of the Trace. \\
We consider the case of a path $\mathcal P$ of $\mathcal P _{m,l}$ with an unmarked origin $i_o$.
The case where the origin is marked leads to similar computations (actually, its suffices to
consider, in all the following computations, the summation over $p \geq 0$ instead of $p \geq 1$).
Thanks to the transformation built in Subsection 4.2, given the instant $2k$ of the first odd edge
of $\mathcal P$, we can associate to $\mathcal P$ a path $\mathcal P'$ of $\mathcal P _{m,l}$ with
a marked origin and a last step up. $\mathcal P'$ has the same vertices and edges as  $\mathcal
P$. And, given the level $p$ of the first odd edge, the number of uplets $( x', k)$ is at most $T_{m-p,l+2p}$. \\
We denote as before by $r$ the number of non closed vertices of $\mathcal P'$, $z'(t)$ the number
of ways to choose a non closed vertex at the instant $t$, $q$ the number of vertices of ${\cal
N}_2$ associated to an edge passed more than three times, and $M_1$ the number of vertices of
$\mathcal P'$ of type greater than 11. We still note $\mathbb{E}_{m-p,l+2p}$ the expectation with respect
to the uniform distribution on ${\cal T}_{m-p, l+2p}$. The contribution of paths $\mathcal P$ with
unmarked origin is then at most (using the computations of the previous case)
\begin{eqnarray}
&&\sum_{l\leq C_{crit}\sqrt N} \sum_{1\leq p\leq m}\,T_{m-p,l+2p}\mathbb{E}_{m-p,l+2p}\Big [\sum_{r \geq
0} \, \sum_{q \geq 1} \,\sum_{1\leq j_1,j_2,\ldots ,j_{N_2} \leq l+m}\sum_{u_1<u_2<\cdots
<u_r}\sum_{w_1<w_2<\cdots <w_q} \crcr
&&\frac{\sigma^{2(s-q)}}{N^{N_2}}\exp{\{-\frac{(l+m)^2}{2N}+\frac{Cs_N^3}{N^2}\}} \Big
((j_{1}-1)\cdots \prod_{d=1}^r\overline
  {(j_{u_d}-u_d)}\prod_{i=1}^q\overline{(j_{v_i}-v_i)}\cdots (j_{N_2}-
  N_2)\Big ) \cr
&&\prod_{1}^r (3z'(t_{j_{u_d}}))\left (c \nu''_N\right)^{q}\prod_{k= 3}^{10}\frac{1}{N_k!}\left
(\frac{C(l+m)^k}{N^{k-1}}\right)^{N_k}\frac{1}{M_1!}\left (CN^{-2/9}\right)^{M_1}\Big ]
.\label{contri1} \end{eqnarray} Note that $(\ref{contri1})\leq \sum_{l\leq C_{crit}\sqrt N}
\sum_{1\leq p\leq m}\,T_{m-p,l+2p}\sigma^{2s_N}e^{\{C N^{1/3}\}},$ for some constant $C$
independent of $p$ and $m$. Thus, as in Proposition \ref{lem : majol}, one can assume that $p \leq
C' _{crit} \sqrt N.$ Indeed, one has \be \sum_{p \geq C' _{crit} \sqrt N} T_{m-p,l+2p} \leq C_{2
s_N} ^{m-C' _{crit} \sqrt N} \leq \exp{ \{- \frac{{C' _{crit}} ^2  N}{2m} \} } C_{2 s_N} ^{m-1}.
\label{decroissance enk} \ee Thus, choosing $ C' _{crit}$ large enough (such that $ {C'
_{crit}}^2/2 > C $), it is easy to deduce that
the contribution of paths $\mathcal P$ for which $p \geq C' _{crit} \sqrt N $ is negligible in the large limit $N$.\\

We now assume that $p\leq C' _{crit} \sqrt N$ and come back to the estimation of (\ref{contri1}).
Note that $\mathcal P'$ is still a path of $T_{m,l}$ so that the same estimate holds for $\sum
z'(t).$ Now, it is easy to see that if $y'$ (resp. $y''$) is the Dyck path associated to $\mathcal
P'$ (resp. $\mathcal P''$) $\max y'(t)\leq \max y''(t)+p.$ This follows from the fact that there
are $p$ more down steps in $\mathcal P'$ than in $\mathcal P''$, and in $\mathcal P'$, after the
first instant at which $P'$ and $P''$ may differ at level $l+p-1$, one can not go below level
$l-1$. Note also that $m \leq 2(m-p)$ as $p<<m.$ We now turn to the estimation of
$\nu_N''(\mathcal P')$ as in Remark \ref{rem: nun}. Assumption ($A_2$) (and ($A_3$)) still holds.
Observe also that if the trajectory $x'$ of $\mathcal P'$ comes back from above to some given level $\nu_N/K$ times, then the trajectory $x'' \in \mathcal T_{m-p, l+2p}$ comes at least $\nu_N/2K$ times to some level (maybe different) without falling below. The probability of such an event still decreases as $4s_N^2 \exp{\{-C_1 \nu_N/2K\}}.$ In this way, using the same computations as in the analysis of $Z_5$ in Lemma \ref{lem: z5juste}, we shall be able to show that $\nu_N(\mathcal P')\leq s_N^{1/2-\epsilon}$ in typical paths. \\

Let $\mathbb{E}_{Y'',T''}$ denote the expectation with respect to
the uniform distribution on the trajectories of $\mathcal T _{2m-2p}^0$ of class $T''$. From Lemma \ref{Lem : max y}, one has that, given any constant $C$, $\exists \: \tilde C$, independent of $T''$, $p$ and $m$, such
that $\displaystyle{\mathbb{E}_{Y'',T''}\left (\exp{\{\frac{C\max y''(t)}{\sqrt m}\}}\right)\leq
\tilde C.}$ Using this and mimicking (\ref{formuletheta=sigma}), one obtains that the contribution
of paths with unmarked origin is
\begin{eqnarray} \label{contri2}
(\ref{contri1}) & \leq &\sum_{l\leq C_{crit}\sqrt N} \sum_{1\leq p\leq C' _{crit} \sqrt
N}\,T_{m-p,l+2p}\exp{\{\frac{3l^2/2+3ml}{N}+\frac{C'p}{\sqrt m}\}}\sigma^{2s_N}\sum_{q=0}^{10
s_N^2/N}\sum_{N_3, \ldots, N_{10}, M_1}\cr &&\times \frac{1}{T_{m-p,l+2p}}\sum_{(p',l_1,\ldots,
l_{p'}, m_o,m_1,\ldots ,
    m_{p'})}\:\mathbb{E}_{Y'',T''_{in}} \left (\exp{\{\frac{18m^{3/2}}{N} \frac{\max
    y''(t)}{\sqrt m}\}}\right)\: \:\prod_{0\leq i\leq p'/m_i>0} T_{m_i,0}\crcr
&& \frac{1}{q!}\left (\frac{Cs_N \nu''_N}{N}\right)^{q}\prod_{k= 3}^{10}\frac{1}{N_k!}\left
  (\frac{C(l+m)^k}{N^{k-1}}\right)^{N_k}\frac{1}{M_1!}\left
  (CN^{-2/9}\right)^{M_1}\exp{\{\frac{Cs_N^3}{N^2}\}}\crcr
  &&\leq \sum_{l\leq C_{crit}\sqrt N} \sum_{1\leq p\leq C' _{crit}
    \sqrt N }\,T_{m-p,l+2p} \sigma^{2s_N} \sum_{N_3, \ldots, N_{10}, M_1}\sum_{q=0}^{10 s_N^2/N}\exp{\{\frac{C'ml}{N}+\frac{C'p}{\sqrt m}\}} \crcr
  &&\frac{1}{q!}\left (\frac{Cs_N \nu''_N}{N}\right)^{q} \frac{1}{N_3!}\left (\frac{C' s_N^3}{N^2}\right)^{N_3}\prod_{k= 4}^{10}\frac{1}{N_k!}\left (\frac{C(l+m)^k}{N^{k-1}}\right)^{N_k}\frac{\left (N^{-1/9}\right)^{M_1}}{M_1!}\exp{\{\frac{Cs_N^3}{N^2}\}}
\end{eqnarray}
Call $Z_6$ the subsum corresponding to the case where $q=0$ and $M_2=M_1+\sum_{i\geq 4}N_i=0$. We
then show that there exist some constants $D_1, D_2>0$, independent of $N$, such that \be Z_6\leq
D_1(2\sigma)^{2s_N}\exp{\{D_2 s_N^3/N^2\}}.\label{z5}\ee To obtain (\ref{z5}), it is enough to prove
that there exists some constant $D>0$ such that \be \label{somme let k} \sum_{l=0}^{C_{crit} \sqrt
N}\sum_{p=1}^{C'_{crit}\sqrt N} T_{m-p,l+2p}\:\sigma^{2s_N}
\exp{\{\frac{C'lm}{N}\}}\exp{\{\frac{C'p}{\sqrt
    m}\}}\leq D 2^{2s_N}\sigma^{2s_N} .
\ee

Consider first the case where $p\geq k_ol$ where $k_o=\sup \frac{m^{3/2}}{N}$ is such that
$\frac{lm}{N}\leq \frac{p}{\sqrt m}$. Then, we can perform the summation over (even) $l$, yielding
that $\sum_{l=0}^{C_{crit} \sqrt N}T_{m-p,l+2p}\leq C_{2s_N}^{s_N+p}.$ We then deduce, as in any
case $m\geq s_N/4$, that $\exists D>0$ such that
$$\sum_{p=1}^{C'_{crit}\sqrt N}C_{2s_N}^{s_N+p}\exp{\{\frac{4C'p}{\sqrt s_N}\}} \leq
\exp{\{-\frac{4C's_N}{\sqrt s_N}\}}2^{2s_N}\left
  (1+\frac{e^{\frac{4C'}{\sqrt s_N}}-1}{2}\right)^{2s_N}\leq D 2^{2s_N}.$$
The case where $l\geq p/k_o$ is analyzed in a similar fashion. This gives (\ref{z5}).\\
From this result, we can readily deduce that typical paths $\mathcal P$ amongst those associated
to a path $\mathcal P'$ for which $\nu_N(\mathcal P')\leq s_N^{1/2-\epsilon},$ have no edges passed
more than twice and are without loop. Finally the contribution of paths $\mathcal P$ associated to
paths $\mathcal P'$ for which $\nu_N(\mathcal P')\geq s_N^{1/2-\epsilon}$ is analyzed as $Z_5$
(see Lemma \ref{lem: z5juste}). Their contribution is negligible with respect to
$(2\sigma)^{2s_N}.$ Combining the whole ensures that the limiting contribution of paths with last
step down depends only
on $\theta=\sigma. \, \square$ \\

To complete our proof, we shall now prove the following estimate.
 \bl \label{Lem : max y} Given a
constant $C$, there exists $\tilde C>0$, independent of $m$ and $T \in \mathbb{T}(m),$ such that
$$\mathbb{E} _{Y,T}\left (\exp \{ C \frac{\max y(t)}{\sqrt{m}}\}\right)\leq \tilde C.$$
\el  \brem A similar
estimate for general Dyck path $y$ of length $2m$ (i.e. without assuming a particular
decomposition in sub-Dyck paths) was used in \cite{Sos}(Lemma 6), but not proved. \erem
\paragraph{Proof of Lemma \ref{Lem : max y}:}We give here a
proof based on some geometrical considerations.\\
\textit{$1^{rst}$ case: No specified number of sub-Dyck paths.} We first show that there exist constants $C_o, \tilde C _o$ independent of
$m$, such that, under the uniform distribution on $\mathcal T_{2m}^0,$   \be \label{max DIck} \mathbb P(\max y(t)=k)\leq
\frac{\tilde C_o}{\sqrt m}\exp{\{-\frac{C_ok^2}{2m}\}}, \text{ if $k\geq 4\sqrt m C_o$}.\ee
 Consider first a general Dyck path
$y$ of length $2m.$ The probability that its maximum is $k$ is at most, if $\tilde T_{2n}^k$
denotes the number of paths with $2n$ steps ending at level $k$ without going below $0$ or above
$k$,
\begin{equation}
\mathbb P(\max y(t)=k )\leq \frac{1}{T_{2m}^0}\Big (\sum_{k/2\leq n \leq m-k/2} \left (
\frac{\tilde T_{2n}^k\: \tilde T_{2m-2n}^k+\tilde T_{2n+1}^k\: \tilde
T_{2m-2n-1}^k}{T_{2n}^0T_{2m-2n}^0}\right )T_{2n}^0T_{2m-2n}^0\Big ).
\end{equation}
This follows from the fact that such a path is the concatenation of a path $\tilde T_{2n(+1)}^k$
and one of $\tilde T_{2s-2n(-1)}^k.$ Actually, in the previous sum, $2n(+1)$ should be seen as the
first instant one reaches the level $k$. Now, we show that there exists constants $C_1,C'_1, C_o
>0$ such that for $k\geq 4C_o \sqrt m$,
\begin{eqnarray}
&&\tilde  T_{2n}^k\leq C_1 \frac{k^2 }{n}\exp{\Big \{
  -\frac{C_ok^2}{n}\Big\}}T_{2n}^0 +\frac{k}{\sqrt n}T_{2n}^0\exp{\Big
  \{-\frac{C_ok^2}{n}\Big \}}  \leq 2C'_1 \frac{k^2 }{n}\exp{\Big \{
  -\frac{C_ok^2}{n}\Big\}}T_{2n}^0, \cr
&&\tilde  T_{2m-2n}^k  \leq  2C'_1  \frac{k^2 }{m-n}\exp{\Big\{
  -\frac{C_ok^2}{(m-n)} \Big\}}T_{2(m-n)}^0 ,\label{bornesurtilde Tlk}
\end{eqnarray}
We only prove the inequality for $\tilde T_{2n}^k$. Recall first that $T_{2n}^k$ equals to the
number of positive paths $y$ of length $2n$ beginning at $0$ and ending at $k$. On the other hand,
by a simple application of the symmetry principle, it is easy to see that $ C_{2n}^{n+k/2+1}$
counts the number of paths that go from $0$ to $k$ touching -1. So, one can write that:
$$\tilde T_{2n}^k=T_{2n}^k - C_{2n}^{n+k/2+1}+\sharp \hat{\mathcal T} _k ^{-1,k+1}$$
where $\hat{\mathcal T} _k ^{-1,k+1}$ is the set of paths of length $2n$ that go from $0$ to $k$
touching both $-1$ and $k+1$. We shall estimate its cardinal. Note that either such a path goes to
$k+1$ after the first time it
goes to $-1$ or it goes first to $k+1$ then to $-1$ and joins $k$ without reaching $k+1$ afterwards.\\
First, by a simple symmetry principle, it is easy to see that paths which go to $k+1$ after the
first time it goes to $-1$ are in bijection with paths of length $2n$ beginning at
$0$ and ending at $-(k+4)$. So there are exactly $C_{2n}^{n+k/2+2}$ such paths.
Now, as $k\leq 2n\leq 2m-k$, by Lemma \ref{lem : majoexpTml}, one has \be
\frac{T_{2n}^k-C_{2n}^{n+k/2+1}+C_{2n}^{n+k/2+2}}{T_{2n}^0}=\frac{T_{2n}^k-T_{2n}^{k+2}}{T_{2n}^0}
\leq   C_2 \frac{k^2 }{n}\exp{\Big \{ -\frac{C_ok^2}{n}\Big\}}.\label{est nombre dyck1} \ee
where $C_2$ and $C_o$ are two constants independent of $n$. Note that $C_o$ can be chosen equal to $96^{-1}$. \\
The numbering of the paths which reach $k+1$ before $-1$ and joins $k$ without reaching $k+1$
afterwards is quite more subtle and we only obtain an upper bound.
The trajectory of such a path can be described as follows. Call $2t+1$ the last instant where the trajectory is at level $k+1.$ Call then $2t'+1$ the last instant where the trajectory is at level $-1$.\\
Assume first that $2t+2\leq 2n(2/3).$ In Between $0$ and $2t+1$ the trajectory goes from $0$ to
$k+1$ without touching $-1$. Then the trajectory in between $2t+2$ and $2n$ goes from $k$ to $k$
without touching $k+1$ (and reaches $-1$, but we will forget this constraint). Then, using that
for any $n$ and $k\leq 2n,$
 $\displaystyle{T_{2n}^k\leq k\exp{\{-\frac{C_ok^2}{n}\}}T_{2n}^0}$ (c.f. Lemma \ref{lem : majoexpTml}),
 the number of such paths is at most
\begin{eqnarray}
 &&\sum_{2t\leq 2n(2/3)} T_{2t+1}^{k+1}T_{2n-2t-2}^0\cr
 &&\leq
2(k+3)\exp{\{-\frac{C_ok^2}{n}\}}\sum_{2\leq 2t\leq
2n(2/3)}\exp{\{-C_ok^2(\frac{1}{t}-\frac{1}{n})\}}T_{2t}^0T_{2n-2t-2}^0\cr &&\leq
C_3(k+3)\exp{\{-\frac{C_ok^2}{n}\}}\frac{n^{1+3/2}}{n^3}T_{2n}^0\int_0^{2/3}\frac{1}{u^{3/2}(1-u)^{3/2}}\exp{\{-\frac{C_ok^2}{n}(1/u-1)\}}
\, du \label{intaboner}\\ &&\leq C_4\frac{k}{\sqrt {n}}T_{2n}^0\exp{\{-\frac{C_ok^2}{n}\}}
.\label{est nombre dyck2}
\end{eqnarray}
To derive the last line, we have used the fact that $C_ok^2/2n \geq 2C_o^3$ to bound the integral in (\ref{intaboner}).\\
 If now $2t\geq 2n(2/3),$
then $2t'>2n(2/3).$ And the path can be described as follows. Between $0$ and $2t'$, the path goes
from $0$ to $0.$ Then, the path goes from $-1$ to $k$ without touching $k+1$ in $2n-(2t'+1)$ steps
with $2n-(2t'+1)\leq 2n(2/3).$ Thus the number of such paths can be majorized as above.
Formulas (\ref{est nombre dyck1}) and (\ref{est nombre dyck2}) finally imply formula (\ref{bornesurtilde Tlk}).\\

Set now $n=um$. There exists $C_5$ independent of $n$ and $m$, such that \be \label{labonne majo}
\frac{\tilde T_{2n(+1)}^k \:\tilde T_{2m-2n(-1)}^k}{T_{2n}^0\:T_{2m-2n}^0}\leq C_5 \exp {\Big \{
-\dfrac{C_ok^2}{2m}\left
  (\dfrac{1}{u}+\dfrac{1}{1-u}\right )\Big\}}.\ee

Thus by (\ref{labonne majo}), one has
\begin{eqnarray}
&& \mathbb P(\max y(t)=k )\leq C_5\sum_{n=um}\frac{\exp{\Big \{
-\frac{C_ok^2}{2m}(1/u+1/(1-u)\Big\}}}{\exp{\{-\frac{2C_ok^2}{m}\}}}
\frac{T_{2n}^0T_{2m-2n}^0}{T_{2m}^0}\exp{\{-\frac{2C_ok^2}{m}\}}.\label{P max}
\end{eqnarray}

Such sum will be divided in three subsums: according to $m/10 \leq n \leq 9m/10$, then to $n \leq
m/10$ and finally for $n \geq 9m/10$. For $m/10 \leq n \leq 9m/10$, we can use Stirling's formula
to obtain that this subsum can be majorized by a term similar to the announced bound, since \be
\sum_{m/10 \leq n \leq 9m/10} T_{2n}^0T_{2m-2n}^0 \leq \frac{C}{\sqrt m}T_{2m}^0.\label{l=m}\ee
When $n\leq m/10$, we use the fact that $\displaystyle{\frac{\exp{\Big \{
-\frac{C_ok^2}{2}(\frac{1}{n}+\frac{1}{m-n})\Big\}}}{\exp{\{-\frac{2C_ok^2}{m}\}}}\leq \exp{\{-C_o
\frac{k^2}{4n}\}}.}$ Thus \be \sum_{ n\leq m/10}\frac{\exp{\Big \{
-\frac{C_ok^2}{2}(\frac{1}{n}+\frac{1}{m-n})\Big\}}}{\exp{\{-\frac{2C_ok^2}{m}\}}}T_{2n}^0T_{2m-2n}^0
\leq \sum_{n\leq m/10}\exp{\Big \{ -\frac{C_ok^2}{4n}\Big \}}\frac{C_6}{n^{3/2}} T_{2m}^0\leq \frac{C_7}{\sqrt m}T_{2m}^0,\label{l <1/4 m} \ee for some constant $C_6, C_7$ independent of $m$,
by a straightforward comparison with an integral. Here the constant $C_7$ does not depend on $m,$
as $k\geq 4C_o\sqrt m.$ We can obtain by symmetry a similar bound for the sum over $n\geq \frac{9m}{10}$ (as $m-n \leq \frac{m}{10}$). Combining then (\ref{l <1/4 m}), (\ref{l=m}) and (\ref{P max})
leads
to (\ref{max DIck}).\\

\textit{$2^{nd}$ case: A specified number of sub-Dyck paths.} We can now finish the proof of Lemma
\ref{Lem : max y}. Here it is enough to prove that $\exists\: \tilde C'_o, C_o$,
independent of $m$ and $T\in \mathbb{T}(m)$, such that, under the uniform distribution on trajectories of $\mathcal T_{2m}^0$ and class $T=(p'+1, m_o, \ldots,m_{p'})$, then \be
\label{estmax y}\mathbb P(\max y(t)=k)\leq \frac{\tilde C'_o}{\sqrt
  m}\exp{\{-C_o\frac{k^2}{2m}\}}, \text{ if $k\geq 4C_o\sqrt m$. }\ee Here $C_o$ is the same constant as in (\ref{max DIck}).
Set $\alpha_i=m_i/m$ so that $\sum_{i=0}^{p'} \alpha_i=1.$ Obviously, if $k\geq
C_o  \sqrt m$ then $k\geq C_o\sqrt {m_i}$. Then, by the above computations, one has that
\begin{eqnarray}
&&\mathbb P(\max y(t)=k)=\mathbb P(\exists i\leq p', \max y_i(t)=k )\leq \sum_{i=0}^{p'}
\frac{\tilde C_o}{\sqrt m}e^{\Big \{ -\frac{C_ok^2}{2m}\Big \}} \frac{1}{\sqrt
\alpha_i}e^{\Big \{  -\frac{C_ok^2}{2m}(\frac{1}{\alpha_i^2}-1)\Big\}}\cr &&\leq
\frac{C_2}{\sqrt m} \exp{\Big \{ -\frac{C_ok^2}{2m}\Big \}}\sum_{i=0}^{p'} \exp{\Big \{
-\frac{C_ok^2}{4m}(\frac{1}{\alpha_i^2}-1)\Big\}},
\end{eqnarray}
where in the last line we have used that $1/\sqrt {\alpha_i}\leq c\exp{\Big \{
\frac{C_ok^2}{4m}(\frac{1}{\alpha_i^2}-1)\Big\}}$ as $k\geq 4C_o\sqrt m$ and $\alpha_i\leq 1.$
Note that the constant $C_2$ does not depend on $p'$. Then, (\ref{estmax y}) holds since one has
$$\sum_{i=1}^{p'}\exp{\Big \{ -\frac{C_ok^2}{4m}(\frac{1}{\alpha_i^2}-1)\Big\}}\leq
\sum_{q\geq 1}(q+1) \exp{\{4C_o^3(1-q^2)\}} \leq A,$$ by using the fact that the number of $\alpha_i$ in
any interval $[1/(q+1), 1/q]$ with $q\geq 1$ is not greater than $q+1$ (since
$\sum_{i=0}^{p'}\alpha_i=1$). Thus $A$ is a constant independent of $p'$. Note that this estimate holds
for any value of the $\alpha_i$ also. This finishes the proof of Lemma \ref{Lem : max y}. $\square$

\brem The investigation of higher moments is a mimicking of the arguments of Section 6 and
\cite{Sos} (p. 42). This is not detailed further.\erem
 \brem \label{casetheta<sigma} In the case where
$0 < \theta < \sigma$, it suffices to observe that $\rho _{\theta}< 2\sigma$. Thus, all the
results of Sections 4 to 7 show that contribution of paths having at least one unreturned edge
($l>0$) is negligible in the expectation and higher moments at any scale $1<< s_N \leq O(N^{2/3})$
and that the main contribution comes from even paths ($l=0$). As a result, if $0<\theta<\sigma$,
the decentring matrices $A_N$ do not affect the limiting behaviour of the largest eigenvalue of
$M_N$. This behavior is then the same as that of the Gaussian Ensemble of the same symmetry (and
for which $\theta=0$). In particular, the limiting distribution of the largest eigenvalue is given
by the classical GUE or GOE Tracy-Widom distribution. The same conclusion holds for the joint
distribution of the $k$ first largest eigenvalues, for any fixed integer $k \geq 1$. This
completes the proof of Theorems $\ref{UnikPe}$ and Theorem \ref{theo: Unireal}. \erem


\addcontentsline{toc}{chapter}{Bibliographie} \markboth{BIBLIOGRAPHIE}{BIBLIOGRAPHIE}
\begin{thebibliography}{50}
\bibitem {Bai} \, Z. Bai,
\newblock Methodologies in spectral analysis of large-dimensional random matrices, a review,
\newblock{\em Statist. Sinica }\textbf{9}: 611--677 (1999).

\bibitem {BBP} \, J. Baik, G. Ben Arous and S.  P\'ech\'e,
\newblock Phase
    transition of the largest eigenvalue for non-null complex sample
    covariance matrices,
    \newblock{\em Ann. Probab. }\textbf{33} no 5:  1643--1697 (2005).

\bibitem {BS} \,J. Baik and J. Silverstein, forthcoming paper (2005).


\bibitem{Fe}\, D. F\'eral,
\newblock On the extremal eigenvalues of large deformed Wigner matrices,
\newblock to appear in ESAIM,  available at http://www.lsp.ups-tlse.fr/Fp/Feral/DWM.ps (2005).
\bibitem {Fu-Ko} \, Z. Furedi and J. Komlos,
\newblock The
    eigenvalues of random symmetric matrices,
   \newblock{\em  Combinatorica}
    \textbf{1}: 233--241 (1981).
\bibitem{Geman} S. Geman, \newblock A limit theorem for the norm of random matrices
\newblock{\em Ann. Prob.} \textbf {8}:
252--261 (1980).

\bibitem {GW} \, P. W. Glynn and W. Whitt,
\newblock Departures from many
    queues in series,
    \newblock{\em Ann. of Applied Prob.} \textbf{4}: 546--572 (1991).

\bibitem {Debashis} \, D. Paul, \newblock Asymptotics of the leading sample eigenvalues for a spiked covariance model , Technical Report Stanford University
\newblock Available at http://www-stat.stanford.edu/~debashis/ (2004).

\bibitem {Pe} \, S. P\'ech\'e,
\newblock The largest eigenvalues of small rank perturbations of
    Hermitian random matrices,
 \newblock{\em    Prob. Theo. Rel. Fields} \textbf{134} no 1: 127--174 (2006).

\bibitem {Si-So1} \, Y. Sinai and A. Soshnikov,
\newblock Central limit
    theorem for traces of large random symmetric matrices with independent
    matrix elements,
 \newblock {\em   Bol. Soc. Brasil. Mat. (N.S.)} \textbf{29}: 1--24 (1998).

\bibitem {Si-So2} \, Y. Sinai and A. Soshnikov,
\newblock A refinement of
    Wigner's semicircle law in a neighborhood of the spectrum edge for random
    symmetric matrices,
    \newblock {\em Funct. Anal. Appl.} \textbf{32}: 114--131 (1998).

\bibitem {Sos} \, A. Soshnikov,
\newblock Universality at the edge of the
    spectrum in Wigner random matrices,
     \newblock {\em Comm. Math. Phys.} \textbf{207}: 697--733 (1999).

\bibitem{TracyWidom}
C.A.~Tracy and H.~Widom, \newblock Level spacing distributions and the {A}iry kernel.
\newblock {\em Comm. Math. Phys.}, \textbf{159}: 33--72 (1994).

\bibitem {TracyWidomFred}
C.A. Tracy and H. Widom, \newblock Fredholm determinants, differential equations and matrix models.
\newblock {\em Commun.Math.Phys}, \textbf{163}: 33--72 (1994).
\end{thebibliography}
\end{document}